\documentclass[a4paper]{article}

\usepackage{epsf}
\usepackage{amsmath, amsthm}
\usepackage{amsfonts}
\usepackage{amssymb}
\usepackage{graphicx,epsfig,float}
\include{psfig}
\usepackage{lscape}
\usepackage[usenames]{color}
\usepackage{pifont}
\usepackage{yhmath}
\usepackage[latin1]{inputenc}
\usepackage{amssymb,amsmath,amsfonts}
\usepackage[english]{babel}
\usepackage{tabularx}
\usepackage{amsthm}
\usepackage{dsfont}
\usepackage{epsfig}
\usepackage{mathrsfs}
\usepackage{fullpage}

\newcommand{\cl}{\centerline}

\newcommand{\bit} {\begin{itemize} }
\newcommand{\eit} {\end{itemize} }

\newcommand\ho{H\"{o}lder }

\newcommand{\§}{§ }

\newcommand{\bN}{\mathbb{N}}

\newcommand{\one}{\ifmmode {\sf 1}\hspace{-.26em}{\sf
l}\hspace{-.35em}{\sf \_} \else ${\sf 1}\hspace{-.26em}{\sf
l}\hspace{-.35em}{\sf \_}$ \fi}

\def\ben{\begin{enumerate}}
\def\een{\end{enumerate}}

\def\iti{\item[(i)]}

\def\itii{\item[(ii)]}

\def\itiii{\item[(iii)]}

\def\cB{{\cal B}}
\def\cC{C}

\def\cE{{\cal E}}
\def\cF{{\cal F}}
\def\cG{{\cal G}}
\def\cH{{\cal H}}
\def\cI{{\cal I}}
\def\cJ{{\cal J}}
\def\cK{{\cal K}}
\def\cL{{\cal L}}

\def\cN{{\cal N}}

\def\cS{{\cal S}}

\newcommand{\sS}{\mathscr{S}}

\def\bC{{\mathbb C}}
\def\bD{{\mathbb D}}
\def\bE{{\mathbb E}}

\def\bN{{\mathbb N}}

\def\bQ{{\mathbb Q}}
\def\bR{{\mathbb R}}
\def\bZ{{\mathbb Z}}

\def\i1{\mathds{1}}  

\def\si{\text{sign}}  

\def\2sp{\vspace{2cm}}
\def\3sp{\vspace{3cm}}
\def\[ent{[\hskip -1.5pt [}
\def\]ent{]\hskip -1.5pt ]}

\def\2sp{\vspace{2cm}}
\def\3sp{\vspace{3cm}}

\def\h1{\hspace{0.1cm}}

\newcommand{\reals}{\ifmmode {\sf I}\hspace{-.15em}{\sf R} \else ${\sf
    I}\hspace{-.15em}{\sf R}$ \fi}

\renewcommand{\leq}{\leqslant}
\renewcommand{\geq}{\geqslant}

\def\gg{{\textquotedblleft}}

\def\h1{{\hspace{0.1cm}}}

\newtheorem{theo}{Theorem}[section]
\newtheorem{theodef}{Theorem-Definition}[section]
\newtheorem{defi}{Definition}[section]
\newtheorem{prop}[theo]{Proposition}

\newtheorem{propers}{Properties}[section]
\newtheorem{lem}[theo]{Lemma}
\newtheorem{ex}[theo]{Example}
\newtheorem{cor}[theo]{Corollary}
\newtheorem{exa}[theo]{Example}
\newtheorem{rem}[theo]{Remark}

\newenvironment{defappli}[4]{\begin{array}{cccl} %
#1 \, : & #2 & \rightarrow & #3 \\ & #4 & \mapsto &}%
{\end{array}}

\newenvironment{g}[2]{\begin{array}{cl} %
<\hspace{-0.2cm}<\hspace{-0.1cm}#1,\ \hspace{-0.45cm} \hspace{-1cm}&#2\hspace{-0.1cm}>\hspace{-0.2cm}>}%
{\end{array}}

\newenvironment{defappliab}[6]{\begin{array}{cccccl} %
#1 \, : & #2    & \stackrel{\unboldmath{M_{h}}}{\rightarrow}   & #3  & \stackrel{\unboldmath{\zeta}}{\rightarrow}     & #4 \\
        & #5    &  \mapsto      &  #6   & \mapsto        &  }%
{\end{array}}

\def\qed {{
\parfillskip=0pt 
\widowpenalty=10000 
\displaywidowpenalty=10000 
\finalhyphendemerits=0 
\leavevmode 
\unskip 
\nobreak 
\hfil 
\penalty50 
\hskip.2em 
\null 
\hfill 
$\square$
\par}} 

\newenvironment{prz}{{\bfseries \textup{Proof.}}}{}
\newenvironment{pr}{\begin{proof}[\bfseries \textup{Proof.}]}{\end{proof}}
  \newenvironment{praa}{{\bfseries \textup{Proof }}}{}

\newtheorem{proe}{Proof}

\def\mathtitre#1{
\font\tenrm=cmr10   scaled \magstep#1
\font\sevenrm=cmr7  scaled \magstep#1
\font\fiverm=cmr5   scaled \magstep#1
\font\teni=cmmi10   scaled \magstep#1
\font\seveni=cmmi7  scaled \magstep#1
\font\fivei=cmmi5   scaled \magstep#1
\font\tensy=cmsy10  scaled \magstep#1
\font\sevensy=cmsy7 scaled \magstep#1
\font\fivesy=cmsy5  scaled \magstep#1
\font\tenex=cmex10  scaled \magstep#1
\textfont0=\tenrm \scriptfont0=\sevenrm \scriptscriptfont0=\fiverm
\textfont1=\teni  \scriptfont1=\seveni  \scriptscriptfont1=\fivei
\textfont2=\tensy \scriptfont2=\sevensy \scriptscriptfont2=\fivesy
\textfont3=\tenex \scriptfont3=\tenex \scriptscriptfont3=\tenex
}

\makeatletter
\renewcommand\theequation{\thesection.\arabic{equation}}
\@addtoreset{equation}{section}
\makeatother

\def\keywordname{{\bf Keywords:}}
\newcommand{\keywords}[1]{\par\addvspace\baselineskip\noindent\keywordname\enspace\ignorespaces#1}

\renewcommand{\thefootnote}{\Alph{footnote}}

\newcommand{\footnoteremember}[2]{ 
 \footnote{#2}
 \newcounter{#1}
 \setcounter{#1}{\value{footnote}}
 \addtocounter{#1}{0} 
}
\newcommand{\footnoterecall}[1]{
 \footnotemark[\value{#1}]
}

\begin{document}
\setlength{\parindent}{0em}

\makeatletter
\renewcommand\theequation{\thesection.\arabic{equation}}
\@addtoreset{equation}{section}
\makeatother

\title{White Noise-based Stochastic Calculus\\ with respect to multifractional Brownian motion}
\author{Joachim \textsc{Lebovits}\footnoteremember{myfootnote}{I.N.R.I.A. Saclay, Regularity team: jolebovits@gmail.com \h1 and \h1 jacques.levy-vehel@inria.fr.}\setcounter{footnote}{\value{footnote}}\addtocounter{footnote}{2}\hspace{-0.25cm}\footnote{Laboratoire de Probabilités et Modèles Aléatoires, C.N.R.S (UMR 7599), Université Pierre et Marie Curie (Paris VI), case 188, 4, pl. Jussieu, F-75252 Paris Cedex 5, France. }
\\
\and
Jacques \textsc{Lévy-Véhel}\footnoterecall{myfootnote}
}
\maketitle
\begin{abstract}
\par \noindent Stochastic calculus with respect to fractional Brownian motion (fBm) has attracted a lot of
interest in recent years, motivated in particular by applications in finance and Internet traffic modeling. Multifractional
Brownian motion (mBm) is a Gaussian extension of fBm that allows to control the pointwise regularity of the paths of the
process and to decouple it from its long range dependence properties. This generalization is obtained by replacing the constant
Hurst parameter $H$ of fBm by a function $h(t)$. Multifractional
Brownian motion has proved useful in many
applications, including the ones just mentioned.

In this work we extend to mBm the construction of a stochastic integral with respect to fBm. This stochastic integral is
based on white noise theory, as originally proposed in \cite{ell}, \cite{bosw},  \cite{ben1} and in  \cite{ben2}. In that view, a multifractional white noise is defined, which allows to integrate with respect to mBm a large class of stochastic processes using Wick products.
Itô formulas (both for tempered distributions and for functions with sub-exponential growth) are given,
along with a Tanaka Formula. The cases of two specific functions $h$ which give notable Itô formulas are presented.
\end{abstract}

\pagebreak


\tableofcontents


\section{Background and Motivations}
\label{rappels sur le mBm}

Fractional Brownian motion (fBm) $\cite{Kol,MVN}$ is a centered Gaussian process with features that make it a useful model in various applications such as financial and teletraffic modeling, image analysis and synthesis, geophysics and more. These features include self-similarity, long range dependence and the ability to match any prescribed constant local regularity. Fractional Brownian motion depends on a parameter, usually denoted $H$ and called the Hurst exponent, that belongs to $(0,1)$. Its covariance function $R_{H}$ reads:
\begin{equation}
\label{covfBm}
 R_{H}(t,s) :=  \frac{\gamma_{H}}{2} ( {|t|}^{2{H}} + {|s|}^{2{H}} - {|t-s|}^{2{H}}),
\end{equation}

where $\gamma_{H}$ is a positive constant. A normalized fBm is one for which $\gamma_{H}=1$. Obviously, when $H=\frac 1 2$, fBm
reduces to standard Brownian motion.

A useful representation of fBm $B^{(H)}$ of exponent $H$ is the so-called harmonizable one:
\vspace{-0.2cm}
\begin{equation}
\label{fBm}
 B^{(H)} (t) =  \displaystyle{ \frac{1}{c_H}\int_{\bR} \frac{e^{itu} - 1}{{|u|}^{H+1/2 }} \widetilde{W}(du)},
\end{equation}

where $c_x :=  {\bigg(\frac{2\cos(\pi x) \Gamma(2-2x)}{x(1-2x)} \bigg)}^{\frac{1}{2}} = { \bigg(  \frac{2\pi}{ \Gamma(2x+1)\sin(\pi x)}  \bigg )  }^{\frac{1}{2}}$ for $x$ in $(0,1)$ and $\widetilde{W}$ denotes the complex-valued Gaussian measure which can be associated in a unique way to $W$, an independently scattered standard Gaussian measure on $\bR$ (see \cite{StoTaq} p.203-204 and \cite{TaqSam} p.325-326 for more information on the meaning of $\int_{\bR} f(u) \h1 \widetilde{W}(du)$ for a complex-valued function $f$). From \eqref{covfBm} and Gaussianity, it is not hard to prove that fBm is $H$-self-similar.

The fact that most of the properties of fBm are governed by
the single number $H$ restricts its application in some situations. Let us give two examples. The long term
correlations of the increments of fBm decay as $k^{(2H-2)}$, where $k$ is the lag, resulting in long range dependence
when $H > 1/2$ and anti-persistent behavior when $H < 1/2$.  Also, almost surely, for each $t$, its pointwise H\"older
exponent is equal to $H$. Since $H$ rules both ends of the Fourier spectrum,
{\it i.e}. the high frequencies related to the \ho regularity and the low frequencies
related to the long term dependence structure, it is not
possible to have at the same time {\it e.g.} a very irregular local behavior (implying $H$ close to 0)
and long range dependence (implying $H >
1/2$). As a consequence, fBm is not adapted to model phenomena which display both these features,
such as Internet traffic or certain highly textured
images with strong global organization. Another
example is in the field of image synthesis: fBm has frequently been used for generating artificial mountains.
Such a modeling assumes that the regularity of the mountain is everywhere the same. This is not realistic, since
it does not take into account erosion or other meteorological phenomena which smooth some parts of mountains
more than others.

Multifractional Brownian motion (mBm) \cite{PL,BJR} was introduced to overcome these
limitations. The basic idea is to replace in \eqref{fBm} the real $H$ by a function $h(t)$. More precisely, we will use the
following definition of mBm:
\begin{defi} [Multifractional Brownian motion] 
\label{mBm}
Let $h : \bR \rightarrow (0,1)$ be a continuous function and $\alpha : (0,1) \rightarrow \bR$ be a $\cC^1$  function.
A multifractional Brownian motion with functional parameters $h$ and $\alpha$ is defined as:
\vspace{0.2cm}
\begin{equation}
\label{dmBmm}
B^{(h,\alpha)} (t) = \alpha(h(t))  \displaystyle{\int_{\bR} \frac{e^{itu} - 1}{{|u|}^{h(t)+1/2 }} \h1 \widetilde{W}(du)}.
\end{equation}
\end{defi}

Its covariance function reads \cite{ACLVL}:
\vspace{-0.2cm}
\begin{equation}
\label{covmbm}
R_{(h,\alpha)}(t,s) = \alpha(h(t)) \hspace{0.1cm}\alpha(h(s)) \hspace{0.1cm} c^2_{h_{t,s}} \hspace{0.1cm} \bigg[\hspace{0.1cm} \frac{1}{2} \big( {|t|}^{2h_{t,s}} + {|s|}^{2h_{t,s}}  - {|t-s|}^{2h_{t,s}}     \big)\bigg],
\end{equation}
\vspace{-0.25cm}

where $h_{t,s} := \frac{h(t)+h(s)}{2}$ and $c_x$ has been defined in \eqref{fBm}.

It is easy to check that mBm is a zero mean Gaussian process, the increments of which are in general
neither independent nor stationary.

For $T$ in $\bR^{*}_+$,  we will again call $(h,\alpha)$-multifractional Brownian motion on $[0,T]$ the centered Gaussian process  whose covariance function is equal to $R_{(h,\alpha)}$ on $[0,T]\times[0,T]$.

When $\alpha = \alpha_c: x \mapsto  \frac{1}{c_x}$, we get that:

\vspace{-0.3cm}
\begin{equation}
\label{azazazaz}
R_{(h,\alpha_c)}(t,s) =   \hspace{0.1cm}  \frac{ c^2_{h_{t,s}}}{c_{(h(t))}c_{(h(s))}}  \hspace{0.1cm} \bigg[\hspace{0.1cm} \frac{1}{2} \big( {|t|}^{2h_{t,s}} + {|s|}^{2h_{t,s}}  - {|t-s|}^{2h_{t,s}}     \big)\bigg].
\end{equation}

As a consequence, if $h$ is the constant function equal to the real $H$, then $B^{(H,\alpha_c)}$ is a normalized fBm. For this reason, we will call $B^{(h,\alpha_c)}$ a normalized mBm. Since in the sequel we will consider only normalized mBm, we
simplify the notation and write from now on $B^{(h)}$ for $B^{(h,\alpha_c)}$ and $R_h$ for $R_{(h,\alpha_c)}$.

One can show \cite{Her,HLV} that the pointwise Hölder exponent at any point $t$ of $B^{(h)}$ is almost surely equal to $h(t)\wedge \beta_{h}(t)$, where $\beta_{h}(t)$ is the pointwise Hölder exponent of $h$ at $t$.
In addition, the increments of mBm display long range dependence for all non-constant $h(t)$ (the notion of long range dependence must be re-defined carefully for non-stationary increments, see
\cite{ACLVL}). Finally, at least when $h$ is $C^1$, for all $u \in \bR$, mBm locally "looks like" fBm with exponent $h(u)$ in the neighbourhood
of $u$ in the following sense \cite{PL}:
\begin{equation}
\lim_{r \to 0_+}\frac{B^{(h)}(u+rt) -B^{(h)}(u)}{r^{h(u)}} =  B^{(h(u))}(t),
\label{locform1}
\end{equation}
where the  convergence holds in law. These properties show that mBm is a more versatile model that fBm: in particular, it is able to mimic in a more faithful way local properties
of financial records, Internet traffic and natural landscapes \cite{bianchi,LLHF,montagnes} by matching their local regularity. This is
important {\it e.g.} for purposes of detection or real-time control. The price to pay is of course that one has to deal with the added complexity brought by having a functional parameter instead of a single number.

Because of applications, in particular in finance and telecommunications, it has been an important objective in recent years
to define a stochastic calculus with respect to fBm. This was not a trivial matter, as fBm is not a semi-martingale for $H \neq \frac 1 2$. Several approaches have been proposed, based mainly on Malliavin calculus \cite{DecUs,nualart,MNIRV}, pathwise
approaches and rough paths (\cite{Z,COU,FV} and references therein), and white noise theory \cite{ell,bosw,ben2}. Since mBm seems to be a more flexible, albeit more complex, model than fBm, it seems desirable to extend the stochastic calculus defined for fBm to it. This is the aim of the current work. In that view, we will use a white noise approach, as it offers several advantages in our frame. The main
task is to define a {\it multifractional white noise} as a Hida stochastic distribution, which generalizes the fractional white noise of, {\it e.g.},
\cite{ell,bosw}. For that purpose, we use the properties of the Gaussian field
${(B^{(H)}(t))}_{(t,H)\in\bR \times (0,1)}$. In particular, it is
a crucial fact for us that the function $H \mapsto (B^{(H)}(t))$ is almost surely $C^{\infty}$. This entails that multifractional white noise behaves essentially as fractional white noise, plus a smooth term.
We obtain an Ito formula that reads:

\begin{multline*}
f(T,B^{(h)}(T)) = f(0,0) + \int^T_0 \h1 \frac{\partial f}{\partial t}(t,B^{(h)}(t)) \h1 dt + \int^T_0 \h1 \frac{\partial f}{\partial x}(t,B^{(h)}(t)) \h1 dB^{(h)}(t) \\
+ \frac{1}{2} \h1 \int^T_0 \h1 \left( \frac{d}{dt}[R_{h}(t,t)] \right)\h1 \frac{\partial^2 f}{\partial x^2}(t,B^{(h)}(t)) \h1 dt,
\end{multline*}
where the meaning of the different terms will be explained below.

The remaining of this paper is organized as follows. In section \ref{White}, we recall basic facts about white noise theory. We study a family of operators, noted ${(M_H)}_{H \in (0,1)}$, which are
instrumental for constructing the stochastic integral with respect to mBm in section \ref{MH}. Section \ref{WienerInt} defines the Wiener integral with respect to mBm.
We build up a stochastic integral with respect to mBm in section \ref{sjjjj}. Various instances of Ito formula are proved in section \ref{Ito}. Finally, section \ref{Tanaka} provides a Tanaka formula, along with
the study of two particular $h$ functions that give notable results. Readers familiar with white noise theory may skip the next section.

\section{White noise theory}
\label{White}

We recall in this section the standard set-up for classical white-noise theory. We refer {\it e.g.} to \cite{Kuo2, HKPS} for more details.

\subsection{White noise measure}

Let  $\sS(\bR):=\{ f \in C^{\infty}(\bR)  \h1: \h1 \forall  (p,q) \in \bN^2, \h1 \lim\limits_{|x| \to +\infty}  \h1| {x}^p \h1 f^{(q)}(x)| = 0     \}$ be the Schwartz space.  A family of functions ${(f_n)}_{n \in \bN}$ of ${(\sS(\bR))}^{\bN}$ is said to converge to $0$ as $n$ tends to $+\infty$ if for all $(p,q)$ in ${\bN}^2$ we have $\lim\limits_{n \to +\infty}  \h1 \underset{x \in \bR }{\sup}{\h1| x^p \h1 f^{(q)}_{n}(x)|} = 0$. The topology hence given on $\sS(\bR)$ is called the usual topology. Let $\sS'(\bR)$ denote the space of tempered distributions, which is the dual space of $\sS(\bR)$. The Fourier transform of a function $f$ which belongs to $L^1(\bR)\cup L^2(\bR)$ will be denoted $\widehat{f}$ or $\cF(f)$:

 \begin{equation}
\label{eigtth}
 \cF(f)(\xi) := \widehat{f}(\xi) := \displaystyle{\int_{\bR} e^{-ix\xi}  f(x) dx, \hspace{0.25cm} \xi \in \bR. }
\end{equation}

Define the probability space as $ \Omega:= {\sS}^{'}(\bR)$ and let $\cF := \cB({\sS}^{'}(\bR))$ be the $\sigma$-algebra of Borel sets. The Bochner-Minlos theorem ensures that there exists a unique probability measure on $\Omega$, denoted $\mu$, such that:

\begin{equation}\label{first}
\int_{{\sS}^{'}(\bR)} e^{i<\omega,f>} \mu(d\omega) = e^{-\frac{1}{2} {||f||}^2_{L^2(\bR)}}, \hspace{0.5cm}  \forall f \in \sS(\bR),
\end{equation}

where $<\omega,f>$ is by definition   $\omega(f)$, \textit{i.e} the action of the distribution $\omega$ on the function $f$. For $f$ in ${\sS}(\bR)$ the map, noted $<.,f>$, from $\Omega$ to $\bR$ defined by $<.,f>(\omega) = <\omega,f>$ is thus a centered Gaussian random variable with variance equal to ${||f||}^2_{L^2(\bR)}$ under the probability measure $\mu$, which is called the {\it white-noise probability measure}. In other words,
$$\bE[<.,f>]= 0,  \hspace{2cm} \bE[{<.,f>}^2]= {||f||}^2_{L^2(\bR)},$$
for all $f$ in $\sS(\bR)$. Besides, for a measurable function $F$, from  ${\sS}^{'}(\bR)$ to $\bR$, the expectation of $F$ with respect to $\mu$ is defined, when it exists, by $E[F] := E_{\mu}[F] := {\int_{\Omega}F(\omega) \mu(d\omega)}$. Equality (\ref{first}) entails  that the map $\zeta$ defined on $\sS(\bR)$ by


%


%
%


\begin{equation}\label{fourth}
\begin{defappli}{\zeta}{(\sS(\bR),{< , >}_{L^2(\bR)})}{{(L^2(\Omega,\cF,\mu), {< , >}_{(L^2(\Omega,\cF,\mu)})}}{f}
\zeta(f):= {<., f >}
\end{defappli}
\end{equation}

is an isometry. Thus, it extends to $L^2{(\bR)}$ and we still note $\zeta$ this extension. For an arbitrary $f$ in $L^2(\bR)$, we then have $<., f > := \lim\limits_{n \to +\infty} <., f_n > $ where the convergence takes place in $L^2(\Omega,\cF,\mu)$ and where ${(f_n)}_{n \in \bN}$ is a sequence of functions which  belongs to $\sS(\bR)$ and  converges to $f$ in $L^2{(\bR)}$. In particular, define for all $t$ in $\bR$, the indicator function $\i1_{[0,t]}$ by

%
%



\begin{equation*}
 \i1_{[0,t]}(s) =
\begin{cases}
1  & {\text if} \hspace{0.25cm} 0 \leq s \leq t,\\
-1 &  {\text if}\hspace{0.25cm}  t \leq s \leq 0 \hspace{0.2cm} \text{except if}\hspace{0.1cm} t = s = 0 \\
0  &  {\text{otherwise}},
\end{cases}
\end{equation*}


Then the process ${({\widetilde{B}}_t)}_{t \in \bR}$, defined for $t \in \bR$, on $\Omega$ by:


\begin{equation}
\label{fifthufhceufdhvdfhv}
{\widetilde{B}}_t (\omega):= {\widetilde{B}}(t,\omega):= {<\omega, \i1_{[0,t]}>}
\end{equation}


 is a standard Brownian motion with respect to $\mu$. It then admits a continuous version which will be denoted $B$.
 Thanks to  (\ref{fifthufhceufdhvdfhv}) we see that, for all  functions  $f$ in $L^2(\bR)$,


\begin{equation}
\label{lala}
 I_1(f)(\omega) = {<\omega, f>} = \displaystyle{\int_{\bR} f(s) dB_s (\omega) \hspace{0.2cm}  \mu-{\text {a.s}}.  }
\end{equation}

\subsection{Properties of Hermite functions and space $\sS'(\bR)$}
For every $n$ in $\bN$,  define

\begin{equation}
\label{iudezdzedrhiuzer}
h_n(x):={(-1)}^n e^{x^2}\frac{d^n}{dx^n}( e^{-x^2}) \h1 \hspace{0.2cm} \text{the} \h1 n^{th} \h1  \text{Hermite polynomial}
\end{equation}

and

\begin{equation}
\label{iudezedzedzezdzedrhiuzer}
e_n(x):=  {\pi}^{-1/4} {(2^n n!)}^{-1/2} e^{-x^2/2} h_{n}(x)  \h1 \hspace{0.2cm} \text{the} \h1 n^{th} \h1  \text{Hermite function}
\end{equation}

We will need the following properties of the Hermite functions:

\begin{theo}
\label{ozdicjdoisoijosidjcosjcsodijqpzeoejcvenvdsdlsiocfuvfsosfd}
\ben
\item The family ${(e_k)}_{k\in \bN}$ belongs to ${(\sS(\bR))}^{\bN}$ and  forms an orthonormal basis of $L^2(\bR)$ endowed with its usual inner product.

\item There exists a real constant $\widetilde{C}$ such that, for every $k$ in $\bN$, \h1  $\underset{x \in \bR }{\max}{|e_k(x)|} \leq \widetilde{C} \h1 {(k+1)}^{-1/12}$. More precisely, there exist positive constants $C$ and $\gamma$ such that,  for every $k$ in $\bN$,
\begin{equation}
\label{fkorep}
 {|e_k(x)|  } \leq
\begin{cases}
C \h1 {(k+1)}^{-1/12}   & {\text if} \hspace{0.25cm} |x| \leq 2\sqrt{k+1},\\
C e^{-\gamma x^2} &  {\text if} \hspace{0.25cm} |x| > 2\sqrt{k+1}. \\
\end{cases}
\end{equation}

\een
\end{theo}

See $\cite{Tha}$ for proofs.

In order to study precisely  $\sS(\bR)$ and its dual $\sS'(\bR)$ it is desirable to have a family of  norms on the space $\sS(\bR)$ which gives us the usual topology.

\begin{defi}
\label{oziqdorereredscqfegvelkunvcrveirvuyieuvynitrnyvtytvtriuyvntuiryvntrynv}
Let  ${({|\h1 \h1 |}_p)}_{p\in\bZ}$ be the family of norm defined by

\begin{equation}
\label{suidhcsdiuhsdihsdiuhcsdi}
 {|f|}^2_p:=   \sum^{+\infty}_{k=0} {(2k + 2)}^{2p} \h1 {<f,e_k>}^2_{L^2(\bR)}, \hspace{0.75cm}  \forall (p,f) \in \bZ \times L^2(\bR).
\end{equation}


 For $p$ in $\bN$, define the spaces $\sS_p(\bR):=\{f \in L^2(\bR), \h1 {|f|}_{p} <+\infty \}$ and  $\sS_{-p}(\bR)$ as being the completion of the space  $L^2(\bR)$ with respect to the norm ${{|\h1 \h1 |}_{-p}}$.
\end{defi}

\begin{rem}
For a function $f$  which is not in $L^2(\bR)$, we may still define ${|f|}_p$ by allowing ${|f|}_p$ to be infinite.
\end{rem}

It is well known (see \cite{Kuo2}) that the Schwartz space $\sS(\bR)$ is the projective limit of the sequence ${(\sS_{p}(\bR))}_{p \in \bN}$ and that the space $\sS'(\bR)$ of tempered distributions is the inductive limit of the sequence  ${(\sS_{\hspace{-0.15cm}-p}(\bR))}_{p \in \bN}$. This means firstly that both equalities $\sS(\bR) = \underset{p \in \bN }{\cap}{\sS_{p}(\bR)}$
  and $\sS'(\bR) = \underset{p \in \bN }{\cup}{{\sS}_{\hspace{-0.15cm}-p}(\bR)}$ hold. Secondly that convergence in ${\sS}(\bR)$ is nothing but convergence in ${\sS}_{p}(\bR)$ for every $p$ in $\bN$ and that convergence in ${\sS'}(\bR)$ is convergence in ${\sS_{\hspace{-0.15cm}-p}}(\bR)$ for some  $p$ in $\bN$.
Moreover one can show that, for any $p$ in $\bN$, the dual space ${\sS'}_{\hspace{-0.15cm}p}(\bR)$ of ${\sS}_{p}(\bR)$ is ${\sS_{\hspace{-0.15cm}-p}(\bR)}$. This is the reason why we will write ${\sS_{\hspace{-0.15cm}-p}(\bR)}$ in the sequel to denote the space  ${\sS'}_{\hspace{-0.15cm}p}(\bR)$.  Finally one can show that the usual topology of the space $\sS(\bR)$ and the topology given by the family of norms ${({|\h1 \h1 |}_p)}_{p\in\bN}$ are the same (see \cite{Hi} appendix A.3 for example).  Moreover, convergence in the inductive limit topology coincides with both convergence in the strong and the $\text{weak}^*$ topologies of ${\sS'}(\bR)$.

In view of definition \ref{oziqdorereredscqfegvelkunvcrveirvuyieuvynitrnyvtytvtriuyvntuiryvntrynv}, it is convenient to have an operator defined on $\sS(\bR)$ whose eigenfunctions are the sequence ${(e_n)}_{n \in \bN}$ and eigenvalues are the sequence $({{2n+2}})_{n \in \bN}$. It is easy to check that the operator  $A$, defined  on $\sS(\bR)$, by

\vspace{-0.4cm}
\begin{equation}
\label{frijsfdiojsdjsfdiojfsdoijsfd}
A:= -\frac{d^2}{dx^2} + x^2 +1
\end{equation}

verifies these conditions. Moreover $A$ is an isometry from $(\sS(\bR), {|\h1|}_{1})$ to $(L^2(\bR),{|\h1|}_{0})$.

 It may thus be extended to an operator, still denoted $A$, from  $\overline{\sS(\bR)}^{{|\h1|}_1} = \sS_1(\bR)$ to $L^2(\bR)$.  It is then easy to show (see \cite{Kuo2} p.17-18 for instance) 
 that $A$ is invertible and  that its inverse $A^{-1}$ is a bounded operator on $L^2(\bR)$. Let us note ${|g|}^2_0:={||g||}^2_{L^2(\bR)}$ for any  $g$ in $L^2(\bR)$. For $p$ in $\bZ$, let  $A^p$  denote the $p^{\text{th}}$ iteration of the operator $A$, if $p$ belongs to $\bN$, and of $A^{-1}$ otherwise. Then $\bD\text{om} (A^p) = \sS_p(\bR)$ and   $\bD\text{om} (A^{-p}) = L^2(\bR)$ where $\bD\text{om}(U)$ denotes the domain of the operator $U$ and $p$ belongs to $\bN$. Moreover, for every $q$ in $\bZ$ and every  $f:=\sum^{+\infty}_{k=0}  \h1 {<f,e_k>}_{L^2(\bR)} e_k$  
 in $\bD\text{om}  (A^q)$, the equality ${A^{q} f} =  \sum^{+\infty}_{k=0} {(2k + 2)}^{q} \h1 {<f,e_k>}_{L^2(\bR)} e_k $ holds. Hence,

\begin{equation}
 {|f|}^2_q = {|A^{q} f|}^{2}_{0} =  \sum^{+\infty}_{k=0} {(2k + 2)}^{2q} \h1 {<f,e_k>}^2_{L^2(\bR)}, \hspace{0.75cm}  \forall q \in \bZ.
\end{equation}

\subsection{Space of Hida distributions}

From now on  we will denote as is customary $(L^2)$ the space $L^2(\Omega,\cG,\mu)$ where $\cG$ is the $\sigma$-field generated by ${(<.,f>)}_{f \in L^2(\bR)}$.
Neither Brownian motion  nor fractional Brownian motion, whatever the value of $H$, are differentiable (see \cite{MVN} for a proof). However,
it occurs that the mapping $t\mapsto B^{(H)}(t)$ is differentiable from $\bR$ into a space, noted ${(\cS)}^*$, called the space of Hida distributions, which contains $(L^2)$. In this section we recall the construction of ${(\cS)}^*$.

 For every random variable $\Phi$ of $(L^2)$ there exists, according  to the Wiener-Itô theorem,  a unique sequence ${(f_n)}_{n \in \bN}$ of functions  $f_n$ in ${\widehat{L}}^2(\bR^n)$ such that $\Phi$ can be decomposed as

\begin{equation}
\label{oleole}
 \Phi =  \displaystyle{\sum^{+ \infty}_{n = 0} I_n(f_n)}
\end{equation}

where ${\widehat{L}}^2(\bR^n)$ denotes the set of all symmetric functions $f$ in $L^2(\bR^n)$ and $I_n(f)$ denotes the $n^{\text{th}}$ multiple Wiener-Itô integral of $f$ defined by

\begin{equation}
 I_n(f):= \int_{\bR^n} f(t) dB^{n}(t)= n! \int_{\bR} ( \int^{t_n}_{-\infty} \cdots \big(    \int^{t_2}_{-\infty} f(t_1,\cdots, t_n) dB(t_1)  \big) dB(t_2) \cdots dB(t_n)),
\end{equation}

 with the convention that $I_0(f_0) = f_0$ for constants $f_0$.  Moreover we have the isometry $E[\Phi^2]= \displaystyle{\sum^{+ \infty}_{n = 0} }n! \hspace{0.05cm}{{||f_n||}^2}_{L^2(\bR^n)}$. 
 For any $\Phi:=  \sum\limits^{+ \infty}_{n = 0}  I_n(f_n)$ satisfying the condition $\sum\limits^{+ \infty}_{n = 0}  n!\h1 {|A^{\otimes n}f_n|}^2_{0} < +\infty$,  define the element $\Gamma(A)(\Phi)$ of $(L^2)$ by

\begin{equation}
\Gamma(A)(\Phi):= \sum\limits^{+ \infty}_{n = 0} I_n(A^{\otimes n}f_n),
\end{equation}

where $A^{\otimes n}$ denotes the $n^{th}$ tensor power of the operator $A$ defined in \eqref{frijsfdiojsdjsfdiojfsdoijsfd} (see \cite{janson} appendix E for more details about tensor products of operators). The operator $\Gamma(A)$ is densely defined on $(L^2)$ and is called the second quantization operator of $A$. It shares a lot of properties with the operator $A$. In particular it is invertible and its inverse  ${\Gamma(A)}^{-1}$ is bounded (see \cite{Kuo2}). Let us denote ${||\varphi||}^2_0:={||\varphi||}^2_{(L^2)}$ for any random variable $\varphi$ in $(L^2)$ and, for $n$ in $\bN$, let $\bD\text{om}({\Gamma(A)}^n)$ be the domain of the  $n^{\text{th}}$ iteration of $\Gamma(A)$. The space of Hida distributions is defined in a way analogous to the one that allowed  to define the  space $\sS'(\bR)$:

\begin{defi}
Define the  family  of norms  ${({||\h1||}_p)}_{p \in \bZ}$ by:

\begin{equation}
\label{erofij}
  {||\Phi||}_p :=  {||\Gamma(A)^p\Phi||}_{0} =  {||\Gamma(A)^p\Phi||}_{(L^2)},  \hspace{1cm} \forall p \in \bZ,\hspace{0.5cm} \forall\Phi \in (L^2)\cap \bD{\text{om}}({\Gamma(A)}^p).
\end{equation}

For any $p$ in $\bN$, let $({\cS}_{p}):=\{\Phi \in (L^2): \h1  \Gamma(A)^p\Phi \h1 \text{exists and belongs to}  \h1 (L^2) \}$ and define
$({\cS}_{-p})$ as being the completion of the space  $(L^2)$ with respect to the norm ${{||\h1 \h1 ||}_{-p}}$.
\end{defi}

As in \cite{Kuo2}, we let $(\cS)$ denote the projective limit of the sequence ${( (\cS_{p}))}_{p \in \bN}$ and ${(\cS)}^*$ the inductive limit of the sequence  ${(({\cS_{-p}}))}_{p \in \bN}$. Again this means that we have the equalities $(\cS) = \underset{p \in \bN }{\cap}{(\cS_{p})}$ (resp.  ${(\cS)}^* = \underset{p \in \bN }{\cup}{({\cS_{-p}})}$) and  that convergence in $(\cS)$ (resp. in ${(\cS)}^*$) means convergence in ${({\cS}_{p})}^*$ for every $p$ in $\bN$   (resp. convergence in  $({\cS_{-p}})$ for some $p$ in $\bN$ ).
The space  $(\cS)$ is called the space of stochastic test functions and ${(\cS)}^*$  the space of Hida distributions. As previously one can show that, for any $p$ in $\bN$, the dual space ${({\cS}_{p})}^*$ of ${\cS_{p}}$ is  $({\cS_{-p}})$. Thus we  will write $({\cS_{-p}})$, in the sequel, to denote the space  ${({\cS}_{p})}^*$.
Note also that ${(\cS)}^*$ is the dual space of $(\cS)$. We will note  $<\hspace{-0.2cm}<\hspace{-0.1cm} \h1,\h1 \hspace{-0.1cm}>\hspace{-0.2cm}>$  the duality bracket between ${(\cS)}^*$ and $(\cS)$. If $\Phi$ belongs to $(L^2)$ then we have the equality $\begin{g}{\Phi}{\varphi}\end{g} = {<\Phi,\varphi>}_{(L^2)} = \bE[\Phi \h1 \varphi]$. Furthermore, as one can check,  the family  ${({|f|}_p)}_{p \in \bZ}$ is an increasing sequence for every $f$ in $\sS(\bR)$. Thus the family  ${({||<.,f>||}_p)}_{p \in \bZ}$ is an increasing sequence for every $f$ in $\sS(\bR)$.

\begin{rem}
\label{laiusdhv}
$(i)$ A consequence of the previous subsection is that for every element $f:= {\sum^{+ \infty}_{n = 0} \h1 a_n e_n}$ in $\sS'(\bR)$ where  ${(a_n)}_{n \in \bN}$ belongs to ${\bR}^{\bN}$, there exists $p_0$ in $\bN$ such that $f$ belongs to  $\sS_{\hspace{-0.15cm}-p_0}(\bR)$. Moreover if we define $\Phi:={\sum^{+ \infty}_{n = 0} \h1 a_n <.,e_n>}$, then $\Phi$ belongs to  $({\cS_{-p_0}}) \subset {(\cS)}^*$ and we have  ${||\Phi||}^2_{-p_0} = {|f|}^2_{-p_0} = {\sum^{+ \infty}_{n = 0}  \h1 \frac{a^2_n}{{(2n+2)}^{2p_0}} } < +\infty$.
Conversely, every element $\Phi$,  written as $\Phi:={\sum^{+ \infty}_{n = 0} \h1 b_n <.,e_n>} $ where ${(b_n)}_{n \in \bN}$ belongs to ${\bR}^{\bN}$,  belongs to ${(\cS)}^*$ if and only if there exists an integer $p_0$ such that ${\sum^{+ \infty}_{n = 0}  \h1 \frac{b^2_n}{{(2n+2)}^{2p_0}} } < +\infty$. In this case the element $f:= {\sum^{+ \infty}_{n = 0} \h1 b_n e_n}$   belongs to  $\sS_{\hspace{-0.15cm}-p_0}(\bR)$ and hence is a tempered distribution which verifies ${|f|}^2_{-p_0} = {||\Phi||}^2_{-p_0} = {\sum^{+ \infty}_{n = 0}  \h1 \frac{b^2_n}{{(2n+2)}^{2p_0}} }$.

$(ii)$ Let ${\widehat{\sS}}({\bR}^n)$ be the space of symmetric Schwartz functions defined on $\bR^n$. Let $p$ in $\bN^*$ and $\Phi$ be an element of $({\cS}_{-p})$. Then  $\Phi$ can be written $\Phi:=\sum^{+ \infty}_{n = 0}  \h1 I_n(F_n)$, where  ${(F_n)}_{n \in \bN}$ is a sequence in ${({\widehat{\sS}} {\textcolor{white}{c}}^{'}({\bR}^n))}^{\bN}$, where ${\widehat{\sS}} {\textcolor{white}{c}}^{'}({\bR}^n)$ is the dual of ${\widehat{\sS}}({\bR}^n)$.
%
%
%
\end{rem}

Since we have defined a topology given by a family of norms on the space ${(\cS)^*}$ it is possible to define a derivative and an integral in ${(\cS)^*}$ (see \cite{HP} chapter $3$ for more details about these notions). Let $I$ be an  interval of $\bR$ (which may  be equal to $\bR$).

\begin{defi}[{\bfseries stochastic distribution process}]
\label{odifjfjpqjkqpcsjkkpdskpdspokcdspocksdpocskpdok}
\medskip
A function $\Phi:I \rightarrow {(\cS)}^{*}$ is called a stochastic distribution process, or an ${(\cS)}^{*}\hspace{-0.1cm}-$process, or a Hida
process.
\end{defi}

\begin{defi}[{\bfseries derivative in ${(\cS)}^{*}$}]
\label{odifjfjpcskpdok}
\medskip
Let $t_0 \in I$. A stochastic distribution process $\Phi:I \rightarrow {(\cS)}^{*}$ is said to be differentiable at  $t_0$ if the quantity
{$\lim\limits_{r \to 0}~\h1~r^{-1}\h1(\Phi(t_0+r)~-~\Phi(t_0))$} exists in ${(\cS)}^{*}$. We note $\frac{d\Phi}{dt}(t_0)$ the ${(\cS)}^{*}$\hspace{-0.1cm}-derivative at $t_0$ of the stochastic distribution process $\Phi$. $\Phi$ is said to be differentiable over $I$ if it is differentiable at $t_0$ for every $t_0$ in $I$.
\end{defi}

The process $\Phi$ is said to be continuous, $\cC^1,\cdots,\cC^k, \cdots$ in ${(\cS)}^{*}$ if the ${(\cS)}^{*}$-valued function $\Phi$ is,  continuous, $\cC^1,\cdots,\cC^k, \cdots$. We also say that the stochastic distribution process  $\Phi$ is ${(\cS)}^{*}$-continuous and so on.
It is also possible to define an ${(\cS)}^{*}$-valued integral in the following way (\cite{Kuo2,HP}). We first recall that $L^1(\bR, dt)$ denotes the set of measurable complex-valued functions defined on $\bR$ such that ${||f||}_{L^1(\bR)} := \int_{\bR} |f(t)| \h1 dt < +\infty$.

\begin{theodef}[{\bfseries integral in  ${(\cS)}^{*}$}]
Assume that $\Phi:\bR \rightarrow {(\cS)}^{*}$ is weakly in $L^{1}(\bR,dt)$, i.e assume that for all $\varphi$ in  ${(\cS)}$, the mapping $u \mapsto <\hspace{-0.2cm}<\hspace{-0.1cm} \h1 \Phi(u),\h1 \varphi \hspace{-0.1cm}>\hspace{-0.2cm}>$  from $\bR$ to $\bR$ belongs to  $L^{1}(\bR,dt)$. Then there exists an unique element in ${(\cS)}^{*}$, noted $\int_{\bR} \Phi(u) du $ such that

\begin{equation}\label{foijfeorj}
{<\hspace{-0.2cm}<\hspace{-0.1cm} \h1 \int_{\bR} \Phi(u) du,\varphi \h1 \hspace{-0.1cm}>\hspace{-0.2cm}>} = \int_{\bR} <\hspace{-0.2cm}<\hspace{-0.1cm} \h1 \Phi(u),\varphi\h1 \hspace{-0.1cm}>\hspace{-0.2cm}> \h1 du \hspace{0.5cm }\text{for all } \varphi \hspace{0.1cm }\text{in} \hspace{0.1cm }{(\cS)}.
\end{equation}

\end{theodef}

We say in this case that $\Phi$ is ${(\cS)}^{*}$-integrable on $\bR$  in the {\it Pettis sense}. In the sequel, when we do not specify a name for the integral of an ${(\cS)}^{*}$-integrable process $\Phi$ on $\bR$, we always refer to the integral of $\Phi$ in Pettis' sense. 
See \cite{Kuo2} p.$245$-$246$ or \cite{HP} def. $3.7.1$ p.$77$ for more details.

Recall from \eqref{lala} that $I_1(f) = \int_{\bR} \h1 f(s) \h1 dB(s)$.  For $p$ in $\bN$, let  $I^{(p)}_1:(L^2(\bR), {|\h1 \h1|}_{-p}) \rightarrow ((L^2),{||\h1 \h1||}_{-p})$ be the map defined by  $I^{(p)}_1(f):= <.,f>$.
For $f$ in $L^2(\bR)$ and $p$ in $\bN$ we have, using \cite{Kuo2} p $26$,  the equality ${||I^{(p)}_1(f)||}_{-p} = {||{\Gamma(A)}^{-p}(I_1(f))||}_{0} = {||I_1(A^{-p}(f))||}_{0}  = {|A^{-p}(f)|}_{0} =  {|f|}_{-p}$.

Hence $I^{(p)}_1$ is an isometry and we can extend it to  $\sS_{-p}(\bR)$. Since we can do this for every integer $p$ in $\bN$, we may then give a meaning to $I_1(f) = <.,f>$  for every tempered distribution $f$ as being the element $I^{(p)}_1(f)$  where $f$ belongs to  $\sS_{-p}(\bR)$.



\begin{rem}
\label{zpopdokokdkkdkdkdkdkdkdkdk}
We will also note  $\int_{\bR} \h1 f(s) \h1 dB(s)$ the quantity $<.,f>$ when $f$ belongs to $\sS'(\bR)$. Hence we give a meaning to the quantity $\int_{\bR} \h1 f(s) \h1 dB(s)$  for every $f$ in $\sS'(\bR)$.
\end{rem}

\subsection{S-transform and Wick product  }
\label{ksksksks}

For $\eta$ in $\sS(\bR)$, the {\it Wick exponential} of $<.,\eta>$, denoted $:e^{<.,\eta>}:$, is defined as the element of $(\cS)$ given by $:e^{<.,\eta>}:\h1  \stackrel{\text{def}}{=} \sum^{+\infty}_{k = 0}  {k!}^{-1} \h1 I_k(\eta^{\otimes k})$ (equality in $(L^2)$). More
generally, for $f \in L^2(\bR)$, we define  $:e^{<.,f>}:$ as the $(L^2)$ random variable equal to
$e^{<.,f> - \frac{1}{2}{|f|}^2_0}$ (see \cite{janson} theorem 3.33).
We will sometimes note $\exp^\diamond {<.,f>}$ instead of $:e^{<.,f>}:$.
  This random variable belongs to $L^p(\Omega, \mu)$ for every integer $p\geq 1$. We now recall the definition of the $S$-transform of an element $\Phi$ of $(\cS^*)$, noted $S(\Phi)$ or $S[\Phi]$. $S(\Phi)$ is defined as the function from  $\sS(\bR)$ to $\bR$ given by

\begin{equation}
\forall \eta \in \sS(\bR), \hspace{0.6cm} S(\Phi)(\eta):= {\begin{g}{\Phi}{ :e^{<.,\eta>}:\hspace{0.075cm}}  \end{g} }.
\end{equation}

Note that $S\Phi(\eta)$ is nothing but $\bE[\Phi :e^{<.,\eta>}: ] =  e^{ -\frac{1}{2}{|\eta|}^2_0} \h1 \bE[\Phi \h1 e^{<.,\eta>} ]$ when $\Phi $ belongs to $(L^2)$. Following \cite{ben2}, formula $(6)$ and $(7)$, define for $\eta$ in $\sS(\bR)$ the probability measure $\bQ_{\eta}$ on the space $(\Omega,\cF)$ by its  Radon-Nikodym derivative given by $\frac{d\bQ_{\eta}}{d\mu} \stackrel{\text{def}}{=} :e^{<.,\eta>}$:. The probability measures  $\bQ_{\eta}$ and $\mu$ are equivalent. Then, by definition,

\begin{equation}
\label{peropdockdspocdspocksdpokcpsdokcsdpocksdpockdspockds}
\forall \Phi \in (L^2), \hspace{0.6cm}  S(\Phi)(\eta) = {\bE}_{\bQ_{\eta}}[\Phi].
\end{equation}

%

%
%
%

\begin{lem}
\label{dpozepdozkekpdzkeodkpozekpdozkepdozekodezpdkezpd}

\ben
\item
\label{lsmqlqsmlsqqlmslqmlslmjoediedh}
Let $p$ be a positive integer and $\Phi$ be an element of $({\cS}_{-p})$. Then

\begin{equation}
\label{estime}
|S(\Phi)(\eta)| \leq {||\Phi||}_{-p} \h1 e^{\frac{1}{2} {|\eta|}^2_{p}}, \h1 \text{for any} \h1 \eta  \h1 \text{in} \h1  \sS(\bR).
\end{equation}

\item Let $\Phi:=  \sum^{+\infty}_{k = 0} a_k <.,e_k>$ belong to ${(\cS)}^*$. The following equality holds for every $\eta$ in $\sS(\bR)$:

\begin{equation}
S(\Phi)(\eta) = \sum^{+\infty}_{k = 0} a_k {<\eta,e_k>}_{L^2(\bR)}.
\end{equation}
\een

\end{lem}

\begin{pr}
Item 1 is proved in \cite{Kuo2} p.$79$. Item 2 is an easy calculation left to the reader.
\end{pr}

%
%
%
%

Another useful tool in white noise analysis is the Wick product:

\begin{theodef}[\cite{Kuo2} p.$92$]
\label{wickproducte}
For every $(\Phi,\Psi)\in {(\cS)}^{*}\times {(\cS)}^{*}$, there exists a unique element of ${(\cS)}^{*}$, called the Wick product of  $\Phi$ and $\Psi$  and  noted $\Phi \diamond \Psi$,  such that, for every $\eta$ in $\sS(\bR)$,

\begin{equation}
\label{mszzmm}
S(\Phi\diamond \Psi)(\eta) = S(\Phi)(\eta) \h1 S(\Psi)(\eta).
\end{equation}
\end{theodef}

\begin{lem}
\label{dede}
For any $(p,q) \in \bN^2$, $X \in ({\cS}_{-p})$ and $Y\in ({\cS}_{-q})$,

\begin{equation}
|S(X \diamond Y)(\eta)| \leq {||X||}_{-p} \h1 {||Y||}_{-q} \h1 e^{{|\eta|}^2_{\max \{p;q\}}}.
\end{equation}
\end{lem}

\begin{pr}
The proof is easy since, for every $\eta$ in $\sS(\bR)$, using lemma \ref{dpozepdozkekpdzkeodkpozekpdozkepdozekodezpdkezpd},

 $|S(X \diamond Y)(\eta)| = |S(X)(\eta)| \h1 |S(Y)(\eta)| \leq  {||X||}_{-p} \h1  e^{\frac{1}{2}{|\eta|}^2_p} \h1 {||Y||}_{-q} \h1  e^{\frac{1}{2}{|\eta|}^2_q} \leq {||X||}_{-p} \h1 {||Y||}_{-q} \h1  e^{{|\eta|}^2_{\max \{p;q\}}}$.
\end{pr}

For any $\Phi$ in ${(\cS)}^{*}$ and $k$ in $\bN^*$ let $\Phi^{\diamond k}$ denote the element $\overbrace{\Phi \diamond \cdots  \diamond \Phi}^{k \text{ times}}$  of $(\cS)^{*}$. We can generalize the definition of $\exp^{\diamond}$ to the case where $\Phi$ belongs to $(\cS)^{*}$:

\begin{defi}
\label{mqmqm}
For any $\Phi$ in ${(\cS)}^{*}$ such that the sum $\sum^{+\infty}_{k = 0} \frac{\Phi^{\diamond k}}{k!}$ converges in $(\cS)^{*}$, define  the element $\exp^{\diamond}\Phi$ of  ${(\cS)}^{*}$ by $\exp^{\diamond}\Phi:= \sum^{+\infty}_{k = 0} \frac{\Phi^{\diamond k}}{k!}$.
\end{defi}

For $f$ in $L^2(\bR)$ and  $\Phi:= <.,f>$, it is easy to verify that $\exp^{\diamond} \Phi$ given by definition \ref{mqmqm} exists and coincide with $:e^{<.,f>}:$ defined at the beginning of this section.

\begin{rem}
\label{qmqmqmmmmm}
If $\Phi$ is deterministic then, for all $\Psi$ in ${(\cS)}^{*}$, $\Phi \diamond \Psi = \Phi \Psi$. Moreover, let ${(X_t)}_{t \in \bR}$ be a Gaussian process and let $\cH$ be the subspace of $(L^2)$ defined by $\cH := \overline{  \text{vect}_{\bR}\{ X_t;t \in \bR \}}^{(L^2)}$. If $X$ and $Y$ are two elements of $\cH$ then $X\diamond Y = XY - \bE[XY]$.
\end{rem}

%

%
%

We refer to \cite{janson} chapters $3$ and $16$ for more details about Wick product. The following results on the S-transform will be  used in the sequel. See \cite{Kuo2} p.$39$ and \cite{HKPS} p.280-281 for proofs.

\begin{lem}
\label{dkdskcsdckksdksdmksdmlkskdm}
The $S$-transform verifies the following properties:
\bit
\iti The map $S:\Phi\mapsto S(\Phi)$, from ${(\cS)}^{*}$ to ${(\cS)}^{*}$, is injective.

\itii Let $\Phi:\bR \rightarrow {(\cS)}^{*}$  be an ${(\cS)}^{*}$ process. If $\Phi$ is ${(\cS)}^{*}$-integrable over $\bR$ then

$\displaystyle{S(\int_{\bR} \Phi(u) \h1 du )(\eta) = \int_{\bR} S(\Phi(u)) (\eta) \h1 du}$ , for all $\eta$ in $\sS(\bR)$.
\itiii  Let $\Phi:\bR \rightarrow {(\cS)}^{*}$  be an ${(\cS)}^{*}$-process differentiable at $t$. Then, for every $\eta$ in $\sS(\bR)$ the map $u\mapsto [S \Phi(u)](\eta)$ is differentiable at $t$ and verifies $\displaystyle{S[\tfrac{d\Phi}{dt}(t)](\eta) = \tfrac{d}{dt}\big[ S[ \Phi(t)](\eta) \big]}$.
\eit
\end{lem}

It is useful to have a criterion for integrability  in ${(\cS)}^{*}$ in term of the $S$-transform. This is the topic of the next theorem (theorem 13.5 in \cite{Kuo2}).

\begin{theo}
\label{peodcpdsokcpodfckposkcdpqkoq}
Let $\Phi:\bR \rightarrow {(\cS)}^*$ be a stochastic distribution process satisfying:
\bit
\iti The map $t\mapsto S[\Phi(t)](\eta)$, from $\bR$ to $\bR$,  is measurable for all $\eta$ in $\sS(\bR)$.
\itii There is a natural integer $p$, a real $a$ and a  function $L$ in $L^1(\bR,dt)$ such that for all  $\eta$ in $\sS(\bR)$,
 $ |S(\Phi(t))(\eta)| \leq L(t) \h1 e^{a {| \eta |}^2_{p} }$.
\eit
Then $\Phi$ is ${(\cS)}^{*}$-integrable over $\bR$. 
%
\end{theo}

Lastly, when the stochastic distribution process is an $(L^2)$-valued process, the following result holds (see \cite{ben2}):

\begin{theo}
\label{ozeijdoezjdozeijdoizjdoijezoijzdijz}
Let $X:\bR \rightarrow (L^2)$ be such that the function $t\mapsto S(X_t)(\eta)$ is measurable for all $\eta$ in $\sS(\bR)$ and that $t \mapsto {||X_t||}_0$ is in $L^1(\bR,dt)$. Then $X$ is ${(\cS)}^*$-integrable over $\bR$ and

\begin{equation*}
{\left|\left| \int_{\bR} \h1 X_t \h1 dt \right|\right|}_0 \leq \int_{\bR} \h1 {||X_t||}_0  \h1 dt.
\end{equation*}
\end{theo}

\section{The operators $M_H$ and their derivatives}
\label{MH}

\subsection{Study of $M_H$}
\label{msldkldkmldklsmmks}

Let us fix some notations. We will still note $\widehat{u}$ or $\cF(u)$ the Fourier transform of a tempered distribution $u$ and we let $L^{1}_{loc}(\bR)$ denote the set of measurable functions which are locally integrable on $\bR$.
 We also identify, here and in the sequel,  any function $f$ of $L^{1}_{loc}(\bR)$ with its associated distribution, also noted $T_f$. We will say that a tempered distribution $v$ is of function type if there exists a locally integrable function $f$ such that $v=T_f$ (in particular, $<v,\phi> = \int_{\bR} f(t) \h1 \phi(t) \h1 dt$  for $\phi$ in  $\sS(\bR)$).

Let $H \in (0,1)$. Following \cite{ell}, we want to define an operator,  denoted $M_H$, which is specified in the Fourier domain by

\begin{equation}
\label{qxs}
\widehat{M_H(u)}(y) := \frac{\sqrt{2\pi}}{c_H} \hspace{0.1cm}{|y|}^{1/2-H} \widehat{u}(y), \hspace{1cm} y \in \bR^*.
\end{equation}

This operator is well defined on the homogeneous Sobolev space of order $1/2-H$, $L^2_H(\bR)$:

 \label{ninth}
 \begin{equation}
L^2_H(\bR):= \{u \in \mathcal{\sS'(\bR)} \hspace{0.1cm}:\hspace{0.1cm} \widehat{u} = T_{f}; \hspace{0.1cm}  f \in L^{1}_{loc}(\bR) \hspace{0.1cm}\text{and}\hspace{0.1cm}  {|| u ||}_{H} <  +\infty  \},
\end{equation}

where ${{|| u ||}^2}_{H}:= {\frac{1}{c^2_H}\int_{\bR} {|\xi|}^{1-2H} {|\widehat{u \hspace{0.1cm}}(\xi) |}^2  d{\xi}     }$\hspace{0.1cm}   derives from the inner product on $L^2_H(\bR)$, defined by:

%
\begin{equation}
\label{alban}
{<u,v>}_{H} := \frac{1}{c^2_H} \int_{\bR} {|\xi|}^{1-2H} {\widehat{u \hspace{0.1cm}}(\xi) } \overline{{\widehat{v \hspace{0.1cm}}(\xi) }}   d{\xi},
\end{equation}


and $c_{H}$ has been defined right after formula (\ref{covmbm}) (the normalization constant $ \frac{\sqrt{2\pi}}{c_H}$ will be explained in remark \ref{remme}). It is well known - see  \cite{che} p.$13$ for example - that $(L^2_H(\bR),{<,>}_{H})$ is a Hilbert space. The nature of the spaces $L^2_H(\bR)$ when $H$ spans $(0,1)$ is described in the following lemma, the proof of which can be found in \cite{che} p$15$, theorem $1.4.1$ and corollary $1.4.1$.
\begin{lem}
\label{qlkjqlkj}
If $H$ is in  $(0,1/2]$, the space $L^2_H(\bR)$ is  continuously embedded in $L^{1/H}(\bR)$. When $H$ is in  $[1/2,1)$, the space $L^{1/H}(\bR)$ is  continuously embedded in $L^2_H(\bR)$.
\end{lem}

Since $\widehat{M_H(u)}$ belongs to $L^2(\bR)$ for every $u$ in $L^2_H(\bR)$,  $M_H$ is well defined as its inverse Fourier transform, {\it i.e.}:

\begin{equation}
\label{uyguyygijopoji}
{M_H}(u)(x):= \frac{1}{2\pi}\cF(\widehat{M_H(u)})(-x),  \hspace{0.5cm}  \text{for almost every}\h1 x \h1 \text{in} \h1 \bR.
\end{equation}


The following proposition is obvious in view of the definition of $M_H$:

\begin{prop}
\label{a la bonne franquette}
$M_H$ is an isometry from $(L^2_H(\bR),{<,>}_{H})$ to  $(L^2(\bR),{<,>}_{L^2(\bR)} )$.
\end{prop}

Let $\cE(\bR)$ denote the space of simple  functions on $\bR$, which is the set of all finite linear combinations of functions $\i1_{[a,b]}(.)$ with $a$ and $b$ in $\bR$. It is easy to check that both $\sS(\bR)$ and $\cE(\bR)$ are subsets of $L^2_H(\bR)$.

It will be useful in the sequel to have an explicit expression for $M_H(f)$ when $f$ is in $\sS(\bR)$ or in $\cE(\bR)$. To compute this value, one may
use the formulas for the Fourier transform of the distributions  ${|\h1 |}^{\alpha}$, $\alpha$ in $(-1,1)$, given for instance in \cite{GeCh} (chapter 1, \§$3$). This yields, for almost every $x$ in $\bR$,

\begin{equation}\label{referq}
  {\mathtitre1 {M_H(\i1_{[a,b]}) (x) = \tfrac{\sqrt{2\pi}}{2 c_H \Gamma(H+1/2)\cos(\frac{\pi}{2}(H-1/2))} \hspace{0.1cm} \bigg[  \tfrac{b - x}{{|b-x|}^{3/2 - H}} - \tfrac{a - x}{{|a-x|}^{3/2 - H}}\bigg] }}.
\end{equation}

 By the same method, for $f$ in $\sS(\bR)$ one gets, for almost every real $x$:

\begin{align}\label{referque}
	{\text{for}}  \hspace{0.1cm} & 0< H < 1/2	\hspace{1.1cm} & M_H(f) (x) &= \gamma_H \int_{\bR} \frac{f(x-t) - f(x)}{ {|t|}^{3/2 - H}  } dt \\
	\label{referquee}
 	{\text{for}}  \hspace{0.1cm} & H = 1/2  	\hspace{1.1cm} &M_H(f) (x) &= f(x) \\
 	\label{referqueee}
 	{\text{for}}  \hspace{0.1cm} & 1/2< H < 1	\hspace{1.1cm} &M_H(f) (x) &= \gamma_H \int_{\bR} \frac{f(t)}{ {|t-x|}^{3/2 - H}  } dt
\end{align}

where $\gamma_H :=  \tfrac{\sqrt{2\pi}}{2 c_H \Gamma(H-1/2)\cos(\frac{\pi}{2}(H-1/2))}  =\tfrac{{ { \big ( \Gamma(2H+1)\sin(\pi H)}  \big )  }^{\frac{1}{2}}}{2 \Gamma(H-1/2)\cos(\frac{\pi}{2}(H-1/2))}$.  When $f$ belongs to $\sS(\bR)$, equality (A.$1$) of \cite{ell} yields (up to a constant), for almost every $x$ in $\bR$,

 \begin{equation}
  \label{referqueosijz}
M_H(f) (x) = \alpha_H  \frac{d}{dx} \left[\int_{\bR} (t-x){|t-x|}^{H-3/2} f(t) dt\right]
 \end{equation}

\begin{equation}\label{alpha_H}
\textrm{where} \qquad \alpha_H:= \frac{-\gamma_H}{H-1/2} = \frac{-\sqrt{2\pi}}{2 c_H \Gamma(H+1/2)\cos(\frac{\pi}{2}(H-1/2))}.
\end{equation}
In order to extend the Wiener integral with respect to fBm to an integral with respect to mBm (in section \ref{Wiener}) we will need the following equality:

\begin{prop}
$\overline{\cE(\bR)}^{{<,>}_H}  = L^2_H(\bR)$.
\end{prop}

This is a straightforward consequence of the following lemma:

\begin{lem}
\label{bcv}
Let $\sigma:\bR \rightarrow \bC$ be a measurable function, continuous on $\bR^*$, such that ${|\sigma|}^2$ is locally integrable at $0$ and that $x \mapsto {\left|\tfrac{{\sigma(x)}}{x}\right|}^2$ is locally integrable at $+\infty$. Define $L^2_{\sigma}(\bR):= \{u \in \mathcal{\sS'(\bR)} \hspace{0.1cm} : \hspace{0.1cm} \widehat{u} = T_{f}; \h1 f \in L^{1}_{loc}(\bR)\h1 \hspace{0.1cm}  \text{such that} \h1 {|| u ||}_{\sigma} <  +\infty  \}$ where ${<u,v>}_{\sigma} :=  \int_{\bR} {|\sigma(\xi)|}^{2} \h1 \widehat{u}(\xi) \h1 \overline{{\widehat{v} \h1}(\xi) } \h1  d{\xi}$.
If $\cE(\bR) \subset L^2_{\sigma}(\bR)$, define $\overline{\cE(\bR)}^{{<,>}_{\sigma}}$ as the completion of $\cE(\bR)$ for the norm ${||\h1||}_{\sigma}$. Then, the space $(L^2_{\sigma}(\bR),{<,>}_{\sigma})$ is a Hilbert space which also verifies $\overline{\cE(\bR)}^{{<,>}_{\sigma}} = L^2_{\sigma}(\bR)$.
\end{lem}

\begin{pr}
The fact that $(L^2_{\sigma}(\bR),{<,>}_{\sigma})$ is a Hilbert space is obvious.
 One needs only to show that the orthogonal space of $\cE(\bR)$ for the norm ${||\h1||}_{\sigma}$ is equal to $\{0_{\cE(\bR)}\}$.
 Let $u$ in $L^2_{\sigma}(\bR)$ be such that  ${<u,v>}_{\sigma} = 0$ for all  $v$ in $\cE(\bR)$. In particular, for all $t$ in $\bR$, $\int_{\bR} {|\sigma(\xi)|}^2 \h1 \widehat{u}(\xi) \h1 \overline{\widehat{ \i1_{[0,t]}}(\xi)}\h1 d\xi~=~0$. For all  ${\psi}$ in $\sS(\bR)$,

\begin{equation*}
\displaystyle{\int_{\bR} \psi' (t)  \bigg(\int_{\bR}  {|\sigma(\xi)|}^2 \h1 \widehat{u}(\xi) \h1    \overline{\widehat{ \i1_{[0,t]}}(\xi)} d\xi \bigg)}dt  = 0,
\end{equation*}

where $\psi'$  denotes the derivative of $\psi$. Thanks to the assumptions  on ${|\sigma|}^2$ and  $x \mapsto {\left|\tfrac{{\sigma(x)}}{x}\right|}^2$, Fubini theorem applies. Moreover, an integration by parts yields

\begin{equation*}
 0 = -\int_{\bR} \frac{{|\sigma(\xi)|}^2}{i\xi} \h1  \widehat{u}(\xi)     \bigg( \int_{\bR} \psi' (t) (1 - e^{i \xi t})    dt  \bigg) d\xi =  \int_{\bR} {|\sigma(\xi)|}^2 \h1  \widehat{u}(\xi)  \h1    \overline{\widehat{\psi}( \xi)} \h1  d\xi.
\end{equation*}

Thus $<{|\sigma|}^2 \widehat{u}, \psi> = 0$ for all $\psi$ in  $\sS(\bR)$. Since $ \xi \mapsto {|\sigma (\xi)|}^2 \widehat{u}(\xi)$ belongs to $L^1_{loc}(\bR)$, it is easy to deduce that $u$ is equal to $0$.
\end{pr}

\begin{rem}
\label{remme}
 1. Because the space  $\sS(\bR)$ is dense  in $L^2_H(\bR)$ for the norm ${||\h1\h1||}_H$ (see \cite{che} p.$13$), it is also possible to define the operator $M_H$ on the space $\sS(\bR)$ and extend it, by isometry, to all elements of $L^2_H(\bR)$. This is the approach of  \cite{ell} and  \cite{bosw} (with a different normalization constant). This clearly yields the same operator as the one defined by \eqref{qxs}. However this approach does not lend itself to an extension to the case where the constant $H$ is replaced by a function $h$, which is what we need for mBm. \\
2. For the same reasons as in 1. it is  possible to define the operator $M_H$ on the space $\cE(\bR)$ and extend it, by isometry, to all elements of $L^2_H(\bR)$. Again, this extension coincide with \eqref{qxs}. We will use this idea in section \ref{Wiener}.
\end{rem}

In view of \eqref{alban}, we find that

\begin{equation}
\label{seventh}
<\i1_{[0,t]},\i1_{[0,s]}>_H  =\frac{1}{ c^2_H} \int_{\bR}  \frac{{(e^{it\xi} - 1)(e^{-is\xi} - 1) }}{{|\xi|}^{2H+1 }} d\xi =  R_H(t,s).
\end{equation}


 As in the case of standard Brownian motion, one deduces that the process ${({{\widetilde{B}}^{(H)}}(t))}_{t \in \bR}$, defined for all  $(t,\omega)$ in $\bR\times\Omega$ by:

\begin{equation}
\label{fifth}
 {{\widetilde{B}}^{(H)}}(t) (\omega):= {\widetilde{B}^{(H)}}(t,\omega):= {<\omega,  M_H(\i1_{[0,t]})>},
\end{equation}

is a Gaussian process which admits, as the next computation shows, a continuous version noted $B^{(H)}:={(B^{(H)}(t))}_{t \in \bR}$. Indeed, under the probability measure $\mu$, the process $B^{(H)}$ is a fractional Brownian motion since we have, using \eqref{seventh} and proposition \ref{a la bonne franquette},

\begin{align}\label{ezdzedeferf}
\bE[B^{(H)}(t) B^{(H)}(s)] &= \bE[<.,  M_H(\i1_{[0,t]})> <.,  M_H(\i1_{[0,s]})> ] \notag \\
 								 &=  <  M_H(\i1_{[0,t]}), M_H(\i1_{[0,s]})>_{L^2(\bR)} \notag \\
 								 &=  {<\i1_{[0,t]},\i1_{[0,s]} >}_{H} =  R_H(t,s).
\end{align}

\begin{rem}
\label{reme}
 The reason of the presence of the constant $\frac{\sqrt{2\pi}}{{c_H}}$ in formula {\normalfont (\ref{qxs})} is now clear since this constant ensures  that, for all $H$ in $(0,1)$, the process $B^{(H)}$ defined by $(\ref{fifth})$ is a normalized fBm.
  \end{rem}

Because our operator $M_H$ is defined on a distribution space, we can not apply  the considerations of  \cite{ell} p.$323$ff about the links between the operator $M_H$ and Riesz potential operator.
However it is crucial for our purpose that $M_H$ is bijective from $L^2_H(\bR)$ into $L^2(\bR)$:

\begin{theo}[properties of $M_H$]
\label{ozdicjdoisoedzijqpzeoejcvenvdsdlsiocfuvfsosfd}
\ben
\item For all $H$ in $(0,1)$, the operator $M_H$ is bijective from $L^2_H(\bR)$ into $L^2(\bR)$.
 \vspace{-0.25cm}

\item For all $H$ in $(0,1)$ and $(f,g)$ in ${(L^2(\bR) \cap L^2_H(\bR))}^2$,

\begin{equation}
\label{dfufdfdiuhuhuidsfhufdshudsfhdfsuhdsfifsdhdsf}
{  {<f,M_H(g)>}_{L^2(\bR)} = {<M_H(f),g>}_{L^2(\bR)}.}
\end{equation}

Moreover $(\ref{dfufdfdiuhuhuidsfhufdshudsfhdfsuhdsfifsdhdsf})$  remains true when $f$ belongs to $L^1_{loc}(\bR) \cap L^2_H(\bR)$ and $g$ belongs to $\sS(\bR)$ (in this case $(\ref{dfufdfdiuhuhuidsfhufdshudsfhdfsuhdsfifsdhdsf})$  reads  $<f,M_H(g)> = {<M_H(f),g>}_{L^2(\bR)}$, where $<,>$ denotes the duality bracket between $\sS'(\bR)$ and $\sS(\bR)$).

\item  There exists a constant $D$ such that,  for every couple $(H,k)$ in $(0,1)\times \bN^*$,

$\underset{x \in \bR }{\max}{\h1 |M_H(e_k)(x)|} \leq \tfrac{D}{c_H} \h1 {(k+1)}^{2/3}$.

\een

\end{theo}

\begin{pr}
\textit{1.} Since $M_H$ is an isometry, we just have to establish the surjectivity of $M_H$, for all $H$ in $(0,1)$. The case $H = 1/2$ being obvious, let us fix $H$ in  $(0,1)\backslash \{1/2\}$, $g$ in $L^2(\bR)$ and define the complex-valued function $w^g_H$ on $\bR$ by $w^g_H(\xi) = \tfrac{c_H}{\sqrt{2\pi}}{|\xi|}^{H-1/2} \hat{g}(\xi)$ if $\xi$ belongs to $\bR^*$ and $w^g_H(0) := 0$. Define the tempered distribution $v^g_H$ by $v^g_H := \frac{1}{2\pi} \check{\widehat{w^g_H}}$, where for all tempered distribution $T$, by definition, $<\check{T},f> := <T,\check{f}>$ for all  functions $f$ in $\sS(\bR)$ and where  $\check{f}(x) = f(-x)$ for all $x$.
We shall prove that $v^g_H$ belongs to $L^2_H(\bR)$ and that $M_H(v^g_H) =g$. Note first that for all $u$ in  ${\sS}'(\bR)$, $\widehat{\widehat{u}} = 2\pi \check{u}$. It is clear that $\widehat{v^g_H} =w^g_H $ and  that $\widehat{v^g_H}$ belongs to $L^1_{\text{loc}}(\bR)$. Moreover, thanks to formula (\ref{alban}), we see that

\begin{equation*}
{||v^g_H||}^2_{H} =  \tfrac{1}{{c_H}^2}\int_{\bR} {|\xi|}^{1-2H} {|\widehat{v^g_H}(\xi)|}^2 d\xi = \tfrac{1}{2\pi} \int_{\bR}  {|\widehat{g}(\xi)|}^2 d\xi =  {{||g||}^2_{L^2(\bR)}} < + \infty .
\end{equation*}

This shows that $v^g_H$ belongs to $L^2_{H}(\bR)$. We can then compute $\widehat{M_H(v^g_H)}$ and obtain, for almost every $\xi$ in $\bR$,

\begin{equation*}
\widehat{M_H(v^g_H)} (\xi) =   \tfrac{\sqrt{2\pi}}{c_H} {|\xi|}^{1/2-H}  \widehat{v^g_H} (\xi)  =   \tfrac{\sqrt{2\pi}}{c_H} {|\xi|}^{1/2-H} w^g_H(\xi)  =  \widehat{g}(\xi).
\end{equation*}

The previous equality  shows that $M_H(v^g_H)$ is equal to $g$ in $L^2(\bR)$ and then establish the surjectivity of $M_H$.

\textit{2.} See equality (3.12) of \cite{bosw}. The case where $f$ belongs to $L^1_{loc}(\bR) \cap L^2_H(\bR)$ is obvious, in view of (\ref{dfufdfdiuhuhuidsfhufdshudsfhdfsuhdsfifsdhdsf}), using the density of $\sS(\bR)$ in  $L^2_H(\bR)$.

\textit{3.} is shown in lemma $4.1$ of  \cite{ell}.
\end{pr}

Of course if we just consider functions in $L^2_H(\bR)$ instead of all elements of $L^2_H(\bR)$, the map $M_H$ is not bijective any more.

\subsection{Study of $\frac{\partial M_H}{\partial H}$}

\makeatletter
\renewcommand\theequation{\thesection.\arabic{equation}}
\@addtoreset{equation}{section}
\makeatother

We now study the operator  $\frac{\partial M_H}{\partial H}$. It will prove instrumental in defining the integral with respect to mBm in section \ref{sjjjj}.

Heuristically, we wish to differentiate with respect to $H$ the expression in definition \eqref{qxs}, {\it i.e.} differentiate the map
$H\mapsto \widehat{M_H(u)}(y)$ on $(0,1)$ for $(u,y)$ in $L^2_H(\bR) \times \bR^*$, assuming this is possible.
By doing so, we define a new operator, denoted $\frac{\partial M_H}{\partial H}$, from a certain subset of $L^2_H(\bR)$ to $L^2(\bR)$.
Of course, in order to compute the derivative at $H_0$ of $H\mapsto \widehat{M_H(u)}(y)$, we need to consider a neighbourhood
$V_{H_0}$ of $H_0$ in $(0,1)$ and thus consider only elements $u$ which belong to $\underset{ \scriptstyle{ H \in V_{H_0}} }{\bigcap}{L^2_{H}(\bR)}$. However, as will become apparent, the formula giving the derivative makes sense without this restriction.

In order to define in a rigorous manner the operator $\frac{\partial M_H}{\partial H}$, we shall proceed in a way analogous to the one that allowed to define $M_H$ in the previous subsection. It will be shown in remark \ref{msddpsfvokpvkpofpokdopfkvopd} that this construction effectively defines
the derivative, in a certain sense, of the operator ${M_H}$.

We will note $c^{'}_H$ the derivative of the analytic map $H \mapsto c_H$ where $c_H$ has been defined in \eqref{covmbm} and set $\beta_H := \frac{c^{'}_H}{c_H}$. Let $H$ belong to  $(0,1)$. Define:

\begin{equation}
\label{nifvfoijnth}
\Gamma_H(\bR) = \{u \in \mathcal{\sS'(\bR)} \hspace{0.1cm} : \hspace{0.1cm} \widehat{u} = T_{f};  \h1 \h1 f \in L^{1}_{loc}(\bR) \hspace{0.1cm}\text{  and}\hspace{0.1cm}  {|| u ||}_{\delta_H(\bR)} <  +\infty  \},
\end{equation}

where the norm ${||\h1||}_{\delta_H(\bR)}$ derives from the inner product on $\Gamma_H(\bR)$ defined by

\begin{equation}
\label{albaojn}
{<u,v>}_{\delta_H} := \frac{1}{c^2_H} \int_{\bR} {(\beta_H + \ln|\xi|)}^2 \h1 {|\xi|}^{1-2H}   \h1 {\widehat{u \hspace{0.1cm}}(\xi) } \h1 \overline{{\widehat{v \hspace{0.1cm}}(\xi) }}  \h1 d{\xi}.
\end{equation}

By slightly adapting lemma \ref{bcv}, it is easy to check that $(\Gamma_H(\bR), {<,>}_{\delta_H(\bR)})$ is a Hilbert space which verifies the equality $\Gamma_H(\bR) = \overline{\sS(\bR)}^{{<,>}_{\delta_H} } = \overline{\cE(\bR)}^{{<,>}_{\delta_H} }$. Note furthermore that, for every $H$ in $(0,1)$, the inclusion $\Gamma_H(\bR) \subset L^2_H(\bR)$ holds. We may now define the operator $\tfrac{\partial M_H}{\partial H }$ from  $(\Gamma_H(\bR), {<,>}_{\delta_H(\bR)})$ to  $(L^2(\bR), {<,>}_{L^2(\bR)})$, in the Fourier domain, by:

\begin{equation}
\label{oldsfijodirjdfoijgfodijdfiojdfgpoijdf}
\widehat{\tfrac{\partial M_H}{\partial H }(u)}(y) :=  -(\beta_H + \ln|y|) \h1 \tfrac{\sqrt{2\pi}}{c_H}  \h1 {|y|}^{1/2 - H}  \h1 \widehat{u}(y),  \hspace{0.25cm} \text{for every }\h1 y \h1 \text{in} \h1 \bR^*.
\end{equation}

In particular, one can check that, for $f$ in $\sS(\bR)$,  $\widehat{\tfrac{\partial M_H}{\partial H }(f)}(y) = \tfrac{\partial}{\partial H} \widehat{M_H(f)(y)}$ for almost every real $y$.
 Since $\widehat{\tfrac{\partial M_H}{\partial H }(u)}$ belongs to $L^2(\bR)$ for every $u$ in $\Gamma_H(\bR)$, $\tfrac{\partial M_H}{\partial H }$
 is well defined and given by its inverse Fourier transform from $(\Gamma_H(\bR), {<,>}_{\delta_H(\bR)})$ to  $(L^2(\bR), {<,>}_{L^2(\bR)})$:

\begin{equation*}
\tfrac{\partial M_H}{\partial H }(u)(x)= \frac{1}{2\pi}\cF( \widehat{\tfrac{\partial M_H(u)}{\partial H }})(-x),  \hspace{0.5cm}  \text{for almost every}\h1 x \h1 \text{in} \h1 \bR.
\end{equation*}

As in the previous subsection it will be useful to compute $\tfrac{\partial M_H}{\partial H }(f)$ for $f$ in $\sS(\bR)$. We summarize, in following proposition, the main properties of $\tfrac{\partial M_H}{\partial H }$.

\begin{prop}
\label{flmdkkmlkdslkdmlkdsmlksdmlksdmlksdml}
$\tfrac{\partial M_H}{\partial H }$ is an isometry from $(\Gamma_H(\bR), {<,>}_{\delta_H(\bR)})$ to  $(L^2(\bR), {<,>}_{L^2(\bR)})$ which verifies:

\begin{align}
\label{jeveuxx}
&\forall f \h1 \in \h1 \Gamma_H(\bR), &  {||f||}_{\delta_H} &=  {|| \tfrac{\partial M_H}{\partial H }(f) ||}_{L^2(\bR)},\\
\label{jeveux2}
&\forall (f,g) \h1 \in \h1 {(\Gamma_H(\bR) \cap L^2(\bR))}^2, & {<\tfrac{\partial M_H}{\partial H }(f),g>}_{L^2(\bR)} &= {<f,\tfrac{\partial M_H}{\partial H }(g)>}_{L^2(\bR)}, \\
\label{jeveux3}
&\forall f \h1 \in \h1 \sS(\bR) \cup \cE(\bR), \h1 \text{ and for }\textit{ a.e. } x \in \bR, & \tfrac{\partial M_H}{\partial H }(f)(x) &= \tfrac{\partial}{\partial H} [M_H(f)(x)].
\end{align}

%

\end{prop}

\renewcommand{\thefootnote}{\arabic{footnote}}

\begin{pr}

Equality \eqref{jeveuxx} results immediately from the definition of $\tfrac{\partial M_H}{\partial H }$ and from \eqref{oldsfijodirjdfoijgfodijdfiojdfgpoijdf}. For any couple of  functions $(f,g)$ in ${(\Gamma_H(\bR) \cap L^2(\bR))}^2$,
$$
\begin{array}{lll}
{<\tfrac{\partial {M_H} }{\partial H }(f),g>}_{L^2(\bR)} &= \tfrac{1}{2\pi} {<\widehat{\tfrac{\partial M_H }{\partial H }(f)},\widehat{g}>}_{L^2(\bR)} \notag \\
&= \tfrac{1}{2\pi} \int_{\bR} -(\beta_H + \ln|y|) \h1 \tfrac{\sqrt{2\pi}}{c_H}  \h1 {|y|}^{1/2 - H}  \h1 \widehat{f}(y) \overline{\widehat{g}(y)} dy \notag \\
& = \tfrac{1}{2\pi} {<\widehat{f},\widehat{\tfrac{\partial {M_H} }{\partial H }(g)}>}_{L^2(\bR)} =  {<f,\tfrac{\partial {M_H} }{\partial H }(g)>}_{L^2(\bR)}.
\end{array}
$$


It just remains to prove \eqref{jeveux3}. Since we will not use \eqref{jeveux3} in the sequel for $f$  in $\cE(\bR)$, we will just establish it here on $\sS(\bR)$.
Let $f$ be in $\sS(\bR)$ and $H$ in $(0,1)$. Formulas (\ref{referque}), (\ref{referquee}) and (\ref{referqueee}) can be summarized by
\vspace{-1mm}
\begin{equation*}
M_H(f)(x) = \gamma_H <{|y|}^{-(3/2-H)}, f(x+y)> \hspace{0.25cm}  \text{for almost every real } x,
\end{equation*}

where we have written, by abuse of notation, ${|y|}^{-(3/2-H)}$ for the  tempered distribution $y \mapsto {|y|}^{-(3/2-H)}$ and $f(x+y)$ for the map $y\mapsto f(x+y)$. Using theorem \ref{peogi} which is given in subsection \ref{ozieffjioefjoezifjioezzjioefjfeoiz} below, we may write

\begin{equation}\label{fouojtfedozdzezedzeijoiejfoijfezife}
\frac{\partial}{\partial H}[ M_H(f)(x)] = \gamma'_H <{|y|}^{-(3/2-H)}, f(x+y)> \hspace{0.25cm}  + \h1 \gamma_H <{|y|}^{-(3/2-H)} \ln |y|, f(x+y)>.
\end{equation}

Furthermore, for almost every real $x$, we may write:
\vspace{-1mm}
\begin{multline*}
\tfrac{\partial {M_H}}{\partial H }(f)(x) =  \frac{1}{2\pi}  \widehat{\widehat{ \tfrac{\partial {M_H}}{\partial H }(f)}} (-x) =  \frac{1}{2\pi} \cF(y\mapsto -(\beta_H + \ln|y|) \frac{\sqrt{2\pi}}{c_H} {|y|}^{(1/2-H)} \widehat{f}(y))(-x)\\
= - \beta_H M_H(f)(x) - \frac{1}{c_H \sqrt{2\pi}} \cF(y\mapsto {|y|}^{(1/2-H)}  \ln|y| \widehat{f}(y))(-x).
\end{multline*}

Define, for every real $x$,  $I(-x):=\cF(y\mapsto {|y|}^{(1/2-H)}  \ln|y| \widehat{f}(y))(-x)$ and, for every $H$ in $(0,1)$, $\nu_H:= 2 \h1 \Gamma'(3/2-H) \sin( \tfrac{\pi}{2}(1/2-H)) + \pi \Gamma(3/2-H) \cos(\tfrac{\pi}{2}(1/2-H)   $.
Using \cite{GeCh} (p.$173-174$) we get, after some computations,

\begin{equation*}
I(-x) = (-c_H \sqrt{2\pi}) \gamma_H <{|y|}^{-(3/2-H)} \ln |y|,f(x+y)> - \nu_H <{|y|}^{-(3/2-H)},f(x+y)>.
\end{equation*}

We finally have, for almost every real $x$,

\begin{multline*}
\tfrac{\partial {M_H}}{\partial H }(f)(x) = ( \frac{-\gamma_H \h1 c'_H}{c_H} + \frac{\nu_H}{c_H \sqrt{2\pi}} ) <{|y|}^{-(3/2-H)},f(x+y)> + \gamma_H <{|y|}^{-(3/2-H)} \ln |y|,f(x+y)>
\end{multline*}
which is nothing but  \eqref{fouojtfedozdzezedzeijoiejfoijfezife} since  $ \frac{-\gamma_H \h1 c'_H}{c_H} + \frac{\nu_H}{c_H \sqrt{2\pi}} = \gamma'_H$.
\end{pr}

\begin{rem}\label{msddpsfvokpvkpofpokdopfkvopd}
For all positive real $r$ and real $H$ in $(0,1)$ let $\Sigma_{H,r}(\bR)$  be the space defined by\\
$\Sigma_{H,r}(\bR)~:=~\Gamma_H(\bR) \bigcap  \left( \underset{ \scriptstyle{ H' \in [H-r,H+r]} }{\bigcap}{L^2_{H'}(\bR)}\right)$\footnote{ Note that we have the equality $\underset{ \scriptstyle{ H \in [a,b]} }{\bigcap}{L^2_{H}(\bR)} = L^2_{a}(\bR) \cap L^2_{b}(\bR)$ for every $[a,b]\subset(0,1)$.} and $\Sigma_H(\bR):= \underset{r \in (0,\min(H,1-H)) }{\bigcup}{\Sigma_{H,r}(\bR)} $. It is possible to show that, for all $f$ in $\Sigma_H(\bR)$, $\tfrac{\partial M_H}{\partial H }(f)(.)$ (resp. $\widehat{\tfrac{\partial M_H}{\partial H }(f)}(.)$ ) is the derivative, in the $L^2(\bR)$-sense, of $M_H(f)$ (resp. of $\widehat{M_H(f)}$ ). Note moreover that the inclusion $\cE(\bR)\cup\sS(\bR) \subset \underset{r \in \bR^*_+ }{\bigcap}{\Sigma_{H,r}(\bR)}$ holds.
\end{rem}

\section{Wiener integral with respect to mBm on $\bR$}
\label{WienerInt}

\subsection{Wiener integral with respect to fBm}\label{mlqmqlmqlqmlqqlmlkdeoidjziehzhziuh}

 Similarly to what is performed in \cite{ell} and \cite{bosw} (in these works this is done only for functions of  $L^2_H(\bR)$), it is now easy to define a Wiener integral with respect to  fractional Brownian motion. Indeed, for any element $g$ in $L^2_H(\bR)$, define $\cJ^H(g)$ as the random variable $<.,  M_H(g)>$. In other words, for all couples $(\omega,g)$ in $\Omega \hspace{-0.05cm} \times \hspace{-0.05cm} L^2_H(\bR)$:

\begin{equation}\label{sixth}
\cJ^H(g) (\omega) := <\omega,  M_H(g) >  = \int_{\bR} M_H(g)(s) \hspace{0.05cm} dB(s) (\omega),
\end{equation}

where the Brownian motion $B$ has been defined just below formula (\ref{lala}). We call the random variable $\cJ^H(g)$ the Wiener integral of $g$ with respect to fBm.
Once again, when $g$ is a tempered distribution which is not a function, $g(s)$ does not have a meaning for a fixed real $s$ and $\cJ^H(g)$ is just a notation for the centered Gaussian random variable $<.,  M_H(g) >$.

\begin{rem}
Note that, for $H$ in $(0,1)$, we are able to give a meaning to $\cJ^H(g)$ only for elements $g$ which belong to $L^2_H(\bR)$ and not anymore for all elements $g$ in $\sS'(\bR)$ as was the case for $\int_{\bR} g(s) dB(s)$ (see remark \ref{zpopdokokdkkdkdkdkdkdkdkdk}).
\end{rem}

%
\subsection{Wiener integral with respect to mBm}\label{Wiener}
We now consider a fractional Brownian field $\Lambda:={\big(\Lambda(t,H)\big)}_{(t,H) \in \bR \times (0,1)}$,  defined, for all
$(t,H)$ in $\bR \times (0,1)$ and all $\omega$ in $\Omega$, by $\Lambda(t,H)(\omega):= {B}^{(H)}(t,\omega):= {<\omega,  M_H(\i1_{[0,t]})>}$.
We wish to generalize the previous construction of the Wiener integral with respect to fBm to the case of mBm. This amounts to replacing
the constant $H$ by a continuous deterministic function $h$, ranging in $(0,1)$. More precisely, let $R_h$ denote the covariance function of  a normalized mBm with function $h$ (see definition \ref{azazazaz}). Define the bilinear form ${<,>}_h$ on $\cE(\bR)\times \cE(\bR)$ by ${<\i1_{[0,t]},\i1_{[0,s]} >}_{h} = R_{h}(t,s)$.  Our construction of the integral of deterministic elements with respect to mBm requires that ${<,>}_h$ be an inner product:

\begin{prop}\label{rfjiijofsdqoijqfsdoijqdsqsdoufgyq}
${<,>}_h$ is an inner product for every function $h$.
\end{prop}

\begin{pr}
See  appendix \ref{appenB}.
\end{pr}

Define the linear map $M_{h}$ by:

$$
\begin{defappli}{M_{h}}{ (\cE(\bR),{< ,\hss >}_{h} )}{ (L^2(\bR),{< ,\hss >}_{L^2(\bR)})}{\i1_{[0,t]}}
M_{h}(\i1_{[0,t]}) := M_{h(t)}(\i1_{[0,t]}) := {M_{H}(\i1_{[0,t]})|}_{H =h(t)} .
\end{defappli}
$$

Define the  process {$\widetilde{B}^{(h)}(t)~=~\hspace{-0.1cm}{<.,M_{h}(\i1_{[0,t]})>}$}, $t \in \bR$. As Kolmogorov's criterion and the proof of following lemma show, this process admits a continuous version which  will be noted ${B}^{(h)}$. A word on notation: we write ${B}^{(.)}$ both for an fBm and an mBm. This should
not cause any confusion since an fBm is just an mBm with constant $h$ function. It will be clear from the context in the following whether the "$h$" is constant or not. Note that:

\begin{equation}\label{fqkiuhgiuhihj}
a.s.,  \ \forall t \in \bR, \quad B^{(h)}(t) = {B^{(H)}(t)}_{|_{H = h(t)}}
\end{equation}

In view of point \textit{2.} in remark \ref{remme} we may state the following lemma.

\begin{lem}
\label{osidjoifjsoj}
\bit
\iti The process ${B}^{(h)}$ is a normalized mBm.
\itii The map $M_h$ is an isometry from $(\cE(\bR),{<,>}_{h})$ to $(L^2(\bR), {<,>}_{L^2(\bR)})$.
\eit
\end{lem}
\begin{prz}
The process ${B}^{(h)}$ is clearly a centered Gaussian process. Moreover, for all $(s,t)$ in $\bR^2$,
\begin{align*}
\bE[B^{(h)}(t)B^{(h)}(s)] &= \bE[ <.,M_{h}(\i1_{[0,t]})><.,M_{h}(\i1_{[0,s]})>] =  {<M_{h}(\i1_{[0,t]}),M_{h}(\i1_{[0,s]})>}_{L^2(\bR)} \\
 & = \frac{1}{2\pi} {<\widehat{M_{h(t)}(\i1_{[0,t]})},\widehat{M_{h(s)}(\i1_{[0,s]})}>}_{L^2(\bR)} = \tfrac{1}{c_{h(t)}c_{h(s)}}  {\scriptsize \int_{\bR}  } \frac{(1-e^{it\xi})(1-e^{-is\xi})}{{|\xi|}^{1+2h_{t,s}}}   d{\xi} \notag \\
 &=  \frac{c^2_{h_{t,s}}}{c_{h(t)}c_{h(s)}}    \frac{1}{2}\hspace{0.1cm} \big( {|t|}^{2h_{t,s}} + {|s|}^{2h_{t,s}}- {|t-s|}^{2h_{t,s}})
  =  R_{h}(t,s)={<\i1_{[0,t]},\i1_{[0,s]} >}_{h}. \qquad \square
\end{align*}
\end{prz}

By isometry, it is possible to extend $M_h$ to the space $\overline{\cE(\bR)}^{{<,>}_{h}}$ and we shall still note $M_h$ this extension. Define
the isometry $\cJ^{h} := \zeta \circ M_{h}$ on $\overline{\cE(\bR)}^{{<,>}_{h}}$, \textit{i.e}:

$$
\begin{defappliab}{ \cJ^{h}}{ (\overline{\cE(\bR)}^{{< >}_{h}}, {< , >}_{h})  }{ (L^2(\bR),{<,>}_{L^2(\bR)})}{ ((L^2), {< , >}_{(L^2)} )    }{ \i1_{[0,t]} }{  M_{h}(\i1_{[0,t]})}
{<.,M_{h}(\i1_{[0,t]})>.}
\end{defappliab}
$$

\vspace{0.1cm}

We can now define the Wiener integral with respect to mBm in the natural following way:

\begin{defi}
\label{wienerR}
Let $B^{(h)}$ be a normalized multifractional Brownian motion. We call Wiener integral on $\bR$ of an element $u$ in $\overline{\cE{(\bR)}}^{{< >}_{h}}$ with respect to $B^{(h)}$, the element ${\cJ}^{h}(u)$ of $(L^2)$ defined thanks to the isometry ${\cJ}^{h}$ given just above.
\end{defi}

\begin{rem}\label{qzfekdpokfdpof}
It follows from  definition  \ref{wienerR} that the  Wiener integral of a finite linear combination of functions $\i1_{[0,t]}$ is ${\cJ}^{h}( \sum^{n}_{k=1} \alpha_k  \i1_{[0,t_k]} )  ~=~\sum^{n}_{k=1} \alpha_k B^{(h)}_{t_k}  $. Moreover, for any element $u$ in $\overline{\cE(\bR)}^{{< >}_{h}}$ (which may be a tempered distribution), the Wiener integral of $u$ with respect to mBm, still denoted  ${\cJ}^{h}(u)$, is given by ${\cJ}^{h}(u)  \stackrel{\text{def}}{=} \lim\limits_{n \to +\infty} {\cJ}^{h}(u_n)$, for any sequence of functions ${(u_n)}_{n \in \bN}$ in ${\cE(\bR)}^{\bN}$ which converges to $u$ in the norm ${|| \hspace{0.1cm} ||}_{h}$ and where the convergence of ${\cJ}^{h}(u_n)$ holds in $(L^{2})$.
\end{rem}

Since we now have a construction of the Wiener integral with respect to  mBm, it is natural to ask which functions admit  such an integral. In particular, we do not know so far whether $\sS(\bR) \subset \overline{\cE(\bR)}^{{< >}_{h}}$. The next section contains more information about the space $\overline{\cE(\bR)}^{{< >}_{h}}$.

\section{Stochastic integral with respect to mBm}\label{sjjjj}

\subsection{Fractional White Noise}

The following theorem will allow us to give a concrete example of a derivative of an ${(\cS)}^*$-process.

\begin{theo}
\ben
\item For any real $H$ in $(0,1)$, the map $M_H(\i1_{[0,.]}): \bR \rightarrow {\sS'}(\bR)$ defined by $M_H(\i1_{[0,.]})(t):=M_H(\i1_{[0,t]})$ is differentiable over $\bR$ and its derivative, noted $\frac{d}{dt} [M_H(\i1_{[0,t]})]$, is equal to\\ $\sum^{+\infty}_{k = 1} M_H(e_k)(t) \h1 e_k$, where the convergence is in ${\sS'}(\bR)$.
\item  For any interval $I$ of $\bR$ and  any  differentiable map $F:I\rightarrow {\sS'}(\bR)$,  the element $<.,F(t)>$ is a differentiable stochastic distribution process which satisfies the equality \\ $\frac{d}{dt}<.,F(t)> =  <.,\frac{d}{dt}F(t)> $.
\een
\end{theo}

\begin{pr}
The proof of point $1$ is a mere re-writing of the one of lemma $2.15$ of \cite{ben1} by replacing $M^H_{\pm}$ by the operator $M_H$. Point $2$ is theorem $2.17$ of \cite{ben1}.
\end{pr}

Let $H\in (0,1)$. The process ${(B^{(H)}(t))}_{t \in \bR}$ defined in (\ref{fifth}) is an fBm, and $M_H(\i1_{[0,t]})$ belongs to $L^2(\bR)$ for every real $t$. Hence, using  equality (\ref{dfufdfdiuhuhuidsfhufdshudsfhdfsuhdsfifsdhdsf}),
we may write, for every real $t$ and almost surely:

\begin{multline}\label{mdem}
   \hspace{-0.2cm} B^{(H)}(t) = \h1 <.,M_H(\i1_{[0,t]})> \h1 = \h1 <.,\sum^{+\infty}_{k = 0} {<M_H(\i1_{[0,t]}),e_k>}_{L^2(\bR)} \h1 e_k >  \\
				 \hspace{-0.2cm} = \h1 \sum^{+\infty}_{k = 0} {<\i1_{[0,t]},M_H(e_k)>}_{L^2(\bR)} \h1 <.,e_k> \h1 = \h1 \sum^{+\infty}_{k = 0} \left(\int^t_0 M_H(e_k)(u) du\right) <.,e_k>.
\end{multline}

\eqref{mdem} and the previous theorem lead to the definition of fractional white noise \cite{ell,bosw}:

\begin{ex}[{\bfseries Fractional white noise }]
\label{oijoifjsoidjdoijfsoijfsoijdsoifjdoisdjosdjf}
Let:
\begin{equation}
\label{oerijefoijefjfijvoidfovijdfovijdfo}
W^{(H)}(t) := \displaystyle{ \sum^{+\infty}_{k = 0} \h1 M_H(e_k)(t) \h1  <.,e_k>}.
\end{equation}

 Then $({W^{(H)}(t))}_{t \in \bR}$ is a  ${(\cS)}^*$-process and is the ${(\cS)}^*$-derivative
of the process ${(B^{(H)}(t))}_{t \in \bR}$.

The proof of this fact is simple: for any integer $p\geq 2$, using remark \ref{laiusdhv}, the mean value theorem and the dominated convergence theorem,

\begin{align}\label{deezdezezezezm}
J_{p,r}(t) &:= {\left|\left|  \tfrac{B^{(H)}(t+r) - B^{(H)}(t)}{r}  - W^{(H)}(t) \right|\right|}^2_{-p} =  {\left|\left| \sum^{+\infty}_{k = 0} \h1     \left( \frac{1}{r} \int^{t+r}_t   M_H(e_k)(u) \h1 du -  M_H(e_k)(t)  \right)  <.,e_k >  \right|\right|}^2_{-p} \notag \\
 &= {\left| \sum^{+\infty}_{k = 0} \h1     \left( \frac{1}{r} \int^{t+r}_t   M_H(e_k)(u) \h1 du -  M_H(e_k)(t)  \right)  \h1 e_k \h1 \right|}^2_{-p}\notag \\
 &=  \sum^{+\infty}_{k = 0} \h1  {(2k+2)}^{-2p} \h1 { \left( \frac{1}{r} \int^{t+r}_t   (M_H(e_k)(u) -  M_H(e_k)(t) ) \h1 du  \right)}^2   \underset{r \to 0 }{\longrightarrow}{ 0 }.
\end{align}
\end{ex}

\begin{rem}
\label{iooisoidjoissid}
In particular we see that for all $(t,H)$ in $\bR\times(0,1)$, $W^{(H)}(t)$ belongs to $({\cS_{-p}})$ as soon as $p \geq 2$.
\end{rem}

\begin{rem}
There are several constructions of fBm. In particular, operators different from $M_H$ may be considered. $\cite{ben1}$ uses an operator denoted $M^H_+$ on the grounds that fBm as defined here is not adapted to the filtration generated by the driving Brownian motion as soon as $H\neq 1/2$. While this is indeed a drawback, the crucial property for our purpose is that the same probability space $(\sS'(\bR),\cG,\mu)$ is used for all parameters $H$ in $(0,1)$. This allows to consider simultaneously several fractional Brownian motions with $H$ taking any value in $(0,1)$, which is necessary when one deals with mBm. We choose here to work with $M_H$ rather than with
$M^H_+$  and $M^H_-$ of $\cite{ben1}$ as its use is simpler. $M^H_+$  and $M^H_-$ would nevertheless allow for a more general approach
encompassing the whole family of mBm at once. This topic will be treated in a forthcoming paper.
\end{rem}

\subsection{Multifractional White Noise}
	
The main idea for defining a stochastic integral with respect to  mBm is similar to the one used for fBm. We will relate the process $B^{(h)}$ to Brownian motion \textit{via} the family of operators ${(M_H)}_{(H \in (0,1))}$. This will allow to define a multifractional white noise, analogous to the fractional white noise of example \ref{oijoifjsoidjdoijfsoijfsoijdsoifjdoisdjosdjf}. From a heuristic point of view, multifractional white noise
is obtained by differentiating with respect to $t$ the fractional Brownian field $\Lambda(t,H)$ (defined at the beginning of section \ref{Wiener}) along a curve $(t,h(t))$.
Assuming that we may differentiate in the sense of ${(\cS)}^*$ (this will be justified below), the differential of $\Lambda$ reads:

\begin{align}\label{covfBmriujotrjogtpoekffjrjtgrjo}
d\Lambda(t,H) &= \frac{\partial \Lambda}{\partial t }(t,H) \h1 dt + \frac{\partial \Lambda}{\partial H }(t,H) \h1dH
= \frac{dB^{(H)}}{dt }(t) \h1 dt + \frac{\partial \Lambda}{\partial H }(t,H) \h1 dH \notag \\
&= W^{(H)}(t) dt +  \frac{\partial \Lambda}{\partial H }(t,H) \h1 dH,
\end{align}

where the equality will be shown to hold in  ${(\cS)}^*$. With a differentiable function $h$ in place of $H$, this formally yields

\begin{equation}\label{frejo}
d\Lambda(t,h(t))  =  \big( W^{(h(t))}(t)  +  h'(t) \h1 {\frac{\partial \Lambda}{\partial H }(t,H) }{\big|}_{{H = h(t)}}\big) \h1 dt.
\end{equation}

In view of the definition of the stochastic integral with respect to fBm given in \cite{ell}, \cite{HOUZ} and  \cite{ben2}, it then
seems natural to set the following definition for the stochastic integral with respect to  mBm of a Hida process $X:\bR\rightarrow {(\cS)}^*$:

\begin{align}\label{frejoj}
\int_{\bR} X(s) dB^{(h)}(s) &:= \int_{\bR} X(s) d\Lambda(s,h(s)) \notag\\
&:=  \int_{\bR} X(s) \diamond \big( W^{(h(s))}(s)  +  h'(s) \frac{\partial \Lambda}{\partial H }(s,H) |_{H = h(s)}\big) ds.
\end{align}

We shall then say that the process $X$ is integrable with respect to  mBm if the right hand side of (\ref{frejoj}) exists in  ${(\cS)}^*$. Remark that when the function $h$ is constant we recover of course the integral with respect to fBm.

In order to make the above ideas rigorous, we start by writing the chaos expansion of $B^{(h)}$. Since
 $M_H(g)$ belongs to $L^2(\bR)$ for all $(g,H)$ in $\sS(\bR)\times (0,1)$, we may  define, for all $H$ in $(0,1)$, $M_H:\sS'(\bR) \rightarrow \sS'(\bR)$, by

\begin{equation}\label{firsozeijdt}
<M_H(\omega),g> \h1 = \h1 <\omega,M_H(g)>, \hspace{0.25cm} \text{for}\hspace{0.1cm}\mu-a.e. \hspace{0.1cm} \omega \hspace{0.25cm} \text{in} \hspace{0.25cm}\Omega = \sS'(\bR).
\end{equation}

Moreover, in view of remark \ref{remme}, we may extend (\ref{firsozeijdt}) to the case where $g$ belongs to $L^2_H(\bR)$ by writing, for all $g$ in  $L^2_H(\bR)$ and almost every $\omega$ in $\Omega$,

\begin{equation}\label{firsozeijdtopopkpo}
<M_H(\omega),g> := \lim_{n \to +\infty} <M_H(\omega),g_n> \h1=\h1 \lim_{n \to +\infty} <\omega,M_H(g_n)> \h1=\h1 <\omega,M_H(g)>,
\end{equation}

for every sequence ${(g_n)}_{n \in \bN}$ of functions of $\sS(\bR)$ which converges to $g$ in the norm ${||\hspace{0.1cm}||}_{L^2_H(\bR)}$.
For all real $t$  and integer $k$ in $\bN$,  define the element of $L^2_{h(t)}(\bR)$:

\begin{equation}
\label{firsozeedeijdt}
d^{(t)}_{k} := M^{-1}_{h(t)}(e_k).
\end{equation}

It is clear that, for all $t$ in $\bR$, the family of functions ${(d^{(t)}_{k})}_{k \in \bN}$ forms an orthonormal basis of $L^2_{h(t)}(\bR)$. Let us now write the chaos decomposition of mBm. For almost every $\omega$ and every  real $t$ we get, using (\ref{dfufdfdiuhuhuidsfhufdshudsfhdfsuhdsfifsdhdsf}) and (\ref{firsozeijdtopopkpo}),

\begin{align*}
 B^{(h)}(t)(\omega) 	&= <\omega, M_{h}(\i1_{[0,t]}) > = <\omega, M_{h(t)}(\i1_{[0t]}) > = <M_{h(t)}(\omega), \i1_{[0,t]} > \notag \\
						&=  <M_{h(t)}(\omega), \displaystyle{\sum^{+ \infty}_{k = 0} {<\i1_{[0,t]}, d^{(t)}_{k} >}_{L^2_H(\bR)} d^{(t)}_{k}}> =  <M_{h(t)}(\omega), \displaystyle{\sum^{+ \infty}_{k = 0} {<M_{h(t)}(\i1_{[0,t]}), {e}_k >}_{L^2(\bR)} d^{(t)}_{k}}>\notag \\
						&= \sum^{+ \infty}_{k = 0} {<M_{h(t)}(\i1_{[0,t]}), {e}_k >}_{L^2(\bR)}
<M_{h(t)}(\omega),d^{(t)}_{k}>\\
						&= \sum^{+ \infty}_{k = 0} {<M_{h(t)}(\i1_{[0,t]}), {e}_k >}_{L^2(\bR)} <\omega,M_{h(t)}(d^{(t)}_{k})>.
\end{align*}

We get finally:

\begin{equation}\label{fiizeufhirsozeedeijdt}
a.s, \h1 \forall t \in \bR, \hspace{0.5cm}   B^{(h)}(t) =  { \sum^{+ \infty}_{k = 0} \left( \int^t_0  M_{h(t)}({e}_k)(s) ds \right) <.,e_k>}.
\end{equation}

We would then like to define multifractional white noise as the $(\cS^*)$-derivative of $B^{(h)}$, which would be formally defined by:

\begin{equation}\label{fiizeufhirsozeedeijdefrt}
W^{(h)}(t) := \displaystyle{\sum^{+ \infty}_{k = 0} \bigg[\frac{d}{dt}} \bigg( \int^t_0  \h1 M_{h(t)}({e}_k)(s)\h1 ds \bigg) \bigg]  \h1 <.,e_k>,
\end{equation}

assuming $h$ is differentiable. The following theorem states that, for all real $t$,  the right hand side of (\ref{fiizeufhirsozeedeijdefrt}) does
indeed belong to ${(\cS)}^*$ and  is  exactly the ${(\cS)}^*$-derivative  of $B^{(h)}$ at $t$.

\begin{theodef}\label{alalalalalal}
Let $h:\bR \rightarrow (0,1)$ be a  $\cC^1$ function such that the derivative function $h'$ is bounded. The process $W^{(h)} := {(W^{(h)}(t))}_{t \in \bR }$ defined by formula $(\ref{fiizeufhirsozeedeijdefrt})$ is an ${(\cS)}^*$-process which verifies the following equality in ${(\cS)}^*$:

\begin{equation}\label{firsozoijkqsjdefrt}
 W^{(h)}(t) =  \displaystyle{  \sum^{+ \infty}_{k = 0} {  M_{h(t)}({e}_k)(t) } <.,e_k> + h'(t)\h1 \sum^{+ \infty}_{k = 0} \bigg( \int^t_0
\tfrac{\partial M_{H}}{\partial H}({e}_k) (s) \big|_{H=h(t)}  ds \bigg) <.,e_k>  }.
\end{equation}

Moreover the process $B^{(h)}$ is ${(\cS)}^*$-differentiable on $\bR$ and verifies in ${(\cS)}^*$

\begin{equation}\label{firsozkqsjdefrt}
\frac{dB^{(h)}}{dt}(t) = W^{(h)}(t) = \frac{d}{dt}[\Lambda (t,h(t))].
\end{equation}
\end{theodef}

In order to prove this theorem, we will need two lemmas.

\begin{lem}
\label{zekkpoerkpeorperokvpo}
For $H$ in $(0,1)$ and  $f$ in  $\sS(\bR)$,  define $g_f: \bR \times (0,1) \rightarrow \bR$ by $g_f(t,H):= \displaystyle{\int^t_0 M_H(f)(x) dx }$. Then 
\vspace{-0.25cm}
\bit
\iti The function $g_f \hspace{0.1cm} \text{ belongs to} \hspace{0.1cm}  \cC^{\infty}(\bR \times (0,1), \bR)$,
\itii $\forall x \in \bR$, \h1  $M_H(f)(x) = \alpha_H \displaystyle{\int^{+\infty}_0 u^{H-1/2} (f'(x+u) - f'(x-u)) \h1 du}$.
\eit	
where $\alpha_H$ has been defined by $(\ref{alpha_H})$.

In particular, the function $(x,H) \mapsto M_H(f)(x)$ is differentiable on $\bR \times (0,1)$.
\bit
\itiii Assume that $h: \bR \rightarrow (0,1)$ is differentiable.  Then, for any real $t_0$
\eit

\begin{equation}\label{firsodezr'kqsjdefrt}
\frac{d}{dt}[g_f(t,h(t))] \big|_{t = t_0} = M_{h(t_0)}(f)(t_0) + h'(t_0)  \displaystyle{\int^{t_0}_0  \tfrac{\partial M_H}{\partial H}(f) (s) \big|_{ H = h(t_0)} ds }.
\end{equation}
\end{lem}

\begin{pr}
\vspace{-0.25cm}
\bit
\iti
Define $\mu_f$ on  $\bR \times (0,1)$ by ${\mu_f(t,H):=\int_{\bR} (u-t){|u-t|}^{H-3/2} f(u) \h1 du }$, for $f$ in $\sS(\bR)$.
\eit
Using (\ref{referqueosijz}) we get, for all $(t,H)$ in $\bR \times (0,1)$, the equality
\begin{equation}\label{fifroizerfiuheriuehriuerhfiuerhhigurehguierhgieurht}
g_f(t,H) = \alpha_H [{\mu}_f(t,H) - {\mu}_f(0,H)].
\end{equation}

A change of variables yields
\begin{align}\label{fifrt}
\mu_f(x,H) 	&= -\int^{x}_{-\infty} {|t-x|}^{H-1/2} f(t) dt + \int^{+\infty}_{x} {|t-x|}^{H-1/2} f(t) \h1 dt \notag \\
			&= \int^{+\infty}_{0} {u}^{H-1/2} (f(x+u) - f(x-u)) \h1 du.
\end{align}

Thanks to  (\ref{fifroizerfiuheriuehriuerhfiuerhhigurehguierhgieurht}) and to the fact that the map $y \mapsto \alpha_y$ is $\cC^{\infty}$  on $(0,1)$, it is sufficient to show that the function $\mu_f$ belongs to $\cC^{\infty}(\bR \times (0,1))$. In view of applying the theorem of differentiation under the integral sign, define $j(x,H,u):= {u}^{H-1/2} (f(x+u) - f(x-u))$ for  $u$ in $\bR^*_+$. Let $n$ in $\bN$ and $(\alpha_1,\alpha_2)$ in ${\bN}^2$ such that $\alpha_1+\alpha_2 = n$.
For almost every  $u$ in $\bR^*_+$, $(x,H) \mapsto j(x,H,u)$ is  $\cC^n$ on  $\bR \times (0,1)$ with partial derivatives given by

\vspace{-0.7cm}
\begin{align*}
 \frac{{\partial}^{n}j}{\partial x^{\alpha_1} \partial H^{\alpha_2}}(x,H,u) = {(\ln u)}^{\alpha_2} u^{H-1/2} (f^{(\alpha_1)}(x+u) - f^{(\alpha_1)}(x-u)).
\end{align*}
\vspace{-0.5cm}

Fix $(x_0,H_0)$ in $\bR\times(0,1)$. Let us  show that $\mu_f$ is $\cC^n$ in a neighbourhood of $(x_0,H_0)$. Choose $(a,b)$ such that $a< x_0 <b$ and $H_1,H_2$ such that $0<H_1<H<H_2<1$. We have

\begin{multline}\label{fifriuhiht}
\big|\frac{{\partial}^{n}j}{\partial x^{\alpha_1} \partial H^{\alpha_2}}(x,H,u)\big| \leq  {|u|}^{H_1-1/2} {|\ln u|}^{\alpha_2}   \h1 \i1_{ \{ 0 < u < 1\}} \h1 \underset{(x,u) \in [a,b]\times[0,1] }{\sup}{|f^{(\alpha_1)}(x\pm u)|}\\
+ {|u|} {|\ln u|}^{\alpha_2} {|f^{(\alpha_1)}(x\pm u)|} \h1 \i1_{\{1 \leq  |u| \}},	
\end{multline}

where $f^{(\alpha_1)}(x\pm u):= f^{(\alpha_1)}(x+u) - f^{(\alpha_1)}(x-u)$. A Taylor expansion shows that there exists a real $D$ such that, for all  $(u,x)$ in $\bR\times(a,b)$, ${|u|}^4 \h1 {|f^{(\alpha_1)}(u\pm x)|}  \leq D$. As a consequence, there exists a real constant $C$ such that, for almost every $u$ in $\bR^*_+$ and every  $(x,H)$ in $[a,b] \times [H_1,H_2]$,

\begin{equation}\label{fifriuhihteddeded}
\big|\tfrac{{\partial}^{n}j}{\partial x^{\alpha_1} \partial H^{\alpha_2}}(x,H,u)\big| \leq C \h1 \big[ {|u|}^{H_1-1/2} {|\ln u|}^{\alpha_2} \h1 \i1_{\{ 0<u <1 \}}  +   {|\ln u|}^{\alpha_2} \h1 \tfrac{1}{{|u|}^3} \h1 \i1_{\{  1  \leq u \} } \big]. 	
\end{equation}

Since the right hand side of the  previous inequality belongs to $L^1(\bR)$, the theorem of differentiation under the integral sign can  be applied  to conclude that the function $\mu_f$ is of class $\cC^n$ in $[a,b]\times [H_1,H_2]$, for all integer $n$ and all $f$ in $\sS(\bR)$. This entails ${\it(i)}$.

${\it(ii)}$ (\ref{fifroizerfiuheriuehriuerhfiuerhhigurehguierhgieurht}) and (\ref{fifrt}) yield

\begin{equation}\label{pdpdopkvpervkpvkfdkvpdfchezroberta}
M_H(f)(x) = \alpha_H \frac{\partial }{\partial x}[\mu_f(x,H)] = \alpha_H \int^{+\infty}_{0} {u}^{H-1/2} (f'(x+u) - f'(x-u)) du,  
\end{equation}

which establishes  {\bfseries{$(ii)$}} and the fact that  $(x,H)\mapsto M_H(f)(x)$ belongs to $\cC^{\infty}(\bR \times (0,1))$.

${\it(iii)}$ For a differentiable function $h$, we have, for every real $t_0$,

\begin{equation*}
\frac{d}{dt}[g_f(t,h(t))] \big|_{t = t_0} = \frac{\partial g_f}{\partial t} (t_0, h(t_0)) + h'(t_0) \frac{\partial g_f}{\partial H} (t_0, h(t_0)).
\end{equation*}

(\ref{fifroizerfiuheriuehriuerhfiuerhhigurehguierhgieurht}), \eqref{pdpdopkvpervkpvkfdkvpdfchezroberta} and \eqref{jeveux3}  show that, for every $f$ in $\sS(\bR)$,

\begin{align*}
\frac{\partial g_f}{\partial H} (t,H) = \int^t_0  \tfrac{\partial  }{\partial H} \h1   [M_H(f)(x)]       \h1 dx = \int^t_0 \tfrac{\partial M_H} {\partial H}(f) (x) \h1 dx.
\end{align*}

and finally:
\begin{equation*}
\frac{d}{dt}[g_f(t,h(t))] \big|_{t = t_0} = M_{h(t_0)}(f)(t_0) + h'(t_0)  \displaystyle{\int^{t_0}_0  \tfrac{\partial M_H}{\partial H}(f) (s) \big|_{ H = h(t_0)} ds }.
\end{equation*}
\end{pr}

\begin{lem}\label{zoidpsqojxiqsoxijqsojdzj}
The following inequalities hold:

\ben
\iti $\forall \h1 [a,b] \subset (0,1),  \exists \rho \in \bR$:  \h1 $\forall k  \in \bN, \hspace{0.1cm} { \underset{(H,u) \in [a,b]\times \bR }{\sup}{ \hspace{0.1cm} \left|\tfrac{\partial M_{H}}{\partial H}(e_k)(u)\right| } } \leq \rho \h1 {(k+1)}^{2/3} \h1 \ln(k+1)$.
\itii $\forall t \in \bR$, $\forall r \in \bR^*_+$, $\exists \widetilde{D_t}(r) \in \bR$,  $\forall k  \in \bN: \h1 \underset{u \in  [t-r,t+r] }{\sup}{ \h1 \left|\tfrac{d}{du} [g_{e_k}(u,h(u))] \right| } \leq \widetilde{D_t}(r) \h1 {(k+1)}^{2/3}$.
\een

\end{lem}

\begin{pr}
$(i)$
Since $\sS(\bR)$ is a subset of $\Gamma_H(\bR)$, \eqref{oldsfijodirjdfoijgfodijdfiojdfgpoijdf} entails that $\widehat{\frac{\partial M_H}{\partial H}(e_k)}$ belongs to $L^1(\bR) \cap L^2(\bR)$ for every $k$ in $\bN$. Furthermore $ \widehat{e_k}(y) = {(-i)}^{k-1} \sqrt{2\pi} e_k(y) $ for every integer $k$ in $\bN^*$ and for almost every real $y$ (see lemma $1.1.3$ p.$5$ of \cite{Tha}). Thus, for every $H \in [a,b]$ and almost every $u \in \bR$,

\begin{align}\label{zuiddhciuhcidsdhcisduhcdifuhck}
\tfrac{\partial  M_H(e_k)}{\partial H}  (u) &= \tfrac{1}{2\pi} \widehat{  \widehat{  \tfrac{\partial M_H(e_k)}{\partial H} }  }(-u) =  -\tfrac{1}{2\pi}  \int_{\bR} e^{iuy} (\beta_H + \ln|y|) \h1 \tfrac{\sqrt{2\pi}}{c_H}  \h1 {|y|}^{1/2 - H}  \h1 \widehat{e_k}(y)  \h1   dy \notag \\
&= \tfrac{{-\beta}_H}{2\pi} \int_{\bR} e^{iuy}     \widehat{M_H(e_k)}(y) \h1 dy \h1- \h1 \tfrac{1}{c_H}\int_{\bR} e^{iuy} {|y|}^{1/2-H}(\ln|y|) \h1 {(-i)}^{k-1} e_k(y) \h1 dy             \notag \\
&= -{\beta}_H \h1 M_H(e_k)(u) + \tfrac{-{(-i)}^{k-1}}{c_H} \underbrace{ \int_{\bR} e^{iuy} {|y|}^{1/2-H}(\ln|y|) \h1 e_k(y) \h1 dy}_{=:V_k(u)}.
\end{align}

Then, using (\ref{fkorep}), we see that there exists a family of real constants denoted ${(\rho_i)}_{1\leq i \leq 11}$ such that we have, for every couple $(H,k)$ in $[a,b]\times\bN$ and almost every real $u$,

\begin{equation}\label{oeirjoerjforejerjjerojvrorjortjo}
\begin{array}{lll}
|V_k(u)|	&\leq \rho_1 \left[\int_{|y| \leq 2\sqrt{k+1}} {|y|}^{1/2-H} |\ln|y|| |e_k(y)|dy + \int_{|y|> 2\sqrt{k+1}} {|y|}^{1/2-H}(\ln|y|)|e_k(y)|dy \right] \\
			&\leq \rho_2 \bigg[ \frac{1}{{(k+1)}^{1/12}} \big(  \underbrace{\int^{2\sqrt{k+1}}_0  |\ln y| \h1 {y}^{1/2-H} dy }_{:=I^{(k+1)}_1} \big) +  \big( \underbrace{\int^{+\infty}_{2\sqrt{k+1}}  {y}^{1/2-H}(\ln y) e^{-\gamma y^2}  dy  }_{:=I^{(k+1)}_2} \big) \bigg].
\end{array}
\end{equation}

An integration by parts yields

\begin{equation}\label{oeirjoerjfore"'rojvrorjortjo}
I^{(k)}_1 =  \rho_3(( \rho_4 + \ln k) (1 + {k}^{3/4 - H/2}) )   \leq \rho_5 \h1 {k}^{3/4 - H/2}  \ln k.
\end{equation}

Using the change of variables $u = y\sqrt{\gamma}$, we get $I^{(k)}_2 \leq  \rho_6 \h1 {\int^{+\infty}_{2\sqrt{k \gamma}}  {|u|}^{1/2-H}(\ln u) e^{-u^2} du	 } =  \rho_6 \h1 J^{(1/2 - H)}_{2\sqrt{k \gamma}} $ where $J^{(\alpha)}_{\delta}~:=~\int^{+\infty}_{\delta}  {|u|}^{\alpha} \ln u \h1 e^{-u^2} du$.
When $\delta > 3e$, an integration by parts shows that $J^{(\alpha)}_{\delta} <  \delta^{\alpha-1} e^{-\delta^2} \ln \delta$ for $0< \alpha <1/2$, and that $J^{(\alpha)}_{\delta} \leq \delta^{\alpha} \int^{+\infty}_\delta  e^{-u^2} \ln u \h1 du$ for $-1/2< \alpha <0$.
Hence we get

\begin{equation}\label{fkoferrep}
J^{(1/2 - H)}_{2\sqrt{k \gamma}}  \leq \rho_7
\begin{cases}
{k}^{-1/4 - H/2} \ln k  & {\text for} \hspace{0.25cm} 0 < H < 1/2 \\
{k}^{1/4 - H/2}  &  {\text for} \hspace{0.25cm}  1/2 < H < 1 \\
\end{cases}
\end{equation}

and we finally obtain

\begin{equation}\label{aaaoeirjoekhrjfore"'rojvrorjortjo}
I^{(k)}_2 \leq \rho_8 \h1 \frac{\ln k}{{k}^{H/2-1/4}}.
\end{equation}

Using (\ref{zuiddhciuhcidsdhcisduhcdifuhck}) to (\ref{aaaoeirjoekhrjfore"'rojvrorjortjo}), item $3$ of theorem \ref{ozdicjdoisoedzijqpzeoejcvenvdsdlsiocfuvfsosfd} and the fact that both functions $H\mapsto \beta_H$ and $H\mapsto \tfrac{1}{c_H}$ are continuous on $[a,b]$ we get, for every $k$ in $\bN$,

\begin{align}\label{oeirjoekhrjfore"'rojvrorjortjo}
\underset{(H,u) \in [a,b]\times \bR }{\sup}{ \hspace{0.1cm} \left|\tfrac{\partial M_{H}}{\partial H}(e_k)(u)\right| }
		&\leq \rho_9 \h1 {(k+1)}^{2/3} + \rho_{10} \h1 \underset{H \in [a,b] }{\sup}{ \hspace{0.1cm}  \left[ ( {(k+1)}^{2/3 - H/2} +   {{(k+1)}^{1/4-H/2}} ) \ln(k+1)      \right]} \notag \\
		&\leq \rho_{11} \h1 {(k+1)}^{2/3} \ln (k+1).
\end{align}

$(ii)$ Let $t \in \bR$ and $r>0$ be fixed and define $[a,b] := [\underset{u \in  [t-r,t+r] }{\inf}{h(u)},\underset{u \in  [t-r,t+r] }{\sup}{h(u)}]$.  Using (\ref{firsodezr'kqsjdefrt}) we have, for every $k$ in $\bN$,

\begin{equation}
\begin{array}{lll}
\underset{u \in  [t-r,t+r] }{\sup}{  \left|\tfrac{d}{du} [g_{e_k}(u,h(u))] \right| } &\leq  \underset{(H,u) \in [a,b] \times \bR }{\sup}{ |M_H(e_k)(u)| }\notag \\
&+  (|t|+r)  \h1 \underset{u \in  [t-r,t+r] }{\sup}{ \h1|h'(u)|} \h1 \underset{(H,u) \in [a,b] \times \bR }{\sup}{ \left|\tfrac{\partial M_H}{\partial H}(e_k)(u)\right| }.
\end{array}
\end{equation}

The result then follows from $(i)$ above and item $3$ of theorem \ref{ozdicjdoisoedzijqpzeoejcvenvdsdlsiocfuvfsosfd}.
\end{pr}

\begin{rem} In the sequel, we will only need the bounds of $(i)$ and $(ii)$ with $(k+1)$ in lieu of ${(k+1)}^{2/3}$.
\end{rem}
We may now proceed to the  proof of theorem \ref{alalalalalal}.

\begin{praa}{\bfseries of  theorem \ref{alalalalalal}.}
From equality \eqref{fiizeufhirsozeedeijdefrt} defining multifractional white noise and equality (\ref{firsodezr'kqsjdefrt}), we have formally
\begin{equation}
W^{(h)}(t) = \displaystyle{ \sum^{+ \infty}_{k = 0}  M_{h(t)}(e_k)(t)} <.,e_k> + \h1 h'(t) \h1 \displaystyle{ \sum^{+ \infty}_{k = 0}  \int^{t}_0  \tfrac{\partial M_H(e_k)}{\partial H} (s) \big|_{ H = h(t)} ds <.,e_k> }.
\end{equation}

In order to establish that $W^{(h)}(t)$ is well defined in ${(\cS)}^{*}$ and that equality (\ref{firsozoijkqsjdefrt}) holds,
it is sufficient to show that both members on the right hand  side of the previous equality are in ${(\cS)}^{*}$.

For $t$ in $\bR$, definition (\ref{oerijefoijefjfijvoidfovijdfovijdfo}) of fractional white noise shows that $\displaystyle{ \sum^{+ \infty}_{k = 0}  M_{h(t)}(e_k)(t)} <.,e_k> = W^{(h(t))}(t)$ and thus belongs to ${(\cS)}^{*}$.

Let us show that $V_H(t) := \displaystyle{ \sum^{+ \infty}_{k = 0}  \int^{t}_0  \tfrac{\partial M_H(e_k)}{\partial H} (s) ds <.,e_k> }$ belongs to $(L^2)$. Using (\ref{jeveux2}) we may write

\begin{equation*}
\begin{array}{lll}
\bE[V^2_H(t)]	&= \sum^{+ \infty}_{k = 0} { {< \tfrac{\partial {M_H}}{\partial H }(e_k), \i1_{[0,t]} >}^2_{L^2(\bR)}  } = \sum^{+ \infty}_{k = 0}{ {<  e_k, \tfrac{\partial {M_H}}{\partial H }(\i1_{[0,t]}) >}^2_{L^2(\bR)}  } \\
				&= {||\tfrac{\partial {M_H}}{\partial H }(\i1_{[0,t]})||}^2_{L^2(\bR)} = {||\i1_{[0,t]}||}^2_{\delta_H} < +\infty .
\end{array}
\end{equation*}

As a consequence, $W^{(h)}(t)$ is the sum of an ${(\cS)}^{*}$ process and an $(L^2)$ process, and thus belongs to ${(\cS)}^*$.
We are left with proving equality (\ref{firsozkqsjdefrt}), \textit{i.e.} that $W^{(h)}(t)$ is indeed the ${(\cS)}^{*}$ derivative of
$B^{(h)}(t)$ for any real $t$.

Let $r \neq 0$ and  $t \geq 0$ (the case $t <0$ follows in a similar way). The equality $W^{(h)}(t) = W^{(h(t))}(t) + h'(t) \h1 V_{h(t)}(t)$ and remark \ref{iooisoidjoissid} entail
that $W^{(h)}(t)$ belongs to $({\cS_{-p}})$ as soon as $p \geq 2$. For such a $p$, one computes:

\begin{align}
J_{p,r}(t) 	&:= {\left|\left|  \tfrac{B^{(h)}(t+r) - B^{(h)}(t)}{r}  - W^{(h)}(t) \right|\right|}^2_{-p} \notag \\
			&=  {\left|\left| \sum^{+\infty}_{k = 0} \h1     \left[ \left( \frac{ g_{e_k}(t+r,h(t+r)) -  g_{e_k}(t,h(t)) }{r}  \right) - \frac{d}{dt}[g_{e_k}(t,h(t))]   \right] <.,e_k>   \right|\right|}^2_{-p} \notag \\
			&= \sum^{+\infty}_{k = 0} \h1  \underbrace{\frac{1}{{(2k+2)}^{2p}} \h1     {\left[ \left( \frac{ g_{e_k}(t+r,h(t+r)) -  g_{e_k}(t,h(t)) }{r}  \right) - \frac{d}{dt}[g_{e_k}(t,h(t))]   \right]}^2 }_{ :=J_{p,r,k}(t) }.
\end{align}

Using lemma \ref{zoidpsqojxiqsoxijqsojdzj} and the Mean-Value theorem we obtain, for $r$ in $(-1/2,1/2)\backslash\{0\}$:

\begin{equation}\label{zoidjdzjdeoiz}
\begin{array}{lll}
J_{p,r,k}(t)	&\leq  \h1  \frac{1}{{(2k+2)}^{2p}} \h1 { \left(2 \h1    \underset{u \in [t-1/2,t+1/2]  }{\sup}{ \h1 \left|   \tfrac{d}{du} [g_{e_k}(u,h(u)) ] \right|  } \right) }^2 \leq 4 \h1  \widetilde{D} \h1   \frac{{(k+1)}^2}{{(2k+2)}^{2p}}  \notag \\
				&\leq \frac{\widetilde{D}}{2^{2(p-1)}}  \h1  \frac{1}{{(k+1)}^{2(p-1)}},
\end{array}
\end{equation}

where $\widetilde{D}:= \widetilde{D}_t(1/2)$.
Since $J_{p,r,k}(t) \underset{r \to 0 }{\longrightarrow}{ 0 }$, equality (\ref{firsozkqsjdefrt}) follows from the dominated convergence theorem.
\end{praa}

\begin{rem}
In  $(ii)$ of lemma $\ref{zoidpsqojxiqsoxijqsojdzj}$, the real constant $\widetilde{D_t}(r)$ can be taken independent of $t$ if the function $t\mapsto t\h1 h'(t)$ is bounded over $\bR$.
\end{rem}

We note that multifractional white noise is a sum of two terms: a fractional white noise that belong to $({\cS_{-p}})$ as soon as $p \geq 2$,
and a "smooth" term which corresponds to the derivative in the "$H$" direction. This is a direct consequence of the fact that the fractional Brownian
field $\Lambda(t,H)$ is not differentiable in the $t$ direction (in the classical sense) but infinitely smooth in the $H$ direction.

\begin{prop}\label{mwncont}
For $p\geq 2$, the map $t\mapsto {||W^{(h)}(t)||}_{-p}$ is continuous.
\end{prop}

\begin{pr}
By definition,
$${||W^{(h)}(t)||}^2_{-p} = \sum^{+ \infty}_{k = 0} \frac{{( \tfrac{d}{dt}[g_{e_k}(t,h(t))])}^2}{{(2k+2)}^{2p}}.$$
Using the estimate given in lemma \ref{zoidpsqojxiqsoxijqsojdzj} (ii), we see that ${||W^{(h)}(t)||}^2_{-p}$ is the sum of a series of continuous functions that converges normally on any compact.
\end{pr}

\subsection{Generalized functionals of mBm}
\label{eprofkerpokveprvkerpovepovk}

In the next section, we will derive various Itô formulas for the integral with respect to mBm. It will be useful to obtain such
formula for tempered distributions. In that view, we define generalized functionals of mBm as in \cite{ben1}.
		
\begin{theodef}\label{lmqsdlmflkjdg}
Let $F$ be a tempered distribution. For $t$ in $\bR^*_+$, define
\begin{equation}
F(B^{(h)}(t)):= \frac{1}{\sqrt{2\pi} t^{h(t)}} \sum\limits^{+\infty}_{k=0} {(k!)}^{-1} t^{-2 k h(t)} <F, {\xi}_{(t,h(t),k)}> I_k\left(   {(\hspace{0.1cm}M_{h(t)}(\i1_{[0,t]})\hspace{0.1cm})}^{\otimes k} \right)
\end{equation}

where the functions ${\xi}_{t,H,k}$ are defined for $(x,H,k)$ in $\bR\times(0,1)\times \bN$ by

\begin{equation}
{\xi}_{t,H,k}(x) :={(\sqrt{2})}^{-k} \h1 t^{k H} \h1 h_{k}(x/{(\sqrt{2} t^H)})  \h1 {\exp}{\{-\frac{x^2}{2 \h1 t^{2H}}\}} = {\pi}^{1/4} {(k!)}^{1/2} t^{k H} {\exp}{\{-\frac{x^2}{4t^{2H}}\}} {e}_{k}{(x/(\sqrt{2} t^H))}.
\end{equation}

Then for all $t$ in $\bR^*_+$, $F(B^{h}(t))$ is a Hida distribution, called {\it generalized functional of $B^{(h)}(t)$}.
\end{theodef}

\begin{pr}
This is an immediate consequence of \cite{Kuo2} p.$61$-$64$ by taking $f := M_{h}(\i1_{[0,t]})$.
\end{pr}

\begin{rem}
As shown in \cite{ben1}, when $F=f$ is of function type, $F(B^{h}(t))$ coincides with $f(B^{h}(t))$.
\end{rem}

The following theorem yields an estimate of ${||F(B^{(h)}(t)||}^2_{-p}$ which will be useful in the sequel.

\begin{theo}
Let $h:\bR \rightarrow [H_1,H_2] \subset (0,1)$ be a continuous function, $B^{(h)}$ an mBm, $p\in \bN$ and $F \in {\sS}_{\hspace{-0.15cm}-p}(\bR)$. Then there is a constant $C_p^{(H_1,H_2)}$, independent of $F$, such that

\begin{equation}\label{spdsdpokpovkopsdkvspodvksopsdvpoksvdokdsvkopsvdopqbschfgsdvhjerhbtkljtyopfdhvlkdfhu}
\forall t>0, \hspace{1cm} {||F(B^{(h)}(t))||}^2_{-p} \leq \max\{t^{-4ph(t)},t^{4ph(t)}\} t^{-h(t)} \h1 C_p^{(H_1,H_2)} {|F|}^2_{-p}.
\end{equation}
\end{theo}

\begin{pr}
For $H \in (0,1)$ and $p \in \bN$, Theorem $3.3$ p.$92$ of \cite{ben1} ensures that there exists $C^{(H)}_p$ such that, $\forall t>0$,

\begin{equation*}
{||F(B^{H}(t))||}^2_{-p} \leq \max\{t^{-4pH},t^{4phH}\} t^{-H} \h1 C_p^{(H)} {|F|}^2_{-p}.
\end{equation*}

Now if $H$ belongs to $[H_1,H_2]$, it is easy to show, by examining closely the iteration of  $(23)$ p.$94$ in \cite{ben1}, that one can choose a constant $C^{(H_1,H_2)}$ independent of $H$. We hence have

\begin{equation}\label{ssqssqu}
\forall t>0,  \hspace{0.25cm} \forall H \in [H_1,H_2], \hspace{1cm}  {||F(B^{H}(t))||}^2_{-p} \leq \max\{t^{-4pH},t^{4phH}\} t^{-H} \h1 C_p^{(H_1,H_2)} {|F|}^2_{-p}.
\end{equation}

For $t>0$, one only needs to set $H=h(t)$ in (\ref{ssqssqu}) to get (\ref{spdsdpokpovkopsdkvspodvksopsdvpoksvdokdsvkopsvdopqbschfgsdvhjerhbtkljtyopfdhvlkdfhu}).
\end{pr}

\subsection{S-Transform of mBm and multifractional white noise}

The following theorem makes explicit the $S$-transforms of 
mBm, multifractional white noise and generalized functionals of mBm.

We denote by $\gamma$ the heat kernel density on $\bR_+\times\bR$ \textit{i.e} $\gamma(t,x):= \frac{1}{\sqrt{2\pi t}}\exp{\{ \frac{-x^2}{2t}  \}}$ if $t \neq 0$ and $0$ if $t = 0$.

\begin{theo}\label{tardileomalet}
Let $h:\bR \rightarrow (0,1) $ be a $\cC^1$  function and ${(B^{(h)}(t))}_{t \in \bR}$ (resp. ${(W^{(h)}(t))}_{t \in \bR}$) be an mBm
(resp. multifractional white noise). For $\eta \in \sS(\bR)$ and $t \in \bR$,

\bit
\iti $S[B^{(h)}(t)](\eta) = {<\eta,M_{h}(\i1_{[0,t]})>}_{L^2(\bR)} = g_{\eta}(t,h(t)) $, where $g_{\eta}$ has been defined in lemma $\ref{zekkpoerkpeorperokvpo}$.
\itii  $S[W^{(h)}(t)](\eta) = \frac{d}{dt}[g_{\eta}(t,h(t))] =  M_{h(t)}(\eta)(t) + h'(t)  \displaystyle{\int^{t}_0  \tfrac{\partial M_H}{\partial H}(\eta)(s) \big|_{ H = h(t)} ds } $.
\itiii
For $p \in \bN$ and $F \in \sS_{-p}(\bR)$,

\cl{$S[F(B^{(h)}(t))](\eta) = \left< F, \gamma \left(t^{2h(t)}, .-\displaystyle{\int^t_0 \h1  M_{h(t)} (\eta)\h1(u) \h1 du}  \right) \right>$.}

Furthermore, there exists a constant $C_p$, independent of $F,t$ and $\eta$, such that
\begin{equation}\label{periojv}
 {|S[F(B^{(h)}(t))](\eta)|}^2 \leq \max\{t^{-4ph(t)},t^{4ph(t)}\} t^{-h(t)} \h1 C_p \h1 {|F|}^2_{-p}  \h1 \exp\{{|A^p \eta |}^2_{0}   \}.
\end{equation}
\eit
\end{theo}

\begin{pr}

$(i)$ Thanks to  (\ref{fiizeufhirsozeedeijdt}) and lemma \ref{dpozepdozkekpdzkeodkpozekpdozkepdozekodezpdkezpd} we have, for every $\eta$ in $\sS(\bR)$ and $t$ in $\bR$,
\begin{multline*}
S(B^{(h)}(t))(\eta) =  \sum^{+ \infty}_{k = 0} {<M_{h(t)}(\i1_{[0,t]}),e_k>}_{L^2(\bR)}  {<\eta,e_k>}_{L^2(\bR)} = {<M_{h(t)}(\i1_{[0,t]}),\eta>}_{L^2(\bR)}\\
  = {<\i1_{[0,t]},M_{h(t)}(\eta)>}_{L^2(\bR)} = g_{\eta}(t,h(t)) .
\end{multline*}

$(ii)$ This is a straightforward consequence of lemma \ref{dkdskcsdckksdksdmksdmlkskdm}, (\ref{firsozkqsjdefrt}) and $(i)$.

$(iii)$ The first equality results from theorem $7.3$ p.$63$ in \cite{Kuo2} with $f = M_{h(t)}(\i1_{[0,t]})$ and from the equality $(i)$. Equality (\ref{periojv}) results from $(\ref{spdsdpokpovkopsdkvspodvksopsdvpoksvdokdsvkopsvdopqbschfgsdvhjerhbtkljtyopfdhvlkdfhu})$ as in theorem $3.8$ p.$95$ of  \cite{ben1}.
\end{pr}

\begin{rem}\label{eoifieroij}
Using lemma $\ref{dpozepdozkekpdzkeodkpozekpdozkepdozekodezpdkezpd}$ and  $(\ref{fiizeufhirsozeedeijdefrt})$  we may also write:
\begin{equation}\label{oeqqqqiefjrodif}
\forall (t,\eta) \in \bR \times \sS(\bR),\hspace{0.75cm}  S(W^{(h)}(t))(\eta) =  \sum^{+ \infty}_{k = 0} \h1 \tfrac{d}{dt}[g_{e_k}(t,h(t))]\h1 {<\eta,e_k>}_{L^2(\bR)}.
\end{equation}
\end{rem}

\subsection{The multifractional Wick-Itô integral}\label{ozieffjioefjoezifjioezzjioefjfeoiz}						
We are now able to define the Multifractional Wick-Itô integral, in a way analogous  to the  definition of the fractional Wick-Itô integral.
\textit{In the sequel of this work, we will always assume that $h$ is a $\cC^1$ function on $\bR$ with bounded derivative.}			
			
\begin{defi}[The multifractional Wick-Itô integral]\label{oezifhherioiheroiuh}
Let $Y:\bR \rightarrow {(\cS)}^*$ be a process such that the process $t \mapsto Y(t)\diamond W^{(h)}(t)$ is  ${(\cS)}^*$-integrable on $\bR$. We then say that the process $Y$ is $dB^{(h)}$-integrable on $\bR$ or integrable  on $\bR$ with respect to  mBm $B^{(h)}$. The integral of $Y$ with respect to $B^{(h)}$ is defined by

\begin{equation}\label{eigfrretth}
\int_{\bR} Y(s) \h1 dB^{(h)}(s) :=  \int_{\bR} Y(s)\diamond W^{(h)}(s) ds.
\end{equation}
For a Borel set $I$ of $\bR$, define $\int_{I} Y(s) dB^{(h)}(s):= \int_{\bR} \h1 {\i1}_{I}(s) \h1 Y(s) dB^{(h)}(s)$.
\end{defi}

When the function $h$ is constant, the multifractional Wick-Itô integral coincides with the fractional Itô integral defined in \cite{ell}, \cite{bosw}, \cite{ben1} and \cite{ben2}.  In particular, when the function $h$ is identically $1/2$, \eqref{eigfrretth} is nothing
but the classical Itô integral with respect to  Brownian motion, provided of course $Y$ is Itô-integrable. The multifractional Wick-Itô integral verifies the following properties:

\begin{prop}\label{zpofkpsokfepsofkspkod}

\bit
\iti Let $(a,b)$ in $\bR^2$, $a<b$. Then $ \int^b_a \h1  dB^{(h)}(u) = B^{(h)}(b) - B^{(h)}(a)$ almost surely.
\itii Let $X:I\rightarrow ({\cS}^*)$ be a  $dB^{(h)}$-integrable process over $I$, a Borel subset of $\bR$. Assume
$\int_{I} X(s) dB^{(h)}(s)$ belongs to $(L^2)$. Then $\bE[\int_{I} X(s) dB^{(h)}(s)] = 0$.
\eit
\end{prop}

\begin{pr}
\textit{$(i)$} From $(ii)$ of theorem \ref{tardileomalet}, $t\mapsto S(\i1_{[a,b]}(t) \h1 W^{(h)}(t))(\eta)$ is measurable on $\bR$ for any $\eta$ in $\sS(\bR)$. Moreover, for any integer $p_0 \geq 2$, we have

\begin{equation*}
|S(\i1_{[a,b]}(t) \h1 W^{(h)}(t))(\eta)| \leq {||W^{(h)}(t)||}_{-p_0}  e^{\frac{1}{2}{|\eta|}^2_{p_0}}
\end{equation*}

thanks to lemma \ref{dpozepdozkekpdzkeodkpozekpdozkepdozekodezpdkezpd}. By proposition \ref{mwncont}, $t \mapsto {||W^{(h)}(t)||}_{-p_0}$
is continuous thus integrable on $[a,b]$. Theorem \ref{peodcpdsokcpodfckposkcdpqkoq} then entails  that $t\mapsto \i1_{[a,b]}(t) \h1 W^{(h)}(t)$ is $(\cS^*)$-integrable over $\bR$.
It is easily seen that the S-transforms of $ \int^b_a \h1  dB^{(h)}(u)$ and $B^{(h)}(b) - B^{(h)}(a)$ coincide.
The result then follows from the injectivity of the S-transform.

\par \noindent \textit{$(ii)$} One computes:
\vspace{-7mm}
\begin{equation*}
S(\int_{I} X(s) dB^{(h)}(s))(0) = \int_{I} S(X(s))(0)  \overbrace{S(W^{(h)}(s))(0)}^{ = \frac{d}{ds}[g_{0}(s,h(s))] = 0 }   \h1 ds = 0.
\end{equation*}
Now, when $\int_{I} X(s) dB^{(h)}(s)$ belongs to $(L^2)$,  $\bE[\int_{I} X(s) dB^{(h)}(s)]=S(\int_{I} X(s) dB^{(h)}(s))(0)$.
\end{pr}

\begin{theo}\label{criteredintegrabilite}
Let $I$ be a compact subset of $\bR$ and  $X:t\mapsto X(t)$ be a  process from $I$ to $(L^2)$ such that $t \mapsto S(X(t))(\eta)$ is measurable on $I$ for all $\eta$ in $\sS(\bR)$ and $t \mapsto {||X(t)||}_{0}$ belongs to $L^1(I)$. Then $X$ is $dB^{(h)}$-integrable on $I$ and there exist a natural integer $q$ and a constant $C_I$ such that,

\begin{equation}\label{fgsdqhdhqsfg}
{\left|\left|\int_I X(t)      \h1 dB^{(h)}(t)  \right|\right|}_{-q} \leq C_I \h1 \int_I {||X(t)||}_{0} \h1 dt.
\end{equation}
\end{theo}
\begin{pr}
For $\eta \in \sS(\bR)$, the measurability on $I$ of  $t \mapsto S(X(t)\diamond W^{(h)}(t))(\eta)$ is clear since\\
$S(X(t)\diamond W^{(h)}(t))(\eta) =  S(X(t))(\eta) \h1 \frac{d}{dt}[g_{\eta}(t,h(t))]$.
By lemma \ref{dede}, we have, for any integer $q \geq 2$,
\begin{equation*}
|S(X(t) \h1 \diamond \h1 W^{(h)}(t))(\eta)|  \leq   {||X(t)||}_{0} \h1 {||W^{(h)}(t)||}_{-q}  \h1 e^{ {|\eta|}^2_{q}}
\end{equation*}

%
for every $t$ in $I$. Since $t\mapsto {||W^{(h)}(t)||}_{-q}$ is continuous by proposition \ref{mwncont} and $t \mapsto {||X(t)||}_{0}$ belongs to $L^1(I)$ by assumption, the result follows from theorem \ref{peodcpdsokcpodfckposkcdpqkoq}. We refer to theorem 13.5 of \cite{Kuo2} for the upper bound.
\end{pr}

\begin{rem}
\label{allahdchgdgdg}
One can show, using appendix $\ref{appendiceB}$, that inequality $(\ref{fgsdqhdhqsfg})$ is true for every integer $q \geq 2$.
\end{rem}

It is of interest to have also a criterion of integrability for generalized functionals of mBm. In that view, we set up the following notation: for $p \in \bN$, $0<a<b$, we consider a map $F:[a,b] \rightarrow \sS_{-p}(\bR)$ (hence  $F(t)$ is a tempered distribution for all $t$). We then define  $F(t,B^{(h)}(t)) := F(t)(B^{(h)}(t))$.
Recall the following theorem (see \cite{GeCh2}, lemma $1$ and $2$ p.$73$-$74$):

\begin{theo}
\label{peogi}
Let $I$ be  an interval of $\bR$, $t \mapsto F(t)$ be a map  from $I$ into ${\sS}_{-p}(\bR)$,  $t \mapsto \varphi(t,.)$ be a map from $I$ into $\sS(\bR)$ and $t_0 \in I$. If both maps $t \mapsto F(t)$ and  $t \mapsto \varphi(t,.)$ are continuous (respectively differentiable) at $t_0$,  then the function $t \mapsto <F(t),\varphi(t,.)>$ is continuous (respectively differentiable) at $t_0$.
\end{theo}

\begin{theo}
\label{dlsdqsdhuqsgdudguqsgdqudgquygdsudgqsuydgqsuydgqsuyqdgq}
Let $p \in \bN$, $0<a<b$ and let $F:[a,b] \rightarrow \sS_{-p}(\bR)$ be a continuous map. Then the stochastic distribution process $F(t,B^{(h)}(t))$ is both ${(\cS)^*}-$integrable and $dB^{(h)}$-integrable over $[a,b]$.
\end{theo}

\begin{pr}
We shall apply theorem \ref{peodcpdsokcpodfckposkcdpqkoq}.

The measurability of $t \mapsto  S[ F(t,B^{(h)}(t))](\eta)$ results from $(iii)$ of theorem \ref{tardileomalet}, the continuity of the two
maps $t \mapsto F(t)$ and $t \mapsto \gamma \left(t^{2h(t)}, .-{\int^t_0 \h1  M_{h(t)} (\eta)\h1(u) \h1 du}  \right)$ and theorem \ref{peogi}.

Since $h$ is bounded on $[a,b]$, lemma \ref{dpozepdozkekpdzkeodkpozekpdozkepdozekodezpdkezpd} and (\ref{spdsdpokpovkopsdkvspodvksopsdvpoksvdokdsvkopsvdopqbschfgsdvhjerhbtkljtyopfdhvlkdfhu}) yield

\begin{align}\label{aleszozos}
{|S[F(t,B^{(h)}(t))](\eta)|}	&\leq \max\{t^{-2ph(t)},t^{2h(t)}\} t^{-h(t)/2} \h1 \sqrt{C^{(H_1,H_2)}_p} \h1 \underset{s \in [a,b] }{\max}{ {{|F(s)|}_{-p}}} \h1 \exp\{\frac{1}{2}{|A^p \eta |}^2_{0}   \} \notag\\
								&\leq ( {(\tfrac{1}{a})}^{2p H_2} + {b}^{2p H_2} ) ( {(\tfrac{1}{a})}^{ \frac{H_2}{2}} + {b}^{ \frac{H_2}{2}} )    \h1 \sqrt{C^{(H_1,H_2)}_p} \h1 \underset{s \in [a,b] }{\max}{ {{|F(s)|}_{-p}}} \h1 \exp\{\frac{1}{2}{|A^p \eta |}^2_{0}   \},
\end{align}

where  $H_1:=\underset{s \in [a,b] }{\min}{ h(s)}$     and $H_2:= \underset{s \in [a,b] }{\max}{ h(s)}$.
This yields the second condition of theorem \ref{peodcpdsokcpodfckposkcdpqkoq} and shows that $F(t,B^{(h)}(t))$ is ${(\cS)^*}-$integrable over $[a,b]$.

For $dB^{(h)}$-integrability, we first note that, by theorem \ref{tardileomalet} $(ii)$,
$$S[ F(t,B^{(h)}(t))\diamond W^{(h)}(t) ](\eta) = S[ F(t,B^{(h)}(t))](\eta) \h1 \frac{d}{dt}[g_{\eta}(t,h(t))].$$
Since the function $t \mapsto \frac{d}{dt}[g_{\eta}(t,h(t))]$ is continuous (by lemma \ref{zekkpoerkpeorperokvpo}), the measurability of $t \mapsto S[ F(t,B^{(h)}(t))\diamond W^{(h)}(t) ](\eta)$ for every function $\eta$ in $\sS(\bR)$ follows.

Moreover, for every integer $p_0 \geq \max\{p,2\}$, $F(t,B^{(h)}(t))$ and $W^{(h)}(t)$ belong to $({\cS}_{-p_0})$ for all $t$ in $[a,b]$. Using
lemma \ref{dede} and (\ref{spdsdpokpovkopsdkvspodvksopsdvpoksvdokdsvkopsvdopqbschfgsdvhjerhbtkljtyopfdhvlkdfhu}),
we may write, for all $t$ in $[a,b]$,

\begin{equation}\label{uiuhhsuqhxshushuqxsxhxsh}
\begin{array}{ll}
\left| S[ F(t,B^{(h)}(t))\diamond W^{(h)}(t) ](\eta)\right| &\leq {||F(t,B^{(h)}(t))||}_{-p_0} \h1 {||W^{(h)}(t)||}_{-p_0} \h1  \exp{ \{{|\eta |}^2_{p_0}} \}\\
															&\leq K \h1  \exp\{{|\eta |}^2_{p_0} \},
\end{array}
\end{equation}

%
where $K:= \underset{t \in [a,b] }{\sup}{{||W^{(h)}(t)||}_{-p_0}} \h1 ( {(\tfrac{1}{a})}^{2p_0 H_2} + {b}^{2p_0 H_2} ) ( {(\tfrac{1}{a})}^{ \frac{H_2}{2}} + {b}^{ \frac{H_2}{2}} )    \h1 \sqrt{C^{(H_1,H_2)}_{p_0}} \h1 \underset{s \in [a,b] }{\max}{ {{|F(s)|}_{-p_0}}} $.

Theorem \ref{peodcpdsokcpodfckposkcdpqkoq} applies  again and shows that $t \mapsto F(t,B^{(h)}(t))\diamond W^{(h)}(t)$ is integrable  over $[a,b]$.
\end{pr}

\begin{rem}
Recall that a function $f:\bR \rightarrow \bR$  is said to be of polynomial growth if there is an integer $m$ in $\bN$ and a constant $C$ such that for all $x \in \bR$, $|f(x)| \leq C(1+|x|^m)$.
The previous theorem entails in particular that both quantities $\int^b_a \h1 f(B^{(h)}(t)) \h1 dt$ and $\int^b_a \h1 f(B^{(h)}(t)) \h1 dB^{(h)}(t)$ exist in ${(\cS)}^*$ if $f$ is a function of polynomial growth.
\end{rem}

\begin{exa}[Computation of $\int^T_0 \h1 B^{(h)}(t) \h1 dB^{(h)}(t)$]
 Let $T>0$ fixed. Then 

\begin{equation}\label{erfepooqoiazpiesssjoidkeoh}
I:= \int^T_0 \h1 B^{(h)}(t) \h1 dB^{(h)}(t) = \int^T_0 \h1 W^{(h)}(t) \diamond  B^{(h)}(t) \h1 dt =  \int^T_0 \h1 \frac{dB^{(h)}(t)}{dt} \diamond  B^{(h)}(t) \h1 dt.
\end{equation}

Let us prove that the last quantity is equal to $\frac{1}{2} {B^{(h)}(T)}^{\diamond 2} := \frac{1}{2} \h1 (B^{(h)}(T)\diamond B^{(h)}(T))=\frac{1}{2} \h1 {(B^{(h)}(T)}^{2}- T^{2h(T)})$ $($see  remark $\ref{qmqmqmmmmm})$.
It is sufficient to compute the $S$-transforms of  both sides of the equality.\hspace{0.05cm}For $\eta$ in $\sS(\bR)$, 

\begin{align}\label{erfepooqorgtiazpiesssjoidkeoh}
S\left(\int^T_0 \hspace{-0.1cm}B^{(h)}(t) dB^{(h)}(t)\right)(\eta)	&=\int^T_0 S(B^{(h)}(t))(\eta)  S(W^{(h)}(t))(\eta) dt =  \int^T_0 \hspace{-0.1cm} g_{\eta}(t,h(t))  \tfrac{d}{dt}[g_{\eta}(t,h(t))] dt \notag \\
																	&= \frac{1}{2} g^2_{\eta}(T,h(T))  = \frac{1}{2}  {(S(B^{(h)}(T))(\eta))}^2 = \frac{1}{2} \h1  S\left(B^{(h)}(T)\diamond B^{(h)}(T)\right)(\eta)  \notag \\
																	&= S\left(\frac{1}{2} \h1 {(B^{(h)}(T)}^{2}- T^{2h(T)})\right)(\eta).
\end{align}

\end{exa}

To end this section, we present a simple but classical stochastic differential equation in the frame of mBm.

\begin{ex}[The multifractional Wick exponential]
Following  $\cite{ell} $ formula $(4.8)$ and  $\cite{bosw}$ example $3.6$, let us consider the multifractional stochastic differential equation

\begin{align}\label{diopjfcoizdezdzedzezedj}
{  \left\{
\begin{aligned}
dX(t) &=  \alpha(t) X(t) dt + \beta(t) X(t) dB^{(h)}(t)   \\
X(0) &\in   (\cS^*),
\end{aligned}
\right. }
\end{align}

where $t$ belongs to $\bR_+$ and where $\alpha:\bR\rightarrow\bR$  and  $\beta:\bR\rightarrow\bR$ are two deterministic continuous functions. $(\ref{diopjfcoizdezdzedzezedj})$ is a shorthand notation for

\begin{equation}\label{eifjreiofjferojefjfdjodijoif}
X(t) = X(0) + \int^t_0 \h1  \alpha(s) \h1 X(s) \h1 ds + \int^t_0 \h1  \beta(s) \h1 X(s) \h1 dB^{(h)}(s),
\end{equation}

where the previous equality holds in ${(\cS)}^*$. Rewrite the previous equation in terms of derivatives in ${(\cS)}^*$  as:

\vspace{-0.5cm}

\begin{align}
{  \left\{
\begin{aligned}
\label{diopjfcoidj}
\frac{dX}{dt}(t) &=   \alpha(t) \h1 X(t)   + \beta(t) \h1 X(t) \h1 \diamond  W^{(h)}(t) = (\alpha(t)+ \beta(t) W^{(h)}(t) )  \diamond  X(t) \\
X(0) &\in   (\cS^*).
\end{aligned}
\right. }
\end{align}

We thus are looking for an $(\cS^*)$-process, noted $Z$, defined on $\bR_+$ such that $Z$ is differentiable on  $\bR_+$ and verifies equation $(\ref{diopjfcoidj})$ in $(\cS^*)$. As in \cite{bosw}, it is easy to guess the solution of $(\ref{diopjfcoidj})$ if we replace  Wick products $\diamond$ by ordinary products. Once we have a solution of  $(\ref{diopjfcoidj})$, we replace ordinary products by  Wick products. This heuristic reasoning leads to defining the process $Z:=({Z(t)}_{t \in \bR_+})$ by

\begin{equation}\label{msmslslqmls}
Z(t) := X(0) \diamond \exp^\diamond { \bigg(  \int^t_0 \alpha(s) ds +   \int^t_0 \beta(s) dB^{(h)}(s) \bigg)} ,\hspace{0.5cm} t \in \bR_+,
\end{equation}

where $\exp^\diamond$ has been defined in section $\ref{ksksksks}$. 
%
%
%

\begin{theo}
The process $Z$ defined by $(\ref{msmslslqmls})$ is the unique solution in $(\cS^*)$ of $(\ref{diopjfcoidj})$.
\end{theo}

\begin{pr}
This is a straightforward application of theorem $3.1.2$ in \cite{HOUZ}.
\end{pr}
\end{ex}

\begin{rem}
$(i)$  \cite{HOUZ} uses the Hermite transform in order to establish the theorem. However it is possible to start from (\ref{eifjreiofjferojefjfdjodijoif}), take $S$-transforms of both sides and solve the resulting ordinary stochastic differential equation.

$(ii)$ Equation $(\ref{diopjfcoizdezdzedzezedj})$ may be solved with other assumptions on $\alpha, \beta$. We refer to $\cite{JLJLV2}$ for more on stochastic differential equations driven by mBm. 
\end{rem}

\begin{rem}
In particular when $X(0)$ is deterministic, equal to $x$, $\alpha() \equiv \alpha$ and $\beta() \equiv \beta$ are constant functions, the solution $X$ of $(\ref{eifjreiofjferojefjfdjodijoif})$  reads

\begin{equation}
X(t) = x \exp{\{ \beta B^{(h)}(t) + \alpha t  - \tfrac{1}{2} \beta^2 t^{2h(t)} \}}, \hspace{0.5cm} t  \in \bR_+,
\end{equation}

which is analogous to formula $(3.31)$ given in  $\cite{bosw}$ in the case of the fractional Brownian motion.
\end{rem}

\subsection{Multifractional Wick-Itô integral of deterministic elements  versus Wiener integral with respect to mBm}

In section \ref{Wiener}, we have defined a Wiener integral with respect to mBm. It is natural to check whether this definition
is consistent with the multifractional Wick-Ito integral when the integrand is deterministic. More precisely, we wish to
verify that $\int_{\bR} f(s)\diamond W^{(h)}(s) ds =  {\cJ}^{h}(f)$ for all functions $f$ such that both members of
the previous equality exist and that the left-hand side member is in $(L^2)$. In that view we first prove the following theorem.

\begin{theo}
\label{sdsiduhsdiuhfsdiudfhdsufhisduhfisduhfsiuhfsidhfiuhdfisudhfiu}
Let $f:\bR\rightarrow \bR$ be a   deterministic function which belongs to  $L^1_{\text{loc}}(\bR)$. Let $Z:=({Z(t)}_{t \in \bR})$ be the process defined on $\bR$ by $Z(t):= \int^t_{0}  f(s) dB^{(h)}(s)$. Then $Z$ is an $(\cS^*)$-process which verifies the following equality in $(\cS^*)$

\begin{equation}\label{pdivjdofijvdoijfdoij}
\int^t_{0}  f(s) dB^{(h)}(s) =  \sum^{+ \infty}_{k = 0}  \left( \int^t_{0} f(s) \tfrac{d}{ds}[g_{e_k}(s,h(s))] ds  \right)     {<.,e_k>}.
\end{equation}

Moreover $Z$ is a (centered) Gaussian process if and only if $\sum^{+ \infty}_{k = 0}  {\left( \int^t_{0} f(s) \tfrac{d}{ds}[g_{e_k}(s,h(s))] ds  \right)}^2~<~+\infty$, for all $t$. In this  case,
\begin{equation}\label{pngphopkhghkponpokg}
Z(t) = \int^t_{0}  f(s) dB^{(h)}(s) \leadsto \cN\left(0, \sum^{+ \infty}_{k = 0}  {\left( \int^t_{0} f(s) \tfrac{d}{ds}[g_{e_k}(s,h(s))] ds  \right)}^2\right), \hspace{1cm} \forall t \in \bR.
\end{equation}

In particular, the process $Z$ is Gaussian when the function $f$ belongs to $\cC^1(\bR,\bR)$ and is  such that  $\underset{t \in \bR }{\sup}{|f'(t)|} < +\infty$.
\end{theo}

\begin{pr}
We treat only the case $t \in \bR^*_+$. The other case follows similarly. Let $f$ be in $L^1_{\text{loc}}(\bR)$. In order to show (\ref{pdivjdofijvdoijfdoij}) let us establish {\bfseries a), b)} and {\bfseries c)} below.

{\bfseries a)} $s\mapsto f(s) \diamond W^{(h)}(s)$ is ${(\cS)}^*$-integrable over $[0,t]$.

Let $\eta \in \sS(\bR)$ and $s$ in $[0,t]$, using lemma \ref{dpozepdozkekpdzkeodkpozekpdozkepdozekodezpdkezpd}, we get:
\vspace{-5mm}
\begin{equation*}
|S(f(s)\diamond W^{(h)}(s))(\eta)| = |f(s)|\h1  \left| S(W^{(h)}(s))(\eta) \right| \leq \overbrace{|f(s)| \h1  {||W^{(h)}(s)||}_{-p_0}}^{=:L(s)} e^{ {\tfrac{1}{2} |\eta|}^2_{p_0}} \h1
\end{equation*}
for $s$ in $[0,t]$ and for $p_0 \geq 2$. Since $L$ is the product of a continuous function and a function of $L^1_{\text{loc}}(\bR)$, {\bfseries a)} is a consequence of theorem \ref{peodcpdsokcpodfckposkcdpqkoq}.

{\bfseries b)} $\Psi_f:=\sum^{+ \infty}_{k = 0}  \left( \int^t_{0} f(s) \tfrac{d}{ds}[g_{e_k}(s,h(s))] ds  \right) {<.,e_k>}$ belongs to $({\cS_{-p_0}})$  as soon as $p_0 \geq 2$.

Lemma \ref{zoidpsqojxiqsoxijqsojdzj} entails that there exists a real $D$ such that, for every  $p_0\geq 2$, we have

\begin{equation*}
\begin{array}{ll}
{||\Psi_f||}^2_{-p_0}	&= \sum^{+ \infty}_{k = 0} {\left(\frac{ \int^t_0  f(s)  \h1 \tfrac{d}{ds}[g_{e_k}(s,h(s))] \h1 ds}{{(2k+2)}^{p_0}}  \right)}^2 \leq {||f||}^2_{L^1([0,t])} \h1 \sum^{+ \infty}_{k = 0} \frac{ { \h1\underset{s \in [0,t] }{\sup}{|\tfrac{d}{ds}[g_{e_k}(s,h(s))] |}   }^2 }{{(2k+2)}^{2p_0}}\\
						&\leq D^2 \h1 {||f||}^2_{L^1([0,t])} \h1  \sum^{+ \infty}_{k = 0} \frac{ {(k+1)}^2}{{(2k+2)}^{2p_0}} <  \h1 +\infty.
\end{array}
\end{equation*}

{\bfseries c)} $\int^t_{0}  f(s) dB^{(h)}(s) = \sum^{+ \infty}_{k = 0}  \left( \int^t_{0} f(s) \tfrac{d}{ds}[g_{e_k}(s,h(s))] ds  \right)  {<.,e_k>}$  in $(\cS^*)$.

Denote $\Phi_f:= \int^t_{0}  f(s) dB^{(h)}(s) = \int^t_{0}  \left(   \sum^{+ \infty}_{k = 0}  f(s) \tfrac{d}{ds}[g_{e_k}(s,h(s))] {<.,e_k>}  \right)ds $ and define the $(\cS^*)$-process  $\tau:[0,t]\rightarrow (\cS^*)$  by $\tau(s) := \sum^{+ \infty}_{k = 0}  f(s) \tfrac{d}{ds}[g_{e_k}(s,h(s))] \h1 {<.,e_k>}$. Moreover, for $N$ in $\bN^*$, define on $[0,t]$, $\tau_N:s\mapsto \tau_N(s) := \sum^{ N}_{k = 0}  f(s) \tfrac{d}{ds}[g_{e_k}(s,h(s))] {<.,e_k>}$. Obviously we have, in $(\cS^*)$,  $\Phi_f = \int^t_{0} \tau(s) ds$, $\Psi_f = \lim\limits_{N \to +\infty} \int^t_{0}  \tau_N(s)\h1 ds $. It then remains to show that  $ \Phi_f = \lim\limits_{N \to +\infty} \int^t_{0} \tau_N(s) \h1ds$ in $(\cS^*)$. Let us use, for this purpose, theorem \ref{cc}. Let $p_0$ be an integer greater than or equal to $2$. It is easily seen that $\tau_n$  and  $\tau$ are weakly measurable  on $[0,t]$ for every $n$ in $\bN$ (see definition \ref{bb}) and that, $\tau_n(s)$  and  $\tau(s)$   belongs to  $(\cS_{-p_0})$  for every $n$ in $\bN$ and $s$ in $[0,t]$. Moreover, both functions $s \mapsto  {||\tau_n(s)||}_{-p_0}$ and $s \mapsto  {||\tau(s)||}_{-p_0}$ belong to $L^1([0,t],du)$  since ${||\tau_n(s)||}_{-p_0} \leq  {||\tau(s)||}_{-p_0} \leq  |f(s)|     \h1 D  \sqrt{\sum^{+ \infty}_{k = 0} {{(2k+2)}^{-2(p_0 - 1)}}}$ for a certain $D$ given by lemma \ref{zoidpsqojxiqsoxijqsojdzj} $(ii)$. We hence have shown that both functions $\tau_n(.)$  and  $\tau(.)$ are Bochner integrable on $[0,t]$. Besides, for every $(n,m)$ in $\bN^2$ with $n\geq m$, we have

\begin{align*}
&\int^t_{0} \h1 {|| \tau_n(s) - \tau_m(s) ||}_{-p_0} \h1 ds  =  \int^t_{0} \h1 {\left|\left| \sum^{ n}_{k = m+1}  f(s) \tfrac{d}{ds}[g_{e_k}(s,h(s))] {<.,e_k>}\right|\right|}_{-p_0} ds \\
& \leq \underbrace{\int^t_{0} \h1 {\left|\left| \sum^{+ \infty}_{k = m+1}  f(s) \tfrac{d}{ds}[g_{e_k}(s,h(s))] {<.,e_k>}\right|\right|}_{-p_0}  ds }_{ = \int^t_{0} \h1 {|| \tau(s) - \tau_m(s) ||}_{-p_0}  ds } \leq  M { \left( \sum^{+ \infty}_{k = m+1}    \frac{\underset{s \in [0,t] }{\sup}{ {|\tfrac{d}{ds}[g_{e_k}(s,h(s))] |}}^2}{{ {(2k+2)}^{2p_0}}}  \right)       }^{1/2} \notag \\
&\leq M \h1 D  \h1  {\left(\sum^{+ \infty}_{k = m+1} \frac{ 1}{{(2k+2)}^{2(p_0 -1)}}\right)}^{1/2} \xrightarrow[(n,m) \to (+\infty,+\infty)]{ } 0,
\end{align*}

where $M := {||f||}_{L^1([0,t])}$ and $D$ is again given by $(ii)$ of lemma \ref{zoidpsqojxiqsoxijqsojdzj}.  Theorem \ref{cc} then applies and allows to write that $\lim\limits_{N \to +\infty} \int^t_{0} \h1 \tau_N(s) \h1 ds =  \int^t_{0} \h1 \tau(s) \h1 ds$ in $(\cS^*)$. We hence have shown that $\Psi_f = \lim\limits_{N \to +\infty} \int^t_{0} \h1 \tau_N(s) \h1 ds =  \int^t_{0} \tau(s) ds =\Phi_f$ in $(\cS^*)$. This  ends the proof of {\bfseries c)} and establishes formula (\ref{pdivjdofijvdoijfdoij}).

If $ \sum^{+ \infty}_{k = 0}  {\left( \int^t_{0} f(s) \tfrac{d}{ds}[g_{e_k}(s,h(s))] ds  \right)}^2~<~+\infty$, for all $t$, then $Z(t)$ is the
$(L^2)$-limit of a sequence of independent Gaussian variables.
Formula (\ref{pngphopkhghkponpokg}) is then obvious.  When $f$ is of class $\cC^1$ and  such that  $\underset{t \in \bR }{\sup}{|f'(t)|} < +\infty$, an integration by parts yields that  $\displaystyle{\sum^{+ \infty}_{k = 0}  {\left( \int^t_{0} f(s) \tfrac{d}{ds}[g_{e_k}(s,h(s))] ds  \right)}^2 < +\infty}$.
\end{pr}

It is easy to check that definitions \ref{wienerR} and \ref{oezifhherioiheroiuh} coincide on the space $\cE(\bR)$. Indeed for  $f:= \sum^{n}_{k=1} \alpha_k  \i1_{[0,t_k]}$ in $\cE(\bR)$, remark \ref{qzfekdpokfdpof} and equality (\ref{firsozkqsjdefrt}) entail that ${\cJ}^{h}(f) = \sum^{n}_{k=1} \alpha_k B^{(h)}_{t_k}$ almost surely. According to $(i)$ of proposition \ref{zpofkpsokfepsofkspkod}, we have the equality $\int_{\bR} f(s) dB^{(h)}(s) ds = \sum^{n}_{k=1} \alpha_k \int_{\bR} \i1_{[0,t_k]}(s)  W^{(h)}(s) ds  = \sum^{n}_{k=1} \alpha_k B^{(h)}_{t_k}$ almost surely. This implies in particular that ${||\int_{\bR} f(s) dB^{(h)}(s) ds ||}_{(L^2)} = {||f||}_{h}$ for all $f$ in $\cE(\bR)$ since we have ${||{\cJ}^{h}(f)||}_{(L^2)} = {||f||}_{h}$ for such $f$.
Since Wiener integrals with respect to  standard Brownian motion are the elements of the set $  \{ \int_{\bR} f(s) \h1 dB(s), \h1  f  \in     {L^2(\bR)}  \} = \overline{  \{ \int_{\bR} f(s) \h1 dB(s), \h1  f  \in     {\cE(\bR)}     \}}^{(L^2)}$, it seems natural to give the following definition.

\begin{defi}{(Wiener integral with respect to mBm) \\}\label{WieWick}
For an mBm $B^{(h)}$, define the Gaussian space $\Theta_{h}$ by

\begin{equation*}
\Theta_{h}:= \overline{  \{ \int_{\bR} f(s) \h1 dB^{(h)}(s), \h1  f  \in {\cE(\bR)}     \}}^{(L^2)}.
\end{equation*}

We call Wiener integral with respect to $B^{(h)}$ the elements of $\Theta_{h}$.
\end{defi}

\begin{rem}\label{philarmonic}
{\bfseries $(i)$}
Obviously, $\Theta_{h} = \overline{  \{ \sum^{+ \infty}_{k = 0}  \left( \int_{\bR} f(s) \tfrac{d}{ds}[g_{e_k}(s,h(s))] ds  \right)     {<.,e_k>} \hspace{0.15cm}   :\hspace{0.15cm}  f \in \cE(\bR)     \}}^{(L^2)}$.
Thanks to definition $\ref{wienerR}$, theorem $\ref{sdsiduhsdiuhfsdiudfhdsufhisduhfisduhfsiuhfsidhfiuhdfisudhfiu}$  and the fact that ${\cJ}^{h}(f) = \int_{\bR} f(s) \h1 dB^{(h)}(s)$ on $\cE(\bR)$, we have

\begin{equation}
\Theta_{h} = \{ {\cJ}^{h}(u) \h1:\h1 u \in \overline{\cE(\bR)}^{{<,>}_{h}} \}
\end{equation}

In other words, the set of Wiener integrals in the sense of definition $\ref{wienerR}$ and \ref{WieWick} coincide. 

{\bfseries $(ii)$}
When $h$ is a constant function equal to $H$ we find that $\Theta_{h} = \Theta_{H} = \{ <.,M_H(f)> \h1:\h1 f\in   L^2_H(\bR) \}$ since
$\overline{\cE(\bR)}^{{<,>}_{h}} =  \overline{\cE(\bR)}^{{<,>}_{H}} = L^2_H(\bR)$. This is exactly what is expected in view of $(\ref{sixth})$.
\end{rem}

In fact we can be a little more precise in the case of fBm. Let $supp(\cK)$ denote the set of measurable functions $f:\bR\to\bR$ with compact support.

\begin{prop}
Let $H \in (0,1)$. Then:
\bit
\iti Let $f:\bR \mapsto \bR$ be in $L^1_{loc}(\bR) \cap L^2_{H}(\bR)$. Then $\int_{\bR} f(s) \h1 dB^{(H)}(s)$ belongs to $(L^2)$ if and only if\\
$\int_{\bR} f(s) \h1 dB^{(H)}(s) =  {\cJ}^{H}(f)$.
\itii $L^1_{loc}(\bR) \cap L^2_{H}(\bR) \cap supp(\cK) \subset \{f:\bR \to \bR : \h1 \int_{\bR} f(s) \h1 dB^{(H)}(s) \in (L^2)  \}$.
\itiii For $\mu$-almost every $f$ in $L^1_{loc}(\bR) \cap supp(\cK) \cap \{f:\bR \to \bR : \h1 \int_{\bR} f(s) \h1 dB^{(H)}(s) \in (L^2)  \}$, $f$ is in $L^2_H{(\bR)}$ and verifies $\int_{\bR} f(s) \h1 dB^{(H)}(s) =  {\cJ}^{H}(f)$.
\eit
\end{prop}

\begin{pr}
 $(i)$ Let $f \in L^1_{loc}(\bR) \cap L^2_{H}(\bR)$ and define $\Phi_f:=\int_{\bR} f(s) \h1 dB^{(H)}(s)$.  By theorem \ref{sdsiduhsdiuhfsdiudfhdsufhisduhfisduhfsiuhfsidhfiuhdfisudhfiu},
$\Phi_f =  \sum^{+ \infty}_{k = 0} <f,M_H(e_k)><.,e_k>$ where the equality holds in ${(\cS)}^*$. If we assume that $\Phi_f$ belongs to
$(L^2)$, then the equality is valid in $(L^2)$. Besides, since $f$ belongs to $L^1_{loc}(\bR) \cap L^2_{H}(\bR)$ we have, according to theorem \ref{ozdicjdoisoedzijqpzeoejcvenvdsdlsiocfuvfsosfd}, ${\cJ}^{H}(f) = <.,M_H(f)> = \sum^{+ \infty}_{k = 0} {<M_H(f),e_k>}_{L^2(\bR)} <.,e_k> = \sum^{+ \infty}_{k = 0} <f,M_H(e_k)><.,e_k>$ in $(L^2)$. The converse part is obvious since $\int_{\bR} f(s) \h1 dB^{(H)}(s) =  {\cJ}^{H}(f)$ entails that $\int_{\bR} f(s) \h1 dB^{(H)}(s)$ belongs to $(L^2)$.

$(ii)$ Since $f$ is in $L^1_{loc}(\bR)\cap supp(\cK)$ theorem \ref{sdsiduhsdiuhfsdiudfhdsufhisduhfisduhfsiuhfsidhfiuhdfisudhfiu} entails that\\
 $\int_{\bR}  f(s) dB^{(H)}(s) =  \sum^{+ \infty}_{k = 0}  \left( \int_{\bR} f(s) M_H(e_k)(s) ds  \right) {<.,e_k>} = \sum^{+ \infty}_{k = 0}  <f,M_H(e_k)> <.,e_k>$ in ${(\cS)}^*$ (one only needs to replace $\i1_{[0,t]}$ by  $\i1_{supp(f)}$, where $supp(f)$ denotes the support of $f$
in theorem \ref{sdsiduhsdiuhfsdiudfhdsufhisduhfisduhfsiuhfsidhfiuhdfisudhfiu}). Besides, since $f$ belongs to $L^2_{H}(\bR)$, ${\cJ}^{H}(f)$ exists
and is equal to $\sum^{+ \infty}_{k = 0} <f,M_H(e_k)><.,e_k>$ in $(L^2)$.

$(iii)$ Fix $f$ in $L^1_{loc}(\bR) \cap supp(\cK) \cap \{f:\bR \to \bR : \h1 \int_{\bR} f(s) \h1 dB^{(H)}(s) \in (L^2)  \}$ and define $\widetilde{\Omega}$ as subset of $\omega$ in $\Omega$ such that (\ref{firsozeijdtopopkpo}) is true for all $g$ in $\{e_k : \h1 k \in \bN \}$. Clearly $\widetilde{\Omega}$ belongs to $\cG$. As soon as $f$ is in $\widetilde{\Omega}$, we can write $\int_{\bR} f(s) \h1 dB^{(H)}(s) = \sum^{+ \infty}_{k = 0} <M_H(f),e_k><.,e_k>$ in $(L^2)$. This entails that $M_H(f)$ belongs to $L^2(\bR)$ and then, by bijectivity of $M_H$, that $f$ belongs to $L^2_H(\bR)$.
\end{pr}

\begin{rem}
This proposition shows in particular, that for $\mu-$almost every $g$ in $supp(\cK)$:
$$g \in L^2_{H}(\bR) \Leftrightarrow \int_{\bR} g(s) \h1 dB^{(H)}(s) \in (L^2).$$
Moreover, in this case,
$${\cJ}^{H}(g)\overset{(L^2)}{=}\int_{\bR} g(s) \h1 dB^{(H)}(s).$$
\end{rem}
%

\section{Itô Formulas}\label{Ito}
		
\subsection{Itô Formula for generalized functionals of mBm on an interval $[a,b]$ with $0< a < b$}

Let us fix some notations.  For a tempered distribution $G$ and a positive integer $n$, let $G^{(n)}$ denote the $n^{\text{th}}$ distributional derivative of $G$. We also write $G':= G^{(1)}$. Hence, by definition, the equality  $<G',\varphi>=-<G,\varphi'>$ holds for all $ \varphi$ in $\sS(\bR)$.
For a map $t \mapsto F(t)$ from $[a,b]$ to $\sS_{-p}(\bR)$ we will note  $\frac{\partial^n F }{\partial x^n}(t)$ the quantity ${(F(t))}^{(n)}$, that is the $n^{\text{th}}$ derivative in $\sS^{'}(\bR)$, of the tempered distribution $F(t)$.
Hence we may consider the map  $t \mapsto \frac{\partial^n F }{\partial x^n}(t)$ from $[a,b]$ to $\sS'(\bR)$. Moreover for any $t_0$ in $[a,b]$, we will note $\frac{\partial F}{\partial t}(t_0)$ the quantity $\lim\limits_{r \to 0} \frac{F(t_{0}+r)-F(t_0)}{r}$ when it exists in $\sS_{-p}(\bR)$, for a certain integer $p$. When it exists, $\frac{\partial F}{\partial t}(t_0)$ is a tempered distribution, which is said to be the derivative of the distribution $F(t)$ with respect to $t$ at point $t=t_0$.
In line with section \ref{eprofkerpokveprvkerpovepovk}, we then  define, for $t_0$ in $[a,b]$ and a positive integer $n$,  the following quantities:

\begin{equation*}
\frac{\partial^n F }{\partial x^n}(t_0,B^{h}(t_0)) := {(F(t_0))}^{(n)} (B^{h}(t_0)) \hspace{0.5cm} \text{and} \hspace{0.5cm}  \frac{\partial F }{\partial t}(t_0,B^{h}(t_0)) := \left(\frac{\partial F }{\partial t}(t_0)\right) (B^{h}(t_0)).
\end{equation*}

\begin{theo}\label{ziejoezijoez}
Let $p \in \bN$, $a$ and $b$ two reals with $0<a < b$, and let $F$ be an element of $\cC^1([a,b],\sS_{-p}(\bR))$ such that both maps $\frac{\partial F}{\partial x}$ and $\frac{\partial^2 F}{\partial x^2}$, from $[a,b]$ into $\sS_{-p}(\bR)$, are continuous. Then the following equality holds in ${(\cS)}^*$:
\begin{multline}\label{orijvedze2}
F(b,B^{(h)}(b)) - F(a,B^{(h)}(a)) = \int^b_a \h1 \frac{\partial F}{\partial t}(s,B^{(h)}(s)) \h1 ds + \int^b_a \h1 \frac{\partial F}{\partial x}(s,B^{(h)}(s)) \h1 dB^{(h)}(s)\\
								+ \frac{1}{2} \h1 \int^b_a \h1 \left( \frac{d}{ds}[R_{h}(s,s)] \right)\h1 \frac{\partial^2 F}{\partial x^2}(s,B^{(h)}(s)) \h1 ds.
\end{multline}
\end{theo}

\begin{rem}
Recall that for all  $t$ in $[a,b]$, \h1  $ \frac{d}{dt}[R_{h}(t,t)] = 2\h1 t^{2h(t)-1} \h1 (h'(t) \h1 t\ln t + h(t))  $.
\end{rem}

\begin{pr}
We follow closely \cite{ben1} p.$97$-$98$ for this proof. First notice that the three integrals on the right side of (\ref{orijvedze2}) exist since all integrands verify the assumptions of theorem \ref{dlsdqsdhuqsgdudguqsgdqudgquygdsudgqsuydgqsuydgqsuyqdgq}. According to lemma \ref{dkdskcsdckksdksdmksdmlkskdm} it is then sufficient to show equality of the $S$-transforms of both sides of (\ref{orijvedze2}). It is easy to
see that, for every $\eta \in \sS(\bR)$, the function $t \mapsto \gamma \left(t^{2h(t)}, .-{\int^t_0 \h1  M_{h(t)} (\eta)\h1(u) \h1 du}  \right)$ is differentiable from $(0,b]$ into $\sS(\bR)$. Using theorem \ref{tardileomalet} and theorem \ref{peogi} we may write, for $t$ in $[0,b]$:

\begin{align*}
\frac{d}{dt} S(F(t,B^{(h)}(t))(\eta)	&= \frac{d}{dt} \left< F(t), \gamma \left(t^{2h(t)}, .-\displaystyle{\int^t_0 \h1  M_{h(t)} (\eta)\h1(u) \h1 du}  \right) \right>      \notag \\
										&=  \left<\frac{\partial F }{\partial t}(t), \gamma \left(t^{2h(t)}, .-\displaystyle{\int^t_0 \h1 M_{h(t)} (\eta)\h1(u) \h1 du}  \right)   \right> \notag \\
										& \h1+ 2\h1 t^{2h(t)-1} \h1 (h'(t) \h1 t\ln t + h(t)) \h1 \left< F(t), \frac{\partial \gamma}{\partial t} \left(t^{2h(t)}, .-\displaystyle{\int^t_0 \h1  M_{h(t)} (\eta)\h1(u) \h1 du}  \right) \right> \notag \\
										&- \frac{d}{dt}[g_{\eta}(t,h(t))] \h1 \left< F(t), \frac{\partial \gamma}{\partial x} \left(t^{2h(t)}, .-\displaystyle{\int^t_0 \h1  M_{h(t)} (\eta)\h1(u) \h1 du}  \right) \right> =: I_1 + I_2 + I_3.
\end{align*}

Now, $I_1 = S\bigg(  \frac{\partial F}{\partial t}(t)(B^{(h)}(t))  \bigg)(\eta) =  S\bigg( \frac{\partial F}{\partial t}(t,B^{(h)}(t)) \bigg)(\eta)$ using  theorem \ref{tardileomalet} $(iii)$. Besides, since $\gamma$ fulfills the equality $\frac{\partial \gamma}{\partial t} = \frac{1}{2} \frac{\partial^2 \gamma}{\partial x^2}$, we get

\begin{align*}
 I_2 	&=  t^{2h(t)-1} \h1 (h'(t) \h1 t\ln t + h(t)) \h1 \left< F(t), \frac{\partial^2 \gamma}{\partial x^2} \left(t^{2h(t)}, .-\displaystyle{\int^t_0 \h1  M_{h(t)} (\eta)\h1(u) \h1 du}  \right) \right> \notag \\
		& =  \frac{1}{2}\h1\frac{d}{dt}[R_{h}(t,t)]  \h1   S\bigg( \frac{\partial^2 F}{\partial x^2}(t,B^{(h)}(t)) \bigg)(\eta).
\end{align*}

Using $(ii)$ of theorem \ref{tardileomalet},

\begin{equation*}
I_3 = S(W^{(h)}(t))(\eta)  \h1 \left<  \frac{\partial F}{\partial x} (t),\gamma \left(t^{2h(t)}, .-\displaystyle{\int^t_0 \h1  M_{h(t)} (\eta)\h1(u) \h1 du}  \right) \right>. \\
\end{equation*}

Finally, we obtain

\begin{multline*}
\frac{d}{dt} S(F(t,B^{(h)}(t))(\eta	)	= S\bigg( \frac{\partial F}{\partial t}(t,B^{(h)}(t)) \bigg)(\eta) + S\left(W^{(h)}(t) \diamond \frac{\partial F}{\partial x} (t,B^{(h)}(t)) \right)(\eta) \\
+ \frac{1}{2}\h1\frac{d}{dt}[R_{h}(t,t)] \h1  S\bigg( \frac{\partial^2 F}{\partial x^2}(t,B^{(h)}(t)) \bigg)(\eta).
\end{multline*}
\end{pr}

In the proof of theorem \ref{zpzpkpskpsokpdokposdkfpsdokfpsokshfigushfzdiziyuezaiozo} we will need the particular case where the function $F(.)$ is constant, equal to a tempered distribution that we denote $F$. In this case we have the following

\begin{cor}\label{oeoeirjfoerjforefierfroirjeofjerofrijeojizio}
Let $0< a < b$ and $F$ be a tempered distribution. Then the following equality holds in ${(\cS)}^*$:

\begin{equation*}
F(B^{(h)}(b))- F(B^{(h)}(a)) = \int^b_a F'(B^{(h)}(s)) \h1 dB^{(h)}(s) + \frac{1}{2} \h1 \int^b_a \h1 \left( \frac{d}{ds}[R_{h}(s,s)] \right)\h1 F''(B^{(h)}(s)) \h1 ds.
\end{equation*}

\end{cor}

\begin{rem}\label{sdckspdkocpsdkcsdpcksdpkpcokdspcsdkpcodskcpdsk}
Of course when the function $h$ is constant on $\bR$,  we  get the Itô formula for  fractional Brownian motion given in \cite{ben1}.
\end{rem}

\subsection{Itô Formula  in $(L^2)$}
In this subsection, we give two further versions of Itô formula. The first one holds for functions with polynomial growth but weak
differentiability assumptions, whereas the second one deals with $\cC^{1,2}$ functions with sub-exponential growth.

\subsubsection{Itô Formula for certain generalized functionals of mBm on an interval}

Theorem \ref{ziejoezijoez} does not extend immediately to the case $a=0$ because the generalized functional is not defined in this situation, since  $M_{h}(\i1_{[0,t]})$ converges to $0$ a.s and in $L^2(\bR)$ when $t$ tends to $0$ (see theorem-definition \ref{lmqsdlmflkjdg}). As in \cite{ben1}, we now extend the formula to deal with this difficulty.
We will need the following lemma which is a particular case of lemma \ref{podiuhyiuuhiudokdfzdokfdepsoksdpokfsdpokfspokdpfoksdpfskfodkof} below.

\begin{lem}\label{podkspdokdfzdokfdepsoksdpokfsdpokfspokdpfoksdpfskfodkof}
Let $f:\bR\rightarrow \bR$ be a continuous function such that there exists a couple $(C,\lambda)$ in $\bR\times\bR_+$ with $|f(y)| \leq C e^{\lambda y^2}$, for all real $y$. Let $g:\bR\rightarrow \bR^*_+$ be a measurable function such that  $\lim\limits_{t \to 0} g(t) = 0$ and define $L_f$  on $\bR^*_+\times\bR$ by  {$L_f(u,x) := \int_{\bR} \h1 f(y) \gamma(u,x-y) \h1 dy$}. Then $\lim\limits_{(t,x) \to (0^+,x_0)}  \hspace{-0.4cm} L_f(g(t),x) = f(x_0)$, for all real $x_0$.
\end{lem}

\begin{theo}
\label{zpzpkpskpsokpdokposdkfpsdokfpsokshfigushfzdiziyuezaiozo}
Let $F:\bR \rightarrow \bR$ be continuous at $0$ and of polynomial growth. Assume that the first distributional derivative of $F$ is of function type $($defined at the beginning of section $\ref{msldkldkmldklsmmks})$. Then the following equality holds in $(L^2)$:

\begin{equation}\label{ozoisdjcjodsjcdsfjciojciosd}
F(B^{(h)}(b))- F(0) = \int^b_0 F'(B^{(h)}(s)) \h1 dB^{(h)}(s) + \frac{1}{2} \h1 \int^b_0 \h1 \left( \frac{d}{ds}[R_{h}(s,s)] \right)\h1 F''(B^{(h)}(s)) \h1 ds.
\end{equation}
\end{theo}

\begin{pr}
We follow again closely \cite{ben1}.

{\bfseries Step $1$:} $\lim\limits_{t \to 0_+} F(B^{(h)}(t)) =  F(0)$ in $(\cS^*)$.
In order to establish this fact, let us use theorem $8.6$ of \cite{Kuo2}. Since $F$ is of polynomial growth, we may write, thanks to formula $(29)$ of \cite{ben1}, that there exist two reals $C$ and $M$ and a positive integer $m$ such that $ \bE[{F(B^{(h)}(t))}^2] \leq C^2 (1+\frac{(2m)!}{2^m m!}{|t|}^{2mh(t)}) \leq M^2$, for all $t$ in $[0,b]$. Since ${||:e^{<.,\eta>}: ||}_{0}=e^{\frac{1}{2} {|\eta|}^2_{0} }$, Cauchy-Schwarz inequality yields

\begin{align}\label{ozoispokpokpkpkpokokpkpokpokpokpok}
{|S[F(B^{(h)}(t))](\eta)|}	&= |\bE[F(B^{(h)}(t)) \h1 :e^{<.,\eta>}: ]| \leq {||F(B^{(h)}(t))||}_{0} \h1 {||:e^{<.,\eta>}: ||}_{0} \notag \\
							&\leq (1+M) \h1 e^{\frac{1}{2} {|\eta|}^2_{0} }; \hspace{1cm} \text{for all $t$ in $[0,b]$ and $\eta$ in $\sS(\bR)$.}
\end{align}

It then just remains to show that $\lim\limits_{t \to 0+} S[F(B^{(h)}(t))](\eta) = F(0)$. Thanks to theorem \ref{tardileomalet}  and lemma \ref{podkspdokdfzdokfdepsoksdpokfsdpokfspokdpfoksdpfskfodkof}, we get

\begin{align}\label{ozohgffispokpokpkpkpokokpkpokpokpokpok}
\lim\limits_{t \to 0+} S[F(B^{(h)}(t))](\eta)	&= \lim\limits_{t \to 0+}  \int_{\bR} \h1 F(y)   \frac{1}{\sqrt{2\pi t^{2h(t)}}}\exp{\left\{-\frac{1}{2t^{2h(t)}}  {\left(\displaystyle{\int^t_0 \h1 M_{h(t)}(\eta)(u) \h1 du-y}   \right )}^2 \right\}  } \h1 dy \notag \\
												&= \lim\limits_{t \to 0+}  \int_{\bR} \h1 F(y) \h1 \gamma\left(t^{2h(t)},{\int^t_0 \h1 M_{h(t)}(\eta)(u) \h1 du-y}\right)  \h1 dy \notag \\
												&=  \lim\limits_{t \to 0+} L_{F}\left( t^{2h(t)}, \tiny{\int^t_0} \h1 M_{h(t)}(\eta)(u) \h1 du\right) = F(0).
\end{align}

{\bfseries Step $2$:}  $\lim\limits_{a \to 0+} \int^b_a   F'(B^{(h)}(t)) \h1 dB^{(h)}(t) = \int^b_0   F'(B^{(h)}(t)) \h1 dB^{(h)}(t) $ in ${(S)}^*$.

Define $[H_1,H_2]:= [\underset{t \in [0,b] }{\min}{h(t)},\underset{t \in [0,b] }{\max}{h(t)} ]$ and let us prove the two following facts

\bit
\iti There exists a constant $D_1$ which depends only of $F$ such that ${||F'(B^{(h)})(t)||}_{-1} \leq D_1  \h1 \max{\left\{\frac{1}{t^{H_1}}, \frac{1}{t^{H_2}}\right\}}  $ for all $t$ in $(0,b]$.
\eit

Let us first notice that, for all $(x,b,t)$ in $\bR \times \bR^*_+ \times (0,b]$, we have \hspace{0.25cm}

\begin{equation*}
\exp{\{-x^2/{4t^{2h(t)}} \}} \leq e^{-x^2 /4} + \varepsilon(b) \h1 \exp{\{-x^2/{4b^{2 H_2}} \}},
\end{equation*}

where $\varepsilon(b) = 1$ if $b \geq 1$ and $\varepsilon(b) = 0$ if $b < 1$. Note moreover that the function $x\mapsto F'(x) \h1 (e^{-x^2 /4} + \varepsilon(b) \h1 \exp{\{-x^2/{4b^{2 H_2}} \}})$ belongs to $L^1(\bR)$ since $F'$ is of function type and belongs to $\sS'(\bR)$.
Since the operator $A^{-1}$ has a norm operator equal to $1/2$ (see \cite{Kuo2} p.$17$) and using the equality ${|M_{h(t)}(\i1_{[0,t]})|}^2_{0} = t^{2h(t)}$, we get the following upper bound, valid for all $k$ in $\bN$,

\begin{equation*}
\begin{array}{lll}
{|<F', {\xi}_{(t,h(t),k)}>|}^2  = {\left|\int_{\bR} \h1 F'(x) \h1 {\pi}^{1/4} {(k!)}^{1/2} t^{k  h(t)} {\exp}{\{-\frac{x^2}{4t^{2 h(t)}}\}} {e}_{k}{(x/(\sqrt{2} t^{ h(t)}))}\h1 dx \right|}^2\\
\qquad \leq {\pi}^{1/2}  \h1 \underset{u \in \bR }{\sup}{|{e}^2_{k}(u)|} \h1 {\left(\int_{\bR} \h1 |F'(x)| \h1 {\exp}{\{-\frac{x^2}{4t^{2 h(t)}}\}} \h1 dx \right)}^2 \h1 t^{2k  h(t)} \h1 {k!} \ \\
\qquad \leq  \h1 \underbrace{ {\pi}^{1/2} \sup\left\{\h1 \underset{u \in \bR }{\sup}{|{e}^2_{k}(u)|}: k \in \bN\right\} \h1 {\left(\int_{\bR} \h1 |F'(x)| \h1 (e^{-x^2 /4} + \varepsilon(b) \h1 \exp{\{-x^2/{4b^{2 H_2}} \}}) \h1 dx \right)}^2}_{=:D_0}  \h1 t^{2k  h(t)} \h1  {k!}.
\end{array}
\end{equation*}

Using $(ii)$ of remark \ref{laiusdhv} and again the fact that the operator $A^{-1}$ has a norm operator equal to $1/2$ (see $(2)$ p.$17$ of \cite{Kuo2} ) we can write,  for all real $t$ in $(0,b]$, that

\begin{align}\label{ozmlpopioiuiuhjgf}
{||F'(B^{(h)}(t))||}^2_{-1} &= {\left|\left|\frac{1}{\sqrt{2\pi} t^{h(t)}} \sum\limits^{+\infty}_{k=0} {(k!)}^{-1} \h1 t^{-2 k h(t)} \h1<F', {\xi}_{(t,h(t),k)}> I_k\left(   {(\hspace{0.1cm}M_{h(t)}(\i1_{[0,t]})\hspace{0.1cm})}^{\otimes k}   \right)  \right| \right|}^2_{-1} \notag \\
							&=  \frac{1}{2\pi t^{2h(t)}} \sum\limits^{+\infty}_{k=0} {(k!)}^{-1} \h1 t^{-4 k h(t)} \h1 {|<F', {\xi}_{(t,h(t),k)}>|}^2 \h1  \underbrace{ {|{(A^{-1})}^{\otimes k} \left({(\hspace{0.1cm}M_{h(t)}(\i1_{[0,t]})\hspace{0.1cm})}^{\otimes k}  \right)|}^2_0}_{ = {\left| A^{-1} \left(M_{h(t)}(\i1_{[0,t]})\right) \right|\h1}^{2k}_0} \notag \\
							&\leq \frac{{D_0}}{2\pi t^{2h(t)}}  \sum\limits^{+\infty}_{k=0} {(k!)}^{-1} \h1 t^{-4 k h(t)} \h1 t^{2k  h(t)} \h1  {k!} \h1 {\left(\frac{1}{2}\right)}^k  t^{2k  h(t)} \leq  \frac{{D_0}}{\pi} \h1 \max{\left\{\frac{1}{t^{2H_1}}, \frac{1}{t^{2H_2}}\right\}}.
\end{align}

\bit
\itii $\int^b_{0} F'(B^{(h)}(t)) \h1 dB^{(h)}(t) $ exists in ${(\cS)}^*$ and is equal to $\lim\limits_{a \to 0+} \int^b_{a} F'(B^{(h)}(t)) \h1 dB^{(h)}(t)$ in the sense of ${(\cS)}^*$.
\eit

In order to establish the existence of $\int^b_{0} F'(B^{(h)}(t)) \h1 dB^{(h)}(t) $ in ${(\cS)}^*$, let us use  theorem \ref{peodcpdsokcpodfckposkcdpqkoq}. From theorem-definition \ref{wickproducte}, we know that $F'(B^{(h)}(t)) \diamond W^{(h)}(t))$ belongs to ${(\cS)}^*$ for every $t$ in $(0,b]$. Moreover using lemma \ref{dede} we get, for $\eta \in \sS(\bR)$ and $t \in (0,b]$,

\begin{align}\label{ozmlpopqjshxqshkjhsqhjqshdjhdsqhjsdqhdqshdsqioiuiuhjgf}
\left| S(F'(B^{(h)}(t)) \diamond W^{(h)}(t))(\eta) \right| &\leq   {||F'(B^{(h)}(t))||}_{-1} \h1 {||W^{(h)}(t)||}_{-2} \h1  \exp\{{|\eta |}^2_{2} \}   \notag \\
															&\leq \underbrace{ \widehat{K}    \h1 \max{\left\{\frac{1}{t^{H_1}}, \frac{1}{t^{H_2}}\right\}}}_{ =: \cL(t)  } \h1 \exp\{{|\eta |}^2_{2} \},
\end{align}

where we have defined $\widehat{K}:= \frac{{D_0}}{\pi} \underset{t \in [0,b] }{\sup}{ \h1 {|| W^{(h)}(t)||}_{-2}  } $. The function $t\mapsto S(F'(B^{(h)}(t)) \h1 \diamond \h1 W^{(h)}(t))(\eta)$ is  measurable on $[0,b]$ since $S(F'(B^{(h)}(t)) \h1 \diamond \h1 W^{(h)}(t))(\eta) = S(F'(B^{(h)}(t)))(\eta) S(W^{(h)}(t))(\eta)$ using theorems \ref{tardileomalet} and \ref{peogi}. Moreover, since $\cL$ belongs to $L^1([0,b])$,
theorem \ref{peodcpdsokcpodfckposkcdpqkoq} applies and shows that $\int^b_{0} F'(B^{(h)}(t)) \h1 dB^{(h)}(t) $ is in ${(\cS)}^*$.
It then just remains to use theorem $8.6$ in \cite{Kuo2} to show the convergence, in the sense of  ${(\cS)}^*$, of $\int^b_{a} F'(B^{(h)}(t)) \h1 dB^{(h)}(t) $ to $\int^b_{0} F'(B^{(h)}(t)) \h1 dB^{(h)}(t) $ as $a$ tends to $0_+$. Let ${(a_n)}_{n \in \bN}$ be a decreasing sequence of real numbers which tends  to $0$ when $n$ tends to $+\infty$ and $\Psi_n:= \int^b_{0} F'(B^{(h)}(t)) \h1 dB^{(h)}(t) - \int^b_{a_n} F'(B^{(h)}(t)) \h1 dB^{(h)}(t)$. For every $\eta \in \sS(\bR)$ and every $n \in \bN$, $S(\Psi_n)(\eta)= \int^{b}_{0}  \i1_{[0,a_n]}(t) \h1 S(F'(B^{(h)}(t)) \diamond W^{(h)}(t))(\eta) \h1 dt$. Using (\ref{ozmlpopqjshxqshkjhsqhjqshdjhdsqhjsdqhdqshdsqioiuiuhjgf}) and the dominated convergence theorem, it is easy to show that $\lim\limits_{n \to +\infty}S(\Psi_n)(\eta) =0$. Hence  theorem $8.6$ of \cite{Kuo2} applies  and shows that $\lim\limits_{a \to 0+} \int^b_a   F'(B^{(h)}(t)) \h1 dB^{(h)}(t) = \int^b_0   F'(B^{(h)}(t)) \h1 dB^{(h)}(t) $ in ${(S)}^*$.

{\bfseries Step $3$:} Proof of $(\ref{ozoisdjcjodsjcdsfjciojciosd})$

For any real $a$ such that $0<a<b$, we have, thanks to  corollary \ref{oeoeirjfoerjforefierfroirjeofjerofrijeojizio},
\begin{equation*}
F(B^{(h)}(b))- F(B^{(h)}(a)) - \int^b_a F'(B^{(h)}(s)) \h1 dB^{(h)}(s) = \frac{1}{2} \h1 \int^b_a \h1 \left( \frac{d}{ds}[R_{h}(s,s)] \right)\h1 F''(B^{(h)}(s)) \h1 ds.
\end{equation*}

Steps $1$ and $2$ ensure that the left hand side has a limit in ${(\cS)}^*$ when $a$ tends to $0$.
Using theorem \ref{peodcpdsokcpodfckposkcdpqkoq}, it is easy to see that $\int^b_0 \h1 \left( \frac{d}{ds}[R_{h}(s,s)] \right)\h1 F''(B^{(h)}(s)) \h1 ds$ belongs to  ${(\cS)}^*$. Hence using the dominated convergence theorem  and lemma \ref{dkdskcsdckksdksdmksdmlkskdm} we deduce that $\lim\limits_{a \to 0+}  \h1 \int^b_a \h1 \left( \frac{d}{ds}[R_{h}(s,s)] \right)\h1 F''(B^{(h)}(s)) \h1 ds$ is equal to $ \int^b_0 \h1 \left( \frac{d}{ds}[R_{h}(s,s)] \right)\h1 F''(B^{(h)}(s)) \h1 ds$ in   ${(\cS)}^*$. Since we have
proved that, for all $t$ in $[0,b]$, $F(B^{(h)}(t))$ belongs to $(L^2)$, the same holds for the right hand side of $(\ref{ozoisdjcjodsjcdsfjciojciosd})$ and then this equality holds also in $(L^2)$.
\end{pr}

\begin{rem}
As in the case of fBm (see \cite{ben1}), the fact that both sides of the equality $(\ref{ozoisdjcjodsjcdsfjciojciosd})$ are in $(L^2)$ does not imply that every single element of the right hand side is in $(L^2)$. This will be true if, for instance, $F''(B^{(h)}(t))$ belongs to $(L^2)$ and {$\displaystyle{\int^b_0 \h1 \left| \frac{d}{ds}[R_{h}(s,s)] \right|\h1 {||F''(B^{(h)}(s))||}_{0} \h1 ds < +\infty}$}.
\end{rem}

\subsubsection{Itô Formula in $(L^2)$ for $\cC^{1,2}$ functions with sub-exponential growth}

Let us begin with the following lemma:

\begin{lem}\label{podiuhyiuuhiudokdfzdokfdepsoksdpokfsdpokfspokdpfoksdpfskfodkof}
Let $T>0$ and $f:[0,T]\times \bR \rightarrow \bR$ be a continuous function such that there exists  a couple $(C_T,\lambda_T)$ of
$\bR\times\bR^*_+$ such that $\underset{t \in [0,T] }{\max}{|f(t,y)|} \leq C_T e^{\lambda_T y^2}$ for all real $y$. Define  $a>\lambda_T$, $I_a:=(0,\frac{1}{4a})$ and $J_f:\bR \times \bR_+ \times I_a \rightarrow \bR$ by $J_f(x,t,u):= \int_{\bR} \h1 f(t,y) \gamma(u,x-y) \h1 dy$. Then $J_f$ is well defined and moreover $\lim\limits_{(x,t,u) \to (x_0,0^+,0^+)}  \hspace{-0.45cm} J(x,t,u) = f(0,x_0)$.
\end{lem}

\begin{pr}

This is an immediate consequence of theorems $1$ p.$88$ and $2$ p.$89$ of \cite{Wid}.
\end{pr}

Let us now give an Itô formula for functions with subexponential growth.

\begin{theo}\label{ergoigfjerfijoijgfoifjgoidfjgofd}
Let $T >0$  and $h:\bR \rightarrow (0,1)$ be a $\cC^1$ function  such that  $h'$ is bounded on $\bR$. Let $f$ be a $\cC^{1,2}([0,T]\times\bR,\bR)$ function. Furthermore, assume that there are constants $C \geq 0$ and $\lambda < \frac{1}{4 {\underset{t \in [0,T] }{\max}{t^{2h(t)}}}}$ {\vspace{-0.25cm} such that for all $(t,x)$ in $[0,T]\times \bR$},

\begin{equation}\label{orijv}
{\underset{t \in [0,T] }{\max}{\bigg\{ \left|f(t,x)\right|,\left|\frac{\partial f}{\partial t}(t,x)\right|,\left|\frac{\partial f}{\partial x}(t,x)\right|,\left|\frac{\partial^2 f}{\partial x^2}(t,x)\right| \bigg\}} \leq C e^{\lambda x^2}}.
\end{equation}

Then, for all $t$ in $[0,T]$, the following equality  holds in $(L^2)$:

\begin{multline}\label{orijv2}
f(T,B^{(h)}(T)) = f(0,0) + \int^T_0 \h1 \frac{\partial f}{\partial t}(t,B^{(h)}(t)) \h1 dt + \int^T_0 \h1 \frac{\partial f}{\partial x}(t,B^{(h)}(t)) \h1 dB^{(h)}(t) \\
						+ \frac{1}{2} \h1 \int^T_0 \h1 \left( \frac{d}{dt}[R_{h}(t,t)] \right)\h1 \frac{\partial^2 f}{\partial x^2}(t,B^{(h)}(t)) \h1 dt.
\end{multline}

\end{theo}

\begin{pr}
Our proof is similar to the one of theorem 5.3 in \cite{ben2}.
Let $T>0$  and $t \in [0,T]$. Formula (\ref{orijv2}) may be rewritten as

\begin{multline}\label{odzerijv2}
 \int^T_0 \h1 \frac{\partial f}{\partial x}(t,B^{(h)}(t)) \h1 dB^{(h)}(t) = f(T,B^{(h)}(T)) - f(0,0) - \int^T_0 \h1 \frac{\partial f}{\partial t}(t,B^{(h)}(t)) \h1 dt \\
- \frac{1}{2} \h1 \int^T_0 \h1 \left( \frac{d}{dt}[R_{h}(t,t)] \right)\h1 \frac{\partial^2 f}{\partial x^2}(t,B^{(h)}(t)) \h1 dt.
\end{multline}

In order to show that all members on the right hand side of (\ref{odzerijv2}) are in $(L^2)$, let us use theorem \ref{ozeijdoezjdozeijdoizjdoijezoijzdijz}. If $G$ belongs to $\left\{f, \frac{\partial f}{\partial t}, \frac{\partial f}{\partial x}, \frac{\partial^2 f}{\partial x^2}  \right\}$ we may write, thanks to (\ref{orijv}),
\vspace{-5mm}
\begin{equation}\label{sdklsdlksdklsdlksdsdkloiezioefhik}
\bE\left[{G(t,B^{(h)}(t))}^2\right] \leq C^2 \int_{\bR} \h1 \exp\left\{{-\frac{1}{2} \left(\frac{1}{t^{2h(t)}}- 4\lambda\right)v^2}\right\} \h1 \frac{dv}{\sqrt{2\pi t^{2h(t)}}} \leq  \overbrace{\h1 { C^2 } \h1 {({  {1-4\lambda \h1 {\underset{t \in [0,T] }{\max}{t^{2h(t)}}}   }})}^{-1/2}}^{=:M^2}.
\end{equation}

Since $B^{(h)}(t) = B^{H}(t)_{| H= h(t)}$ \textit{a.s}, it is easy to see (in view of \cite{ben2} p.978) that $B^{(h)}(t)$ is a  Gaussian variable with mean equal to $\int^t_0 \h1 M_{h(t)}(\eta)(u) \h1 du = g_{\eta}(t,h(t))$ and variance equal to $t^{2h(t)}$  under the probability $\bQ_{\eta}$ which has been defined in (\ref{peropdockdspocdspocksdpokcpsdokcsdpocksdpockdspockds}). Hence, for every $t$ in $(0,T]$ and $\eta$ in $\sS(\bR)$,

\begin{align}\label{hiuiuiuiuiu}
S(G(t,B^{(h)}(t))(\eta) &= \bE_{\bQ_{\eta}}[G(t,B^{(h)}(t))] = \int_{\bR} \h1  G(\h1 t,u + g_{\eta}(t,h(t)) \h1) \h1 \h1\gamma(t^{2h(t)},u) \h1 du \notag\\
						&=  \int_{\bR} \h1 G\left(t,u \h1 t^{h(t)} +  \int^t_0 \h1 M_{h(t)}(\eta)(x) \h1 dx  \right)      \frac{1}{\sqrt{2\pi}} \h1 e^{-u^2/2} \h1 du.
\end{align}

Using the theorem of continuity under the integral sign, we see that the functions  $t\mapsto S[{G(t,B^{(h)}(t))}](\eta) $ and  $t\mapsto S[G(t,B^{(h)}(t)) \diamond W^{(h)}(t)](\eta) $  are continuous on $[0,T]$. Moreover, in view of (\ref{sdklsdlksdklsdlksdsdkloiezioefhik}),
$t \mapsto {||G(t,B^{(h)}(t))||}_0$ belongs to $L^1([0,T])$. Furthermore, note that

\begin{equation*}
\int^T_0 \h1 {\left|\left| \h1 \left( \frac{d}{dt}[R_{h}(t,t)] \right)\h1 \frac{\partial^2 f}{\partial x^2}(t,B^{(h)}(t)) \h1 \right|\right|}_0 \h1 dt \leq \h1 2 M \h1 \int^T_0 \h1  t^{2h(t)-1} \h1 \left| h'(t) \h1 t\ln t + h(t)\right| \h1 dt < +\infty.
\end{equation*}

Thus, theorem \ref{ozeijdoezjdozeijdoizjdoijezoijzdijz} applies and shows that all members on the right side of (\ref{odzerijv2}) are in $(L^2)$.

Let us now show that $t \mapsto \frac{\partial f}{\partial x}(t,B^{(h)}(t)) \h1 \diamond W^{(h)}(t)$ is  ${(\cS)}^*$-integrable over $[0,T]$.

Reasoning as in the estimate (\ref{ozmlpopqjshxqshkjhsqhjqshdjhdsqhjsdqhdqshdsqioiuiuhjgf}), we note that there exists an integer $q\geq 2$ such that $\frac{\partial f}{\partial x}(t,B^{(h)}(t)) \h1 \diamond W^{(h)}(t)$ belongs  to $({\cS}_{-q})$ for every $t$ in $[0,T]$. Moreover, for every $\eta$ in $\sS(\bR)$ and every $t$ in $(0,T]$ we have, using lemma \ref{dede},

\begin{align*}
\left|S\left(\frac{\partial f}{\partial x}(t,B^{(h)}(t)) \h1 \diamond W^{(h)}(t)\right)(\eta)\right|	&\leq {\big|\big|\frac{\partial f}{\partial x}(t,B^{(h)}(t)) \big|\big|}_{0} \h1 {||W^{(h)}(t)||}_{-2} \h1  \exp\{{|\eta |}^2_{2} \} \notag \\
																										&\leq \left(\underset{t \in [0,T] }{\sup}{{||W^{(h)}(t)||}_{-2}}\right) \h1 {\big|\big|\frac{\partial f}{\partial x}(t,B^{(h)}(t)) \big|\big|}_{0} \h1 \exp\{{|\eta |}^2_{2} \},
\end{align*}


Since $t \mapsto {\big|\big|\frac{\partial f}{\partial x}(t,B^{(h)}(t)) \big|\big|}_{0} $belongs to $L^1([0,T])$,  theorem \ref{peodcpdsokcpodfckposkcdpqkoq} applies on $[0,T]$ to the effect that\\ $\int^T_0 \h1 \frac{\partial f}{\partial x}(t,B^{(h)}(t)) \h1 dB^{(h)}(t)$ belongs to ${(\cS)}^*$. It then just remains to show the following equality for all $t$ in $[0,T]$ and $\eta$ in $\sS(\bR)$:

\begin{multline}\label{lafinestla}
S\left(\int^T_0 \h1 \frac{\partial f}{\partial x}(t,B^{(h)}(t)) \h1 dB^{(h)}(t)\right)(\eta) = S\left(f(T,B^{(h)}(T)) - f(0,0) - \int^T_0 \h1 \frac{\partial f}{\partial t}(t,B^{(h)}(t)) \h1 dt\right)(\eta)\\
- S\left(\frac{1}{2} \h1 \int^T_0 \h1 \left( \frac{d}{dt}[R_{h}(t,t)] \right)\h1 \frac{\partial^2 f}{\partial x^2}(t,B^{(h)}(t)) \h1 dt\right)(\eta).
\end{multline}

Using (\ref{hiuiuiuiuiu}) with $G=f$ and applying the  theorem of differentiation under the integral sign, we get

\begin{equation*}
\begin{array}{l}
\frac{d}{dt} \left[S(f(t,B^{(h)}(t))(\eta)  \right] = \int_{\bR} \h1  \frac{d}{dt} \h1 \left[ f(\h1 t,u + g_{\eta}(t,h(t)) \h1) \h1 \h1\gamma(t^{2h(t)},u) \right] \h1 du \notag \\
\quad = \underbrace{\int_{\bR} \h1 \h1\gamma(t^{2h(t)},u) \h1 \frac{d}{dt} \h1 \left[ f(\h1 t,u + g_{\eta}(t,h(t)) \h1) \right]  \h1  du }_{=:U_1(t)} + \underbrace{\int_{\bR} \h1   f(\h1 t,u + g_{\eta}(t,h(t)) \h1)  \h1 \h1\frac{d}{dt} \h1 \left[ \gamma(t^{2h(t)},u) \right]\h1 du }_{=:U_2(t)}.
\end{array}
\end{equation*}

Now,

\begin{align*}
U_1(t)	&= \int_{\bR} \h1 \h1\gamma(t^{2h(t)},u) \h1 \frac{\partial f}{\partial t}( t,u + g_{\eta}(t,h(t))) \h1 du + \int_{\bR} \h1 \h1\gamma(t^{2h(t)},u) \h1 \frac{\partial f}{\partial x}( t,u + g_{\eta}(t,h(t)))   \h1 \frac{d}{dt} g_{\eta}(t,h(t))   du \notag \\
		&= S\left(\frac{\partial f}{\partial t}\big( t,B^{(h)}(t) \big)\right)(\eta) + S\left(\frac{\partial f}{\partial x}\big(t,B^{(h)}(t))\big)\right)(\eta)\h1 \h1 S\left(W^{(h)}(t)\right)(\eta).
\end{align*}

Besides, using the equality $\frac{\partial \gamma}{\partial t} = \frac{1}{2} \frac{\partial^2 \gamma}{\partial x^2}$   and  an integration by parts, we get

\begin{align*}
U_2(t) = \frac{1}{2} \frac{d}{dt}[t^{2h(t)}] \h1  \h1  \int_{\bR} \h1  \frac{\partial^2 f}{\partial x^2}( t,u + g_{\eta}(t,h(t))) \h1\gamma(t^{2h(t)},u) \h1 du  = \frac{1}{2} \frac{d}{dt}[t^{2h(t)}]  \h1    S\left(\frac{\partial^2 f}{\partial x^2}\big( t,B^{(h)}(t) \big)\right)(\eta).\notag \\
\end{align*}

Hence we obtain, for any  $\varepsilon > 0$, upon integrating $t \mapsto U_1(t) + U_2(t)$ between $\varepsilon$ and $T$,

\begin{multline}\label{pdfofdpokdpffkgdpfkpdogkdpfgkdfpogkp}
S(f(t,B^{(h)}(t))(\eta) - S(f(\varepsilon,B^{(h)}(\varepsilon))(\eta) = \int^{T}_{\varepsilon} \h1   S\left(\frac{\partial f}{\partial t}\big( t,B^{(h)}(t) \big)\right)(\eta)   \h1 dt\\
+\int^{T}_{\varepsilon} \h1  S\left(\frac{\partial f}{\partial x}\big(t,B^{(h)}(t))\diamond W^{(h)}(t)\right)(\eta)   \h1 dt+ \frac{1}{2} \int^{T}_{\varepsilon} \h1 \frac{d}{dt}[t^{2h(t)}]  \h1    S\left(\frac{\partial^2 f}{\partial x^2}\big( t,B^{(h)}(t) \big)\right)(\eta) \h1 dt.
\end{multline}

Let us now  show that $\lim\limits_{\varepsilon \to 0+} S\left(f(\varepsilon,B^{(h)}(\varepsilon) \right)(\eta) = f(0,0) = S\left(f(0,B^{(h)}(0) \right)(\eta) $. For every $\varepsilon>0$, (\ref{hiuiuiuiuiu}) can be rewritten as $S(f(\varepsilon,B^{(h)}(\varepsilon))(\eta) = \int_{\bR} \h1  f(\varepsilon,y) \h1 \gamma(\varepsilon^{2h(\varepsilon)},g_{\eta}(\varepsilon,h(\varepsilon))-y) \h1 dy$.

For a fixed $T>0$, let $\lambda_T, C_T$ be such that (\ref{orijv}) is fulfilled.   
There exists $b>0$ such that $\varepsilon^{2h(\varepsilon)}$ belongs to $I_a$ (defined in lemma \ref{podiuhyiuuhiudokdfzdokfdepsoksdpokfsdpokfspokdpfoksdpfskfodkof}) as soon as $0<\varepsilon<b$. Hence we may write, for any $\varepsilon$ in $(0,b)$, $S(f(\varepsilon,B^{(h)}(\varepsilon))(\eta) = J_f(g_{\eta}(\varepsilon,h(\varepsilon)),\varepsilon,\varepsilon^{2h(\varepsilon)})$. Since $\lim\limits_{\varepsilon \to 0+} \varepsilon^{2h(\varepsilon)}$ and $\lim\limits_{\varepsilon \to 0+} g_{\eta}(\varepsilon,h(\varepsilon))$ are equal to $0$, lemma \ref{podiuhyiuuhiudokdfzdokfdepsoksdpokfsdpokfspokdpfoksdpfskfodkof} applies with $x_0 = 0$, and yields $\lim\limits_{\varepsilon \to 0+} S\left(f(\varepsilon,B^{(h)}(\varepsilon) \right)(\eta) = f(0,0)$.

Let us now establish (\ref{lafinestla}). Thanks to the fact that both the  functions $t\mapsto S[{G(t,B^{(h)}(t))}](\eta) $ and  $t\mapsto S[G(t,B^{(h)}(t)) \diamond W^{(h)}(t)](\eta)$ are continuous on $[0,T]$ and using the dominated convergence  we can take the limit when $\varepsilon$ tends to $0$ on the right hand side of (\ref{pdfofdpokdpffkgdpfkpdogkdpfgkdfpogkp}) and finally get

\begin{multline}\label{pkdsvdfvifdudfofdpokdpffkgdpfkpdogkdpfgkdfpogkp}
S[f(T,B^{(h)}(T)) - f(0,0)](\eta) = S\left(\int^{T}_{0} \h1 \frac{\partial f}{\partial t}\big( t,B^{(h)}(t) \big) \h1 dt \right)(\eta)   + S\left(\int^{T}_{0} \h1 \frac{\partial f}{\partial x}\big( t,B^{(h)}(t) \big) \h1 dB^{(h)}(t)\right) (\eta)\\
+  S\left(\frac{1}{2} \int^{T}_{0} \h1 \frac{d}{dt}[t^{2h(t)}]  \h1   \frac{\partial^2 f}{\partial x^2}\big( t,B^{(h)}(t) \big)\right)(\eta) \h1 dt.
\end{multline}
\end{pr}

\begin{rem}
We observe that if we take expectations on both sides of Itô's formula $(\ref{orijv})$, we get exactly formula $(1)$ of theorem $2.1$ of \cite{Hry}, which is a general weak Itô formula for Gaussian processes, in the particular case where the Gaussian process is chosen to be an mBm.
\end{rem}

\section{Tanaka formula and examples}\label{Tanaka}

In this section we first give a Tanaka formula as a corollary to theorem \ref{zpzpkpskpsokpdokposdkfpsdokfpsokshfigushfzdiziyuezaiozo} with $F:x \mapsto |x-a|$.
We then consider the case of two particular $h$ functions that give noteworthy results.

\subsection{Tanaka formula}
\begin{theo}[Tanaka formula for mBm]
Let $h:\bR \rightarrow (0,1)$  be of class $\cC^1$, $a \in \bR$ and $T>0$. The following equality holds in $(L^2)$:

\begin{equation}\label{tanfor}
|B^{(h)}(T) - a| = |a| + \int^{T}_{0} \h1  \text{sign}\left(B^{(h)}(t) - a\right) \h1 dB^{(h)}(t) +  \int^{T}_{0} \h1  \frac{d}{dt}[R_{h}(t,t)] \h1 \delta_{\{a\}}(B^{(h)}(t)) \h1 dt,
\end{equation}

where the function sign is defined on $\bR$ by $\si(x) :=\i1_{\bR^*_{+}}(x) - \i1_{\bR_{-}}(x)$.
\end{theo}

\begin{pr}
This is a direct application of theorem \ref{zpzpkpskpsokpdokposdkfpsdokfpsokshfigushfzdiziyuezaiozo} with $F:x \mapsto |x-a|$.
\end{pr}

\begin{rem}
That the previous equality holds true in $(L^2)$ does not imply of course that both integrals above are in $(L^2)$. This last result will be established in a forthcoming paper.
\end{rem}

\subsection{Itô formula for functions $h$ such that $\frac{d}{dt}[R_{h}(t,t)]  = 0$}

If $h$ verifies $\frac{d}{dt}[R_{h}(t,t)]  = 0$, then the second order term $\frac{\partial^2 f}{\partial x^2}(t,B^{(h)}(t))$
disappears in Itô formula. The formula is then formally the same one as in ordinary calculus. In this case, \eqref{tanfor} reads:

\begin{equation*}
|B^{(h)}(T) - a| = |a| + \int^{T}_{0} \h1  \text{sign}(B^{(h)}(t) - a) \h1 dB^{(h)}(t).
\end{equation*}

Note that the \gg local time" part disappears in this equality.

Solving the differential equation $\frac{d}{dt}[R_{h}(t,t)]  = 0$ yields two families of functions, denoted ${(h_{1,\lambda})}_{\lambda \in \bR^*_+}$
and ${(h_{2,\lambda})}_{\lambda \in \bR^*_-}$:

\begin{equation*}
\begin{defappli}{ \forall \lambda>0, \hspace{0.5cm} h_{1,\lambda}}{(-\infty, -e^{\lambda}) \cup (e^{\lambda},+\infty) }{(0,1)}{t}
 \frac{\lambda}{\ln |t|},
 \end{defappli} \\
\end{equation*}
and
\begin{equation*}
\begin{defappli}{ \forall \lambda<0, \hspace{0.5cm} h_{2,\lambda}}{ (-e^{\lambda},0) \hspace{0.25cm} \cup \hspace{0.25cm} (0,e^{\lambda}) \hspace{0.25cm}}{(0,1)}{t}
 \frac{\lambda}{\ln |t|}.
 \end{defappli}
\end{equation*}

In order to obtain an mBm defined on a compact interval, we may choose a compact subset of  $(-\infty, -e^{\lambda}) \cup (e^{\lambda},+\infty)$ when $ \lambda>0$ and a  compact subset of $(-e^{\lambda},0) \hspace{0.25cm} \cup \hspace{0.25cm} (0,e^{\lambda})$ when $ \lambda<0$.

\renewcommand{\thefootnote}{\arabic{footnote}}

Figures 1 and 2 display examples of mBm with functions $h_{1}(t) := \frac{1}{\ln t}$ defined on $[e+10^{-3},100]$ and $h_{2}(t) := \frac{-1}{\ln t}$ defined on $[10^{-3},1/e-10^{-2}]$.

\begin{figure}[!ht]
 	\begin{minipage}[c]{.46\linewidth}
	\includegraphics[height=5.5cm]{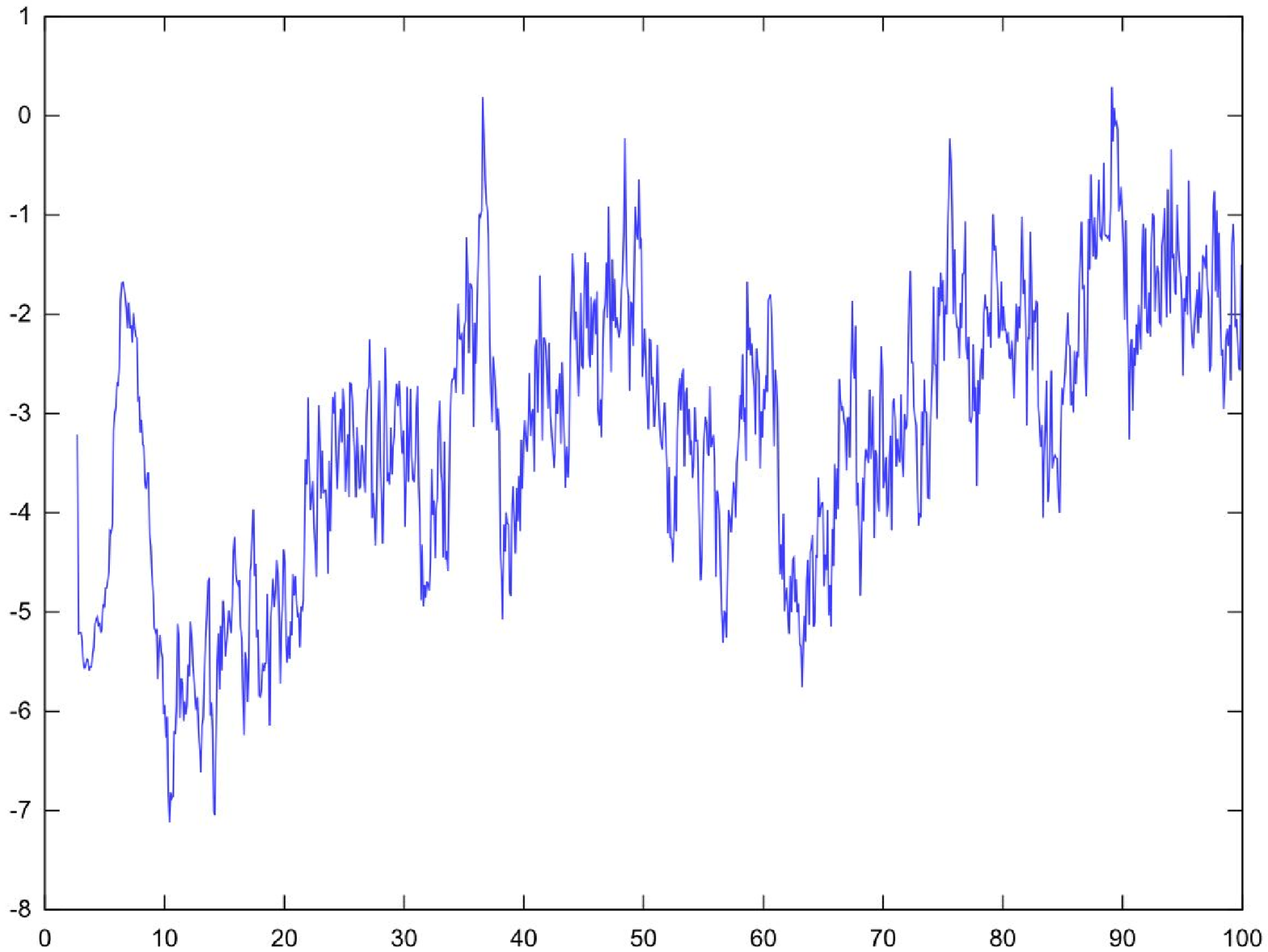}
	\caption{$t \mapsto B^{(h_{1})}(t)$ with $h_1(t) := \frac{1}{\ln t}$ on  $[e+10^{-3},100]$.}
	\end{minipage} \hfill
	\begin{minipage}[c]{.46\linewidth}
	\includegraphics[height=5.5cm]{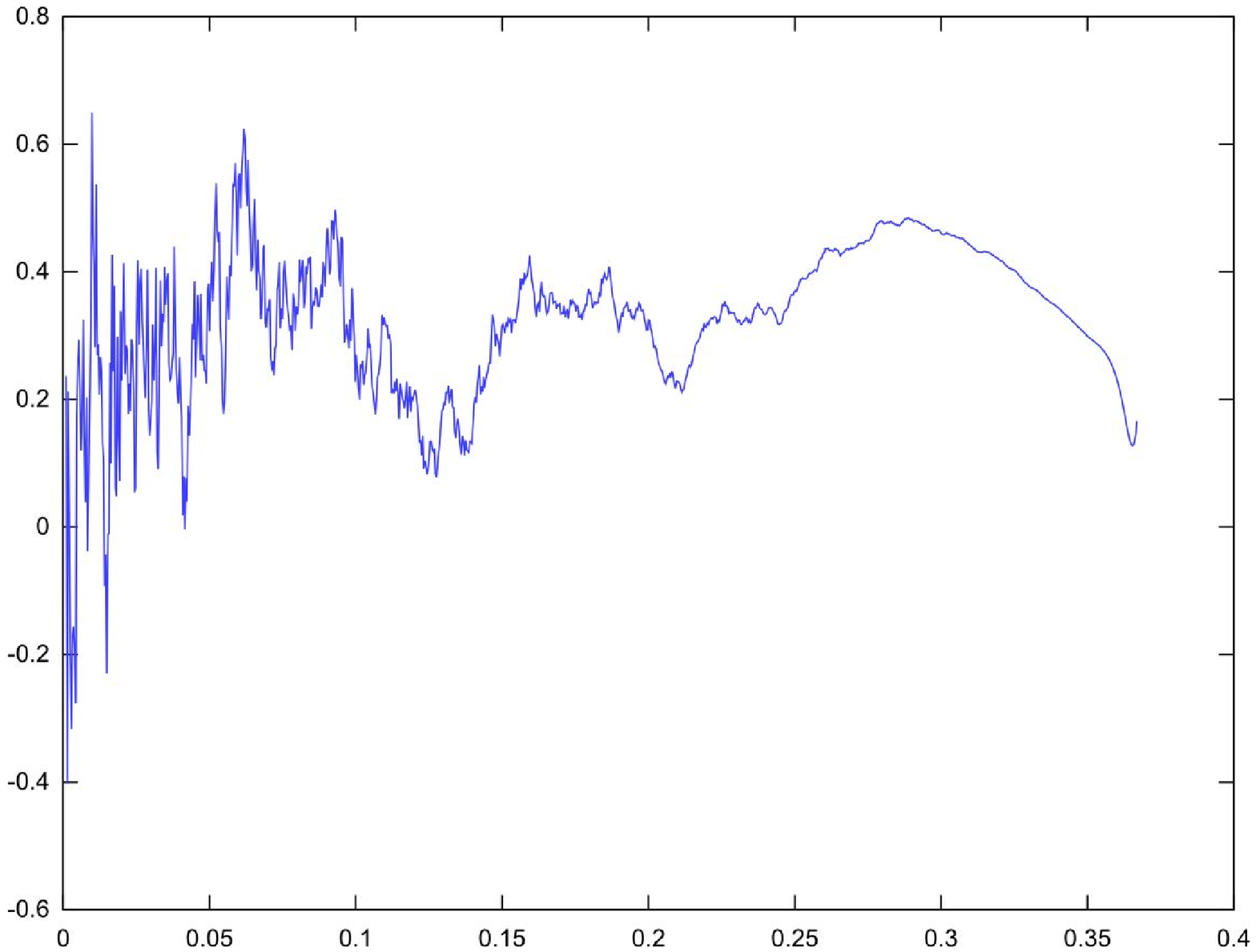}
	\caption{$t \mapsto B^{(h_{2})}(t)$ with $h_2(t) := \frac{-1}{\ln t}$ on  $[1.10^{-3},1/e-10^{-3}]$.}
	\end{minipage}
	\label{fig:figure_1}
\end{figure}

Note moreover that $\lim\limits_{t \to +\infty} h_i(t) = \lim\limits_{t \to 0} h_i(t) =0$ for $i=1,2$.

\subsection{Itô formula for functions $h$ such that $\frac{d}{dt}[R_{h}(t,t)]  = 1$}

The situation where $\frac{d}{dt}[R_{h}(t,t)]  = 1$ is interesting since then Itô formula is formally the same as in the case of standard Brownian motion. As a consequence, Tanaka formula takes the familiar form :

\begin{equation*}
|B^{(h)}(T) - a| = |a| + \int^{T}_{0} \h1  \text{sign}(B^{(h)}(t) - a) \h1 dB^{(h)}(t) +  \int^{T}_{0} \h1  \delta_{\{a\}}(B^{(h)}(t)) \h1 dt.
\end{equation*}

Thus, instead of a "weighted" local time as in \eqref{tanfor}, we get here an explicit expression for the local time of mBm for a family of
$h$ functions that we describe now.

The solutions of the differential equation are given by

\begin{equation}
\begin{defappli}{ h_c}{(c,+\infty) \backslash \{-1,0,1\}}{\bR}{t}
 \frac{1}{2} \frac{\ln(t-c)}{\ln |t|},
 \end{defappli}
\end{equation}

where $c \in \bR$.

Recall that $h_c$ is required to range in $(0,1)$. Denote, for $c \in \bR$, $I_c:=\{ \h1 t \in (c,+\infty) \backslash \{-1,0,1\} \h1 \h1 : \h1 0<h_c(t)<1 \}$. For $c$ in $(-\infty,1/4)$, let $t_1:=t_1(c):=\frac{1-\sqrt{1-4c}}{2}$ and  $t_2:=t_2(c):=\frac{1+\sqrt{1-4c}}{2}$.
Then $I_c$ is explicitly given as follows:

\begin{align*}
&\forall c \in (-\infty,-2], \h1  &I_c = (1+c,t_1)\cup(t_2,+\infty),  \\
&\forall c \in (-2,-1], \h1  &I_c = (t_1,1+c)\cup(t_2,+\infty),  \\
&\forall c \in (-1,0), \h1  &I_c = (t_1,0)\cup (0,1+c)\cup(t_2,+\infty),  \\
&\forall c \in [0,1/4), \h1  &I_c = (t_1,t_2)\cup(1+c,+\infty), \\
&\forall c \geq 1/4, \h1  &I_c = (1+c,+\infty).
\end{align*}

Figures 3 and 4 display examples of multifractional Brownian motion with regularity functions ${h}_{3}(t) := \frac{1}{2} \frac{\ln(t-1)}{\ln t}$ defined on $[2+10^{-3},5]$ and ${h}_{4}(t) := \frac{1}{2} \frac{\ln(t+1)}{\ln |t|}$ defined on $[\frac{1-\sqrt{5}}{2}+10^{-3},-10^{-3}]$.

\begin{figure}[!ht]
 	\begin{minipage}[c]{.46\linewidth}
	\includegraphics[height=5.5cm]{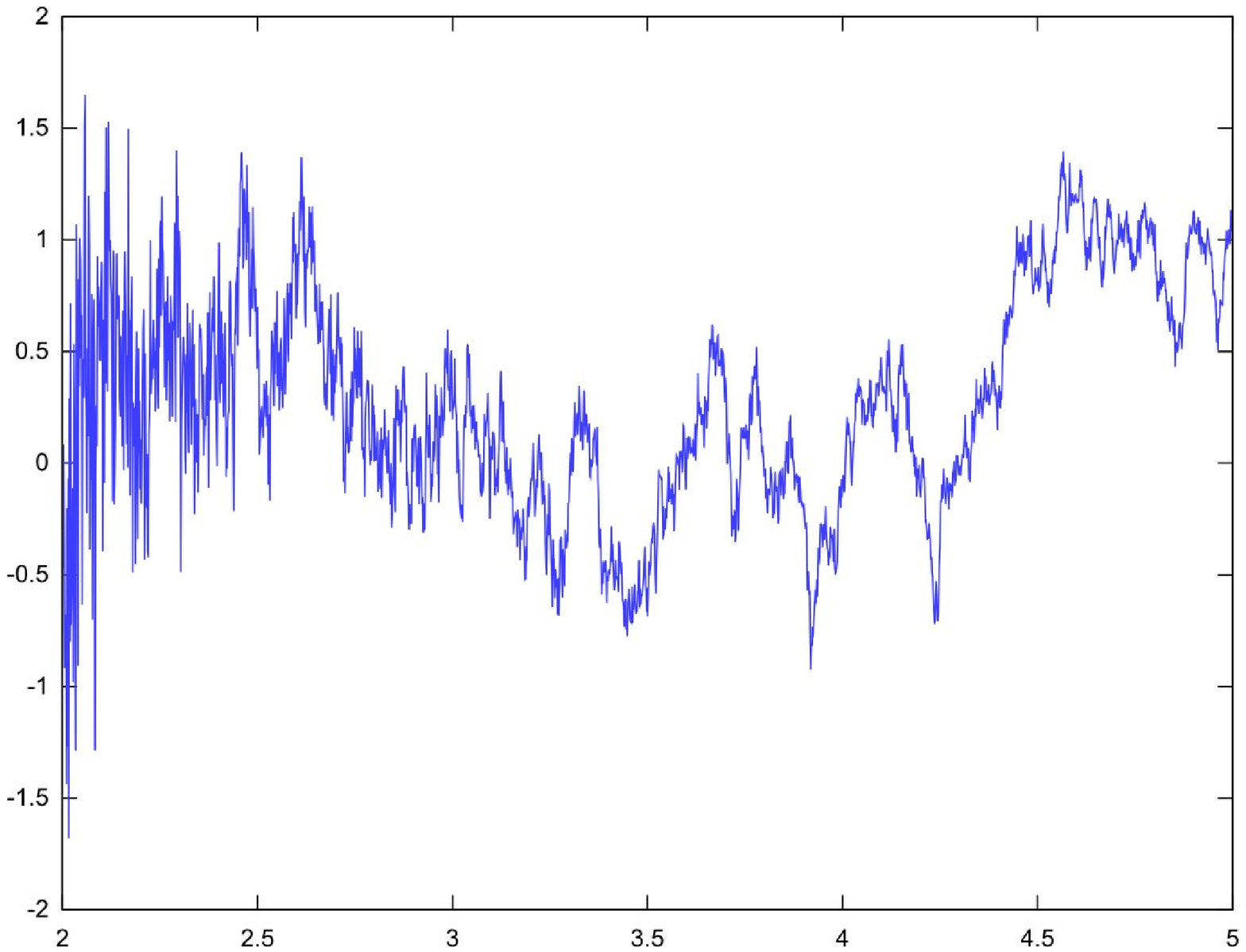}
	\caption{$t \mapsto B^{(h_{3})}(t)$ with $h_3(t) := \frac{1}{2} \frac{\ln(t-1)}{\ln t}$ on  $[2+10^{-3},5]$.}
	\end{minipage} \hfill
	\begin{minipage}[c]{.46\linewidth}
	\includegraphics[height=5.5cm]{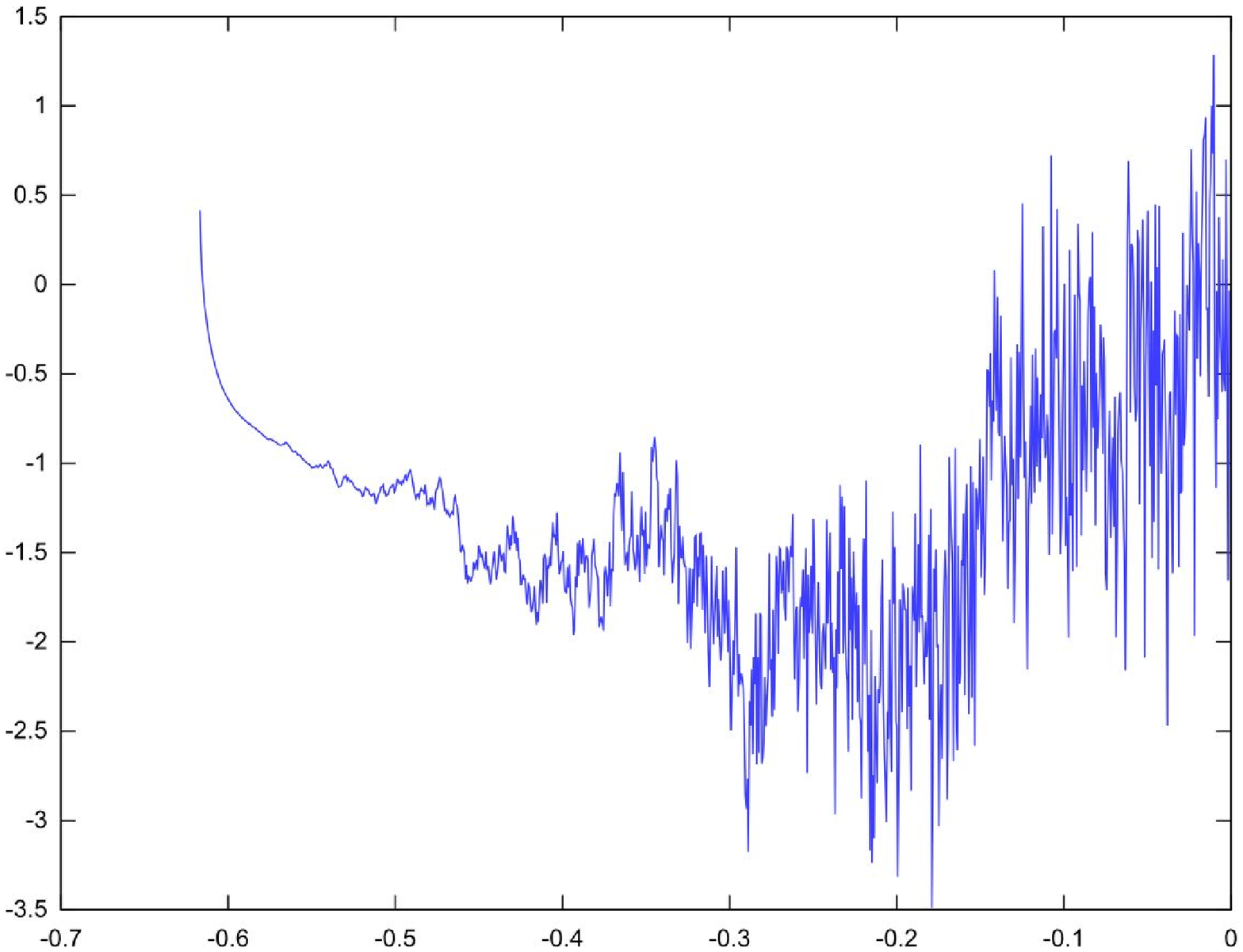}
	\caption{$t \mapsto B^{(h_{4})}(t)$ with $h_4(t) := \frac{1}{2} \frac{\ln(t+1)}{\ln |t|}$  on  $[\frac{1-\sqrt{5}}{2}+10^{-3},-10^{-3}]$.}
	\end{minipage}
	\label{fig:fig:figure_2}
\end{figure}

Note that the case $c = 0$ yields the constant function $h_c \equiv 1/2$, {\it i.e.} standard Brownian motion. Moreover, since $\lim\limits_{t \to +\infty} h_{c}(t) = 1/2$ for every $c$, we see that the family of functions ${h}_{c}$ behaves,	asymptotically, like the constant function equal to $1/2$. However this does not mean that there is convergence in law of $B^{(h_c)}$ to Brownian motion; in fact one needs to scale $B^{(h_c)}$.
For every $t$ in $\bR^*$, define the process $X_t$ on $\bR_+$ by $X_t(u):=\frac{B^{(h_c)}(tu)}{\sqrt{t}}$. Then $\{X_t(u);\h1 u \in \bR_+\} \xrightarrow[t \to +\infty]{\cL} \{B(u); \h1 u\in \bR_+\}$ where $B$ still denotes a Brownian motion and $\cL$ denotes convergence in law.

\section{Conclusion and future work}

In this paper we have used a white noise approach to define a stochastic integral with respect to multifractional Brownian motion which generalizes the one for fBm based on the same approach. This stochastic calculus  allows  to solve some particular stochastic differential equations. We are currently investigating several extensions of this work. In order to apply this calculus to financial mathematics or to physics, it is necessary to study further the theory of stochastic differential equations driven by a mBm. This is the topic of the forthcoming paper \cite{JLJLV2}.

The Tanaka formula we have obtained suggests that one can get several integral representations of local time with respect to mBm. Finally, since mBm is a Gaussian process, it seems also natural to investigate the links between the construction of stochastic integral with respect to mBm that we gave and the one provided by Malliavin calculus. An extension in higher dimension is also desirable.

\section{Acknowledgments}
The authors are thankful to Marc Yor for many helpful remarks and comments. They are also grateful to Erick Herbin for suggesting this research topic and to Sylvain Corlay for the simulations of mBm.


%
%
%
%
%
%
%
%


\appendix

\vspace{1cm}

\cl{ \bfseries \Large Appendix }


\makeatletter
\renewcommand\theequation{\thesection.\arabic{equation}}
\@addtoreset{equation}{section}
\makeatother

\section{Bochner integral}

\label{appendiceB}

All the following notions about the integral in the Bochner sense come from \cite{HP} p.$72$, $80$ and $82$ and from \cite{Kuo2} p.$247$.

%

\begin{defi}[Bochner integral \cite{Kuo2} p.247]
\label{bb}
Let $I$ be a subset of $\bR$ endowed with the Lebesgue measure. One says that  $\Phi:I \rightarrow {(\cS)}^{*}$ is Bochner integrable on $I$ if it satisfies the two following conditions:
\ben
\item $\Phi$ is weakly measurable on $I$ \textit{i.e} $u\mapsto <\hspace{-0.2cm}<\hspace{-0.1cm} \h1 \Phi(u),\varphi\h1 \hspace{-0.1cm}>\hspace{-0.2cm}>$  is measurable on $I$ for every $\varphi$ in $(\cS)$.
\item There exists $p \in \bN$ such that $\Phi(u) \in ({\cS}_{-p})$ for almost every $u \in I$ and $u \mapsto {||\Phi(u)||}_{-p}$ belongs to $L^1(I)$.
\een

The Bochner-integral of $\Phi$ on $I$ is denoted $\int_{I} \Phi(s)  \h1 ds $ .
\end{defi}

\begin{propers}
If $\Phi:I \rightarrow {(\cS)}^{*}$ is Bochner-integrable on $I$  then
\ben
\item there exists an integer $p$ such that

\begin{equation}
{\left|\left| \int_{I} \Phi(s)  \h1 ds \right|\right|}_{-p} \leq \int_{I}  {\left|\left|\Phi(s)\right|\right|}_{-p} \h1 ds.
\end{equation}

\item $\Phi$ is also Pettis-integrable on $I$ and both integrals coincide on $I$.
\een
\end{propers}

\begin{rem}
The previous property shows that there is no risk of confusion by using the same notation for both
the  Bochner integral and the Pettis integral.
\end{rem}

\begin{theo}
\label{cc}

Let $p \in \bN$ and ${(\Phi_n)}_{n \in \bN}$ be a sequence of  processes from $I$ to ${(\cS)}^{*}$ such that  $\Phi_n(u) \in ({\cS}_{-p})$ for almost every $u \in I$ and for every $n$. Assume moreover that $\Phi_n$ is  Bochner-integrable on $I$ for every $n$ and that

\begin{equation*}
\lim\limits_{(n,m) \to (+\infty,+\infty)}  \int_{I}  {\left|\left|\Phi_m(s) - \Phi_n(s)\right|\right|}_{-p} \h1 ds = 0.
\end{equation*}

Then there exists an ${(\cS)}^*$-process (almost surely $({\cS}_{-p})$-valued), denoted $\Phi$, defined and Bochner-integrable on $I$ such that

\begin{equation}\label{iuerhferuhfrufhu}
\lim\limits_{n \to +\infty}  \int_{I}  {\left|\left|\Phi(s) - \Phi_n(s)\right|\right|}_{-p} \h1 ds = 0.
\end{equation}

Furthermore, if there exists an ${(\cS)}^{*}$-process, denoted $\Psi$, which verifies $(\ref{iuerhferuhfrufhu})$, then $\Psi(s) = \Phi(s)$ for almost every $s$ in $I$. Finally the following
equality holds:

\begin{equation*}
\lim\limits_{n \to +\infty}  \int_{I} \Phi_n(s) \h1 ds =  \int_{I} \Phi(s) \h1 ds \hspace{0.25cm} \text{in} \h1(\cS^*).
\end{equation*}
\end{theo}

%

\section{ Proof of proposition \ref{rfjiijofsdqoijqfsdoijqdsqsdoufgyq}}

\label{appenB}

Let $B^{(h)}$ be a normalized mBm on $\bR$ with  covariance function noted $R_h$. It is well known that one can define on the linear space  $\text{span}_{\bR}\{R_{h}(t,.) : t\in \bR \}$ an inner product, denoted ${<,>}_{R_h}$, by  ${<R_h(t,.),R_h(s,.)>}_{R_h} := R_h(t,s)$ (see \cite{janson} p.120 ff.). Define $\Xi_h$ the closure of $\text{span}_{\bR}\{R_{h}(t,.) : t\in \bR \}$ for the norm  ${||\h1||}_{R_h}$.  The space $(\Xi_h,{||\h1||}_{R_h} )$ is called the Cameron-Martin space (or Reproducing Kernel Hilbert Space (R.K.H.S.)) associated to the the Gaussian process $B^{(h)}$.
Let $\widetilde{\cE}(\bR)$ denote the quotient space obtained by identifying all functions of $\cE(\bR)$ which are equal almost everywhere.
On $\widetilde{\cE}(\bR) \times \widetilde{\cE}(\bR)$ define a bilinear form, noted ${< , >}_{h}$, by ${<\i1_{[0,t]},\i1_{[0,s]} >}_{h} := R_{h}(t,s)$. Then ${< , >}_{h}$ is an inner product provided the linear map $\kappa_h:\widetilde{\cE}(\bR) \rightarrow \Xi_h$ defined by $\kappa_h(\i1_{[0,t]}):=R_h(t,.)$, $t \in \bR$, is injective. Define  $\cI_{h}:= \overline{  \text{vect}_{\bR}\{ B^{(h)}(t) : t \in \bR \}}^{(L^2)}$ the first Wiener chaos of $B^{(h)}$. It is a well-known property of R.K.H.S. that the map $\tau_h: (\Xi_h,{||\h1||}_{R_h} ) \rightarrow ( \h1 \overline{  \text{span}_{\bR}\{ B^{(h)}(t) : t \in \bR \}}^{(L^2)}, {||\h1 ||}_{(L^2)} )$, defined for all real $t$ by $\tau_h(R_h(t,.)) = B^{(h)}(t)$ is an isometry. As a result, $\kappa_h$ is injective if and only if $\tau_h \circ \kappa_h$ is injective. The next proposition states that this is indeed the case for any continuous function $h$:

\begin{prop}\label{lmqlqlmqmlqkllkdhueiuzegfyzegfyuz}
Let $h$ be a continuous function defined on $\bR$ and ranging in $(0,1)$. The family ${(B^{(h)}(t))}_{t \in \bR^*}$ is linearly independent on $\bR$, \textit{i.e} for every positive integer $n$, $(\beta_1,\beta_2\cdots,\beta_n)$ in $\bR^n$ and $(t_1,t_2\cdots,t_n)$ in ${(\bR^*)}^n$, such that $t_i\neq t_j$ for $i \neq j$,   the equality

\begin{equation}\label{onyest}
 \sum^{n}_{j = 1} \beta_j \h1 B^{(h)}(t_j) = 0 \hspace{0.25cm} \textit{a.s},
\end{equation}

implies $\beta_1 = \beta_2 = \cdots = \beta_n = 0$.

\end{prop}

The proof of this proposition requires the following lemma, the proof of which is easy and left to the reader.

\begin{lem}\label{dkqskqmlsddqsmklqdskq}

Define, for $t \in \bR$, the function $A_t:\bR \rightarrow \bC$  by $A_t(\xi):= \frac{e^{it\xi}-1}{i\xi}$ if $\xi \neq 0$ and  $A_t(0):= t$. Then,  for all $t$, $A_t$ is $\cC^{\infty}$ on $\bR$  and verifies, for every $n \in \bN$,

\begin{align}\label{jsjsjsjsjsjsjsjsjsjsjsjsj}
A^{(n)}_t(0) = \frac{t \h1 {(it)}^{n}}{n+1}.
\end{align}

where $A^{(n)}_t $ denotes the $n^{\text{th}}$ derivative of $A_t$.
\end{lem}


\begin{pr}{ \bfseries  of proposition \ref{lmqlqlmqmlqkllkdhueiuzegfyzegfyuz}.}
Let us use a proof by contradiction. By decreasing $n$ if necessary we may always assume that $(\beta_1,\beta_2\cdots,\beta_n)$ belongs to ${(\bR^*)}^n$. Besides, thanks to \eqref{fqkiuhgiuhihj} and lemma \ref{osidjoifjsoj}  $(i)$, equality \eqref{onyest} also reads $<. , \h1 \sum^{n}_{j = 1} \beta_j \h1 M_{h(t_j)}(\i1_{[0,t_j]}) > = 0$ \h1 \textit{a.s}. By taking Fourier transforms, we get $\sum^{n}_{j = 1} \beta_j \h1 \widehat{  M_{h(t_j)}(\i1_{[0,t_j]})} = 0$ \h1 \textit{a.e.}. Using \eqref{qxs} this yields

\begin{equation}\label{okdkdzzkpodzkpozkopddz}
\sum^{n}_{j = 1} \alpha_j \h1 {|\xi|}^{1/2-h(t_j)} \h1 \widehat{\i1_{[0,t_j]}}(\xi) = 0, \hspace{0.25cm} \forall \xi \in \bR^*,
\end{equation}

where we have defined, for $j$ in $\{1;2;\cdots;n\}$, $\alpha_j:={\beta_j} \h1 {(c_{h(t_j)})}^{-1}$. By re-arranging if necessary the $(t_i)_i$, we
  may assume without loss of generality that $h(t_1) \geq h(t_2) \geq \cdots \geq h(t_n)$. Let $card(A)$ denote the cardinal of the set $A$. We  distinguish three cases.

\medskip

{\bfseries First case:} card($\{h(t_1);h(t_2);\cdots;h(t_n)\}$) = 1.

Since $h(t_1) = h(t_2) = \cdots = h(t_n) =:H $, we get, by multiplying equality \eqref{okdkdzzkpodzkpozkopddz} by ${|\xi|}^{H-1/2}$ and taking inverse Fourier transform,  $\sum^{n}_{j = 1} \alpha_j \h1 \i1_{[0,t_j]} = 0$ \h1 almost everywhere on $\bR$. This entails that $\{\alpha_1;\alpha_2;\cdots;\alpha_n\}$ and then $\{\beta_1;\beta_2;\cdots;\beta_n\}$ is equal to $\{0\}$.

\medskip

{\bfseries Second case:} $h(t_1) > h(t_2)$.
Using that $\widehat{\i1_{[0,t]}}(\xi) = A_t(\xi)\*$, \eqref{okdkdzzkpodzkpozkopddz} reads:

\begin{equation}\label{okoerorieojfrroedkdzzkpodzkpozkopddz}
\alpha_1 \h1 (\tfrac{e^{it_1 \xi}-1}{i\xi}) = - \sum^{n}_{j = 2} \alpha_j \h1 {|\xi|}^{h(t_1)-h(t_j)} \h1 (\tfrac{e^{it_j \xi}-1}{i\xi}), \hspace{0.25cm} \forall \xi \in \bR^*.
\end{equation}

By lemma \ref{dkqskqmlsddqsmklqdskq} and  taking the limit when $\xi$ tends to $0$ in \eqref{okoerorieojfrroedkdzzkpodzkpozkopddz},   we get $\alpha_1 = 0$ which constitutes a contradiction.

\medskip

{\bfseries Third case: $h(t_1) = h(t_2)$. }

There exists an integer $r$ in $\{2;3;\cdots;n-1\}$,  $(k_1,k_2,\cdots,k_r)$ in ${(\bN^*)}^r$ with $2 \leq k_1 < k_2 <\cdots <k_r = n$, such that

\vspace{-0.6cm}

\begin{equation}\label{onyesoioefjierjit}
\begin{array}{lll}
h(t_1) = h(t_2) = \cdots = h(t_{k_1})							&=:	&H_1\\
h(t_{k_1 + 1}) = h(t_{k_1 + 2}) = \cdots = h(t_{k_2}) 			&=:	&H_2\\
\quad \vdots													&	&\vdots\\
h(t_{k_{r-1} + 1}) = h(t_{k_{r-1} + 2})  = \cdots = h(t_{k_r})	&=:	&H_r,
\end{array}
\end{equation}

where $1 > H_1 > H_2 > \cdots > H_r > 0$. Note that when $r=1$ we have $card(\{h(t_1);\cdots;h(t_n)\}) = 1$ (treated in the first case) and when $r=n$ we have $card(\{h(t_1);\cdots;h(t_n)\}) = n$  and then $h(t_1) > h(t_2)$ (treated in the second case). We hence assume from now that
$2 \leq k_1 \leq n-1$ and $2 \leq r \leq n-1$.

Define the sets $I_1,I_2,\cdots,I_r$ by $I_1:=\{1;2;\cdots;k_1\}$, $I_2:=\{k_1 +1;k_1 +2;\cdots;k_2\}, \cdots, \h1 I_r:=\{k_{r-1} +1; k_{r-1} +2 ;\cdots; k_r\}$. Using lemma \ref {jsjsjsjsjsjsjsjsjsjsjsjsj}, equality  \eqref{okdkdzzkpodzkpozkopddz} can be rewritten  as

\begin{equation}\label{okdkdzzkpodzkpozkopddziuhze}
\sum^{r}_{l = 1} \h1 \sum_{j \in I_l }  \h1 \alpha_j \h1 {|\xi|}^{1/2-H_l} \h1 A_{t_j}(\xi) = 0, \hspace{0.5cm} \forall \xi \in \bR^*.
\end{equation}

\begin{lem}\label{ihuziehuuhdzuhuedz}
For every $p$ in $\bN^*$,
\vspace{-0.25cm}

\begin{equation}\label{okdkdzzkpodzkpozkopddziuhzeuyuyft}\tag{$L_p$}
\sum_{j \in I_1 }  \h1 \alpha_j \h1 t^p_j = 0.
\end{equation}

\end{lem}

Let us admit this lemma for the moment. The equalities~\eqref{okdkdzzkpodzkpozkopddziuhzeuyuyft} for $p$ in $\{1;2; \cdots;k_1\}$ yield the following linear system:

\begin{equation}
\underbrace{
\begin{pmatrix}
t_1 & t_2 & \cdots & t_{k_1-1} &t_{k_1} \\
t^2_1 & t^2_2 & \cdots & t^2_{k_1-1} &t^2_{k_1} \\
\vdots & \vdots & \cdots & \vdots &\vdots\\
t^{k_1}_1 & t^{k_1}_2 & \cdots & t^{k_1}_{k_1-1} &t^{k_1}_{k_1}
\end{pmatrix}
}_{=:D}
\begin{pmatrix}
\alpha_1 \\
\alpha_2 \\
\vdots \\
\alpha_{k_1}
\end{pmatrix}
 =
\begin{pmatrix}
0\\
0 \\
\vdots \\
0
\end{pmatrix}
.
\end{equation}

The determinant of this system is a Vandermonde determinant which is non zero since all the $t_i$ are  distinct from each other. As the result, the
only solution is $\alpha_1 = \alpha_2 = \cdots \alpha_{k_1} = 0$ which constitutes a contradiction and proves the proposition.
\end{pr}
%
%
%

We now present a sketch of proof of lemma \ref{ihuziehuuhdzuhuedz}.

 \begin{praa}{\bfseries of lemma \ref{ihuziehuuhdzuhuedz}.}
By multiplying both sides of equality  \eqref{okdkdzzkpodzkpozkopddziuhze} by ${|\xi|}^{H_1-1/2}$ and then taking the limit when $\xi$ tends to $0$, we get, using lemma \ref{jsjsjsjsjsjsjsjsjsjsjsjsj}, $\sum_{j \in I_1 }  \h1 \alpha_j \h1 t_j = 0$, which is equality $(L_1)$.
Now, fix $p$ in $\bN^*$. Starting from equality  \eqref{okdkdzzkpodzkpozkopddziuhze} we first

- multiply both sides of equality  \eqref{okdkdzzkpodzkpozkopddziuhze} by  ${|\xi|}^{H_1-1/2}$ and call  (\ref{okdkdzzkpodzkpozkopddziuhze} bis)
the resulting equality. For any $\xi$ in $\bR^*$, we then take the $p^{\text{th}}$ derivative of both
sides of equality  (\ref{okdkdzzkpodzkpozkopddziuhze} bis)  at point $\xi$. We call $(E_{1})$ the equality thus obtained. It reads:

\begin{equation}
\sum^r_{l = 1 } \sum_{j \in I_l } \h1 \alpha_j \h1  {[ {|\xi|}^{H_1-H_l}  \h1 A_{t_j}(\xi)   ]}^{(p)}  \h1 = 0,  \hspace{0.5cm} \forall \xi \in \bR^*,
\label{jjsjlkqsjlksqjklqsjklsjklqsqjkl} \tag{$E_{1}$}
\end{equation}

where ${[ g(\xi) ]}^{(p)}$ denotes the $p^{\text{th}}$ derivative of the $p$-times differentiable map $\xi \mapsto g(\xi)$.

Now, starting from \eqref{jjsjlkqsjlksqjklqsjklsjklqsqjkl}, we recursively perform the following operations successively for $l = 2, \ldots, r$:

- multiply both sides of equality $(E_{l-1})$ by ${|\xi|}^{H_l-H_{l-1} + p}$  and call $(E_{l-1} \h1 \text{bis})$ the resulting equality.

- take the $p^{\text{th}}$ derivative of both sides of equality $(E_{l-1} \h1 \text{bis})$  at every point $\xi$ in $\bR^*$ and call $(E_{l})$ the resulting equality.

Equality $(E_r)$ then reads:

\begin{multline}
\sum^r_{l = 1 } \sum_{j \in I_l } \h1 \alpha_j \h1 \times   \\
\underbrace{{\left[
{\left[
\cdots
{\left[
{\left[
 {[ {|\xi|}^{H_1-H_l}  .\h1 A_{t_j}(\xi)   ]}^{(p)} \h1 . {|\xi|}^{H_2-H_1+p}
\right]}^{(p)}  \h1 . {|\xi|}^{H_3-H_2+p} \right]}^{(p)} \cdots
{|\xi|}^{H_{r-1}-H_{r-2}+p}
\right]}^{(p)} . {|\xi|}^{H_r-H_{r-1}+p}
\right]}^{(p)}}_{=:{[\cdots]}_{l,j}^{(p)}(\xi)} \\
  \h1 = 0, \hspace{0.5cm} \forall \xi \in \bR^*,
  \label{jjsjlkqsjlksqjklqsjklsjklqsqjklrrrrr} \tag{$E_{r}$}
\end{multline}

Lemma \ref{jsjsjsjsjsjsjsjsjsjsjsjsj} yields that $ \lim\limits_{\xi \to 0} \alpha_j \h1 A^{(p)}_{t_j}(\xi) =  \frac{i^p}{p+1} \h1 \alpha_j \h1 t^{p+1}_j$. We want to let $\xi$ tend to $0$ in the previous equality. However, for $(l,j)$ in $\{1;2;\cdots;r\} \times I_l$, $ \lim\limits_{\xi \to 0} {[\cdots]}_{l,j}^{(p)}(\xi) = +\infty$. Nevertheless, it is easy to show that, for every $(l,j)$ in $\{1;2;\cdots;r\} \times I_l$,
$\lim\limits_{\xi \to 0} \frac{{[\cdots]}_{l,j}^{(p)}(\xi)}{   A^{(p)}_{t_j}(\xi) \h1 {|\xi|}^{H_r-H_l} }$ exists in $\bC^*$ and is independent of $j$. Define $c_l:=  \lim\limits_{\xi \to 0} \frac{{[\cdots]}_{l,j}^{(p)}(\xi)}{   A^{(p)}_{t_{k_l}}(\xi) \h1 {|\xi|}^{H_r-H_l} }$, $l=1, \ldots, r$. Denote, for $(l,j)$ in $\{1;2;\cdots;r\} \times I_l$, $U_{l,j}: \bR \rightarrow \bC$, the continuous map on $\bR$ such that
${[\cdots]}_{l,j}^{(p)}(\xi) = c_l  \h1 {|\xi|}^{H_r-H_l}  \h1 A^{(p)}_{t_j}(\xi) \h1 (1 + \h1 U_{l,j}(\xi))$. Equality $(E_{r})$ then reads

\begin{equation}\label{qsdqsduiuqsdiusqdfhuqdsh}
\sum^r_{l = 1 }   c_l \h1 \sum_{j \in I_l } \h1 \alpha_j  \h1 {|\xi|}^{H_r-H_l}  \h1 A^{(p)}_{t_j}(\xi) \h1 (1 + \h1 U_{l,j}(\xi)) = 0, \hspace{0.5cm} \forall \xi \in \bR^*.
\end{equation}

Upon multiplying both sides of the previous equality by ${|\xi|}^{H_1-H_r}$ we get, for $\xi$ in $\bR^*$

\begin{equation}\label{qsdqsduiuqsdiusqdfhuqdsddssgyuuydh}
c_1 \h1 \sum_{j \in I_1 } \h1 \alpha_j    \h1 A^{(p)}_{t_j}(\xi) \h1 (1 + \h1 U_{1,j}(\xi)) = - \sum^r_{l = 2 }   c_l \h1 \sum_{j \in I_l } \h1 \alpha_j  \h1 {|\xi|}^{H_1-H_l}  \h1 A^{(p)}_{t_j}(\xi) \h1 (1 + \h1 U_{l,j}(\xi)).
\end{equation}

Since $H_1 > H_l$ for $l$ in $\{2;\cdots;r\}$, taking the limit when $\xi$ tends to $0$ in \eqref{qsdqsduiuqsdiusqdfhuqdsddssgyuuydh} and using
lemma \ref{dkqskqmlsddqsmklqdskq} yields $ c_1 \h1 i^p \h1 {(p+1)}^{-1} \h1  \sum_{j \in I_1 } \h1 \alpha_j  \h1 t^{p+1}_j = 0$, which is nothing but $(L_{p+1})$. This  ends the proof.
  \qed
\end{praa}

\begin{rem}
Another way to establish that $R_H(.,.)$ defines an inner product on $\widetilde{\cE}(\bR)$ for $H$ in $(0,1)$ is to use \eqref{alban} and \eqref{ezdzedeferf}.
\end{rem}


\bibliographystyle{plain}
\bibliography{bib1}

\begin{thebibliography}{10}

\bibitem{nualart}
E.~Alos, O.~Mazet, and D.~Nualart.
\newblock Stochastic {C}alculus with {R}espect to {G}aussian {P}rocesses.
\newblock {\em Annals of Probability}, 29(2):766--801, 2001.

\bibitem{ACLVL}
A.~Ayache, S.~Cohen, and J.~Lévy-Véhel.
\newblock The covariance structure of multifractional {B}rownian motion, with
  application to long range dependence (extended version).
\newblock {\em ICASSP}, Refereed Conference Contribution, 2000.

\bibitem{BJR}
A.~Benassi, S.~Jaffard, and D.~Roux.
\newblock Elliptic {G}aussian random processes.
\newblock {\em Rev. Mat. Iberoamericana}, 13(1):19--90, 1997.

\bibitem{ben1}
C.~Bender.
\newblock An {I}tô formula for generalized functionals of a fractional
  {B}rownian motion with arbitrary {H}urst parameter.
\newblock {\em Stochastic Processes and their Applications}, 104:81--106, 2003.

\bibitem{ben2}
C.~{Bender}.
\newblock An {S}-transform approach to integration with respect to a fractional
  {B}rownian motion.
\newblock {\em Bernouilli}, 9(6):955--983, 2003.

\bibitem{bosw}
F.~Biagini, A.~Sulem, B.~{\O}ksendal, and N.N. Wallner.
\newblock An introduction to white-noise theory and {M}alliavin calculus for
  fractional {B}rownian motion.
\newblock {\em Proc. Royal Society, special issue on stochastic analysis and
  applications}, pages 347--372, 2004.

\bibitem{bianchi}
S.~Bianchi and A.~Pianese.
\newblock Modelling stock price movements: multifractality or
  multifractionality?
\newblock {\em Quantitative Finance}, 7(3):301--319, 2007.

\bibitem{che}
J.Y. Chemin.
\newblock Analyse harmonique et équation des ondes et de {S}chrödinger.
\newblock Polycopié de cours, 2003.

\bibitem{GeCh}
I.M. Chilov and G.E. Gelfand.
\newblock {\em Les distributions}, volume~1.
\newblock Dunod, 1962.

\bibitem{GeCh2}
I.M. Chilov and G.E. Gelfand.
\newblock {\em Les distributions}, volume~2.
\newblock Dunod, 1962.

\bibitem{COU}
L.~Coutin.
\newblock An {I}ntroduction to ({S}tochastic) {C}alculus with {R}espect to
  {F}ractional {B}rownian {M}otion.
\newblock {\em Séminaire de Probabilités XL}, 1899:3--65, 2007.

\bibitem{DecUs}
L.~Decreusefond and A.~S. {\"U}st{\"u}nel.
\newblock Stochastic analysis of the fractional {B}rownian motion.
\newblock {\em Potential Anal.}, 10(2):177--214, 1999.

\bibitem{Do}
V.~Dobric and F.M. Ojeda.
\newblock Fractional {B}rownian fields, duality, and martingales.
\newblock {\em IMS Lecture Notes, Monograph Series, High Dimensional
  Probability.}, 51:77--95, 2006.

\bibitem{montagnes}
A.~Echelard, O.~Barri{\`e}re, and J.~L{\'e}vy-Véhel.
\newblock {T}errain modelling with multifractional {B}rownian motion and
  self-regulating processes.
\newblock In {L.} {B}olc, {R.} {T}adeusiewicz, {L.}~{J}. {C}hmielewski, and
  {K.}~{W}. {W}ojciechowski, editors, {\em {C}omputer {V}ision and {G}raphics.
  {S}econd {I}nternational {C}onference, {W}arsaw, {P}oland, {S}eptember 20-22,
  {P}roceedings, {P}art {I} {ICCVG}}, volume 6374 of {\em {L}ecture {N}otes in
  {C}omputer {S}cience}, pages 342--351. {S}pringer, 2010.

\bibitem{ell}
R.J. Elliott and J.~Van~der Hoek.
\newblock A general fractional white noise theory and applications to finance.
\newblock {\em Math. Finance}, pages 301--330, 2003.

\bibitem{FV}
P.K. Friz and N.B. Victoir.
\newblock {\em Multidimensional {S}tochastic {P}rocesses as {R}ough {P}aths:
  {T}heory and {A}pplications}.
\newblock Cambridge University Press, 2010.

\bibitem{MNIRV}
M.~Gradinaru, I.~Nourdin, F.~Russo, and P.~Vallois.
\newblock {$m$}-order integrals and generalized {I}t\^o's formula: the case of
  a fractional {B}rownian motion with any {H}urst index.
\newblock {\em Ann. Inst. H. Poincar\'e Probab. Statist.}, 41(4):781--806,
  2005.

\bibitem{Her}
E.~Herbin.
\newblock From n-parameter fractional {B}rownian motions to n-parameter
  multifractional {B}rownian motions.
\newblock {\em Rocky Mountain J. of Math.}, 36-4:1249--1284, 2006.

\bibitem{HLV}
E.~Herbin and J.~Lévy-Véhel.
\newblock Stochastic 2 micro-local analysis.
\newblock {\em Stochastic Processes and their Applications}, 119(7):2277--2311,
  2009.

\bibitem{Hi}
T.~Hida.
\newblock {\em Brownian Motion}.
\newblock Springer-Verlag, 1980.

\bibitem{HKPS}
T.~Hida, H.~Kuo, J.~Potthoff, and L.~Streit.
\newblock {\em White Noise. An Infinite Dimensional Calculus}, volume 253.
\newblock Kluwer academic publishers, 1993.

\bibitem{HP}
E.~Hille and R.S. Phillips.
\newblock {\em Functional Analysis and Semi-Groups}, volume~31.
\newblock American Mathematical Society, 1957.

\bibitem{Hry}
F.~Hirsch, B.~Roynette, and M.~Yor.
\newblock From an {I}t\^o type calculus for {G}aussian processes to integrals
  of log-normal processes increasing in the convex order.
\newblock {\em To appear in Journal of the mathematical society of Japan},
  2011.

\bibitem{HOUZ}
H.~Holden, B.~Oksendal, J.~Ubøe, and T.~Zhang.
\newblock {\em Stochastic Partial Differential Equations, A Modeling, White
  Noise Functional Approach}.
\newblock Springer, second edition, 2010.

\bibitem{janson}
S.~Janson.
\newblock {\em Gaussian Hilbert Spaces}, volume 129.
\newblock Cambridge Tracts in Mathematics, 2008.

\bibitem{Kol}
A.~Kolmogorov.
\newblock {W}ienersche {S}piralen und einige andere interessante {K}urven in
  {H}ilbertsche {R}aum.
\newblock {\em C. R. (Dokl.) Acad. Sci. URSS}, 26:115--118, 1940.

\bibitem{Kuo2}
H.H. Kuo.
\newblock {\em White Noise Distribution Theory.}
\newblock CRC-Press, 1996.

\bibitem{Kuo_1}
H.H. Kuo.
\newblock {\em Introduction to Stochastic Integration}.
\newblock Springer, 2006.

\bibitem{JLJLV2}
J.~Lebovits and J.~Lévy-Véhel.
\newblock Stochastic differential equations driven by a multifractional
  {B}rownian motion using white noise theory.
\newblock 2011.
\newblock In preparation.

\bibitem{LLHF}
M.~Li, S.C. Lim, B.J. Hu, and H.~Feng.
\newblock Towards describing multi-fractality of traffic using local {H}urst
  function.
\newblock In {\em Lecture Notes in Computer Science}, volume 4488, pages
  1012--1020. Springer, 2007.

\bibitem{MVN}
B.~Mandelbrot and J.W. Van~Ness.
\newblock Fractional {B}rownian motions, fractional noises and applications.
\newblock {\em SIAM Rev.}, 10:422--437, 1968.

\bibitem{Nu}
D.~Nualart.
\newblock {\em The Malliavin Calculus and related Topics}.
\newblock Springer, 2006.

\bibitem{PL}
R.~Peltier and J.~Lévy-Véhel.
\newblock Multifractional {B}rownian motion: definition and preliminary
  results, 1995.
\newblock rapport de recherche de l'INRIA, $n°$ $2645$.

\bibitem{JP}
J.~{Picard}.
\newblock Representation formulae for the fractional {B}rownian motion.
\newblock {\em Séminaire de Probabilités XLIII}, 43:3--72, 2011.

\bibitem{TaqSam}
G.~Samorodnitsky and M.S. Taqqu.
\newblock {\em Stable Non-Gaussian Random Processes, Stochastic Models with
  Infinite Variance}.
\newblock Chapmann and Hall/C.R.C, 1994.

\bibitem{StoTaq}
S.A. Stoev and M.S. Taqqu.
\newblock How rich is the class of multifractional {B}rownian motions?
\newblock {\em Stochastic Processes and their Applications}, 116:200--221,
  2006.

\bibitem{Tha}
S.~Thangavelu.
\newblock {\em Lectures of Hermite and Laguerre expansions}.
\newblock Princeton University Press, 1993.

\bibitem{Wid}
D.V. {Widder}.
\newblock Positive temperatures on an infinite rod.
\newblock {\em Trans. Amer. Math. Soc.}, 55:85--95, 1944.

\bibitem{Z}
M.~Z{\"a}hle.
\newblock On the link between fractional and stochastic calculus.
\newblock In {\em Stochastic dynamics ({B}remen, 1997)}, pages 305--325.
  Springer, New York, 1999.

\end{thebibliography}
\nocite{*}

\end{document}